\RequirePackage{rotating}
\documentclass[smallextended, nospthms]{svjour3}       
\smartqed  

\usepackage{amsthm}
\usepackage{amsmath}
\usepackage{amssymb}
\usepackage{graphicx}
\usepackage{color}
\usepackage{url}
\usepackage{algorithm}
\usepackage{mathtools}
\newcommand{\subparagraph}{}
\usepackage{algorithmicx}
\usepackage{algpseudocode}
\usepackage{algcompatible}
\usepackage{bbm}
\usepackage{booktabs} 

\renewcommand{\thesection}{ \arabic{section}}
\renewcommand{\thesubsection}{ \arabic{section}.\arabic{subsection}}
\renewcommand{\thesubsubsection}{\arabic{section}.\arabic{subsection}.\alph{subsubsection}}
\theoremstyle{plain}
\newtheorem{theorem}{Theorem}[section]

\newtheorem{lemma}{Lemma}[section]
\newtheorem{example}{Example}[section]

\theoremstyle{definition}

\newtheorem{assumption}{Assumption}[section]
\newtheorem{remark}{Remark}[section]

\newcommand{\R}{\mathbb{R}}
\newcommand{\Rext}{\R\cup\{+\infty\}}

\newcommand{\set}[1]{\left\{#1\right\}}
\newcommand{\norm}[1]{\left\Vert#1\right\Vert}
\newcommand{\norms}[1]{\Vert#1\Vert}

\newcommand{\Eproof}{\hfill $\square$}
\newcommand{\prox}{\mathrm{prox}}

\newcommand{\kproj}[2]{\mathrm{proj}_{#1}\left(#2\right)}
\newcommand{\relint}[1]{\mathrm{ri}\left(#1\right)}
\newcommand{\argmin}{\mathrm{arg}\!\min}
\newcommand{\dom}[1]{\mathrm{dom}(#1)}
\newcommand{\zero}[1]{{\boldsymbol{0}}}

\newcommand{\Xc}{\mathcal{X}}
\newcommand{\Yc}{\mathcal{Y}}
\newcommand{\Zc}{\mathcal{Z}}

\newcommand{\Nc}{\mathcal{N}}
\newcommand{\Pc}{\mathcal{P}}
\newcommand{\Kc}{\mathcal{K}}

\newcommand{\Xopt}{\mathcal{X}^{\star}}
\newcommand{\Yopt}{\mathcal{Y}^{\star}}

\newcommand{\iprods}[1]{\langle #1\rangle}

\newcommand{\kdist}[2]{\mathrm{dist}_{#1}\left(#2\right)}

\newcommand{\BigO}[1]{\mathcal{O}\left(#1\right)}

\newcommand{\beforesubsec}{\vspace{-3ex}}
\newcommand{\aftersubsec}{\vspace{-2ex}}
\newcommand{\beforesec}{\vspace{-3ex}}
\newcommand{\aftersec}{\vspace{-2ex}}
\newcommand{\beforesubsubsec}{\vspace{-1.0ex}}
\newcommand{\aftersubsubsec}{\vspace{-2.0ex}}
\newcommand{\beforeparagraph}{\vspace{-2.5ex}}

\title{An Adaptive Primal-Dual Framework for \\ Nonsmooth Convex Minimization}
\titlerunning{An Adaptive Primal-Dual Framework for Nonsmooth Convex Minimization}        

\author{Quoc Tran-Dinh \and Ahmet Alacaoglu \and Olivier Fercoq \and Volkan Cevher}
\institute{Quoc Tran-Dinh \at
              Department of Statistics and Operations Research, University of North Carolina at Chapel Hill (UNC-Chapel Hill), 333 Hanes Hall, CB\#3260, UNC Chapel Hill, NC 27599-3260. \\
              \email{\texttt{quoctd@email.unc.edu}}           
           \and
           Ahmet Alacaoglu \and Volkan Cevher \at
              Laboratory for Information and Inference Systems (LIONS), \'Ecole Polytechnique F\'ed\'erale de Lausanne (EPFL), CH1015-Lausanne, Switzerland. \\
              \email{\texttt{\{ahmet.alacaoglu, volkan.cevher\}@epfl.ch}}
          \and
           Olivier Fercoq \at
              LTCI, T\'el\'ecom ParisTech,  Universit\'e Paris-Saclay, 75634-Paris, France. \\
              \email{\texttt{olivier.fercoq@telecom-paristech.fr}}
}

\date{Received: date / Accepted: date}

\begin{document}
\maketitle

\begin{abstract}
\small
We propose a new self-adaptive, double-loop smoothing algorithm to solve composite, nonsmooth, and constrained convex optimization problems. 
Our algorithm is based on Nesterov's smoothing technique via general Bregman distance functions.
It self-adaptively selects the number of iterations in the inner loop to achieve a desired complexity bound without requiring the accuracy a priori as in variants of Augmented Lagrangian methods (ALM).
We prove $\BigO{\frac{1}{k}}$-convergence rate on the last iterate of the outer sequence for both unconstrained and constrained settings in contrast to ergodic rates which are common in ALM as well as alternating direction method-of-multipliers  literature.
Compared to existing inexact ALM or quadratic penalty methods, our analysis does not rely on the worst-case bounds of the subproblem solved by the inner loop.
Therefore, our algorithm can be viewed as a restarting technique applied to the ASGARD method in \cite{TranDinh2015b} but with rigorous theoretical guarantees or as an inexact ALM with explicit inner loop termination rules and adaptive parameters. Our algorithm only requires to initialize the parameters once, and automatically update them during the iteration process without tuning.
We illustrate the superiority of our methods via several examples as compared to the state-of-the-art.  

\vspace{1cm}
\end{abstract}
    
\beforesec
\section{Introduction}\label{sec:intro}
\aftersec
We study the following  nonsmooth composite convex minimization template:
\begin{equation}\label{eq:ns_cvx_prob}
P^{\star} := \min_{x \in \R^p}\Big\{ P(x) := f(x) + g(Ax) \Big\},
\end{equation}
where both $f : \R^p \to \Rext$ and $g : \R^n\to\Rext$ are proper, closed, and nonsmooth convex functions, and $A: \mathbb{R}^p \rightarrow \mathbb{R}^n$ is a linear operator. 

Under only convexity and zero duality gap assumptions, the state-of-the-art methods for solving \eqref{eq:ns_cvx_prob} include primal-dual first-order methods (PDFOM) \cite{Bauschke2011,Chambolle2011}, and augmented Lagrangian-based algorithms \cite{Bertsekas1996d,Boyd2011}. While PDFOM directly tackles problem \eqref{eq:ns_cvx_prob}, the augmented Lagrangian-based framework (ALM) and its variants solve \eqref{eq:ns_cvx_prob} via a constrained reformulation as follows:  
\begin{equation}\label{eq:ns_cvx_prob_splitted}
P^{\star} := \min_{x \in \R^p, z\in\mathbb{R}^n}\Big\{ P(x, z) := f(x) + g(z)~~\mathrm{s.t.}~~ Ax - z = 0 \Big\},
\end{equation}
Alternating direction method of multipliers (ADMM) is  another (and perhaps the most) successful method to solve \eqref{eq:ns_cvx_prob_splitted}.
ADMM can be viewed as an approximation to ALM by alternating between $x$ and $z$ to break the computational bottleneck in the primal subproblem. Inexact and linearized variants  enhance the scalability of ALM and ADMM for the same problem template \cite{xu2017accelerated,ouyang2015accelerated,xu2017iteration}. 

While ADMM and PDFOM and their variants work really well in practice, their best-known convergence rate is $\BigO{\frac{1}{k}}$ under only convexity and zero duality gap assumptions, where $k$ is the iteration counter.
Moreover, such a rate is achieved via an ergodic sense (i.e., using an averaging sequence or a weighted averaging sequence) \cite{Chambolle2011,chambolle2016ergodic,Davis2014b,Davis2014,Monteiro2010,Monteiro2012b,Shefi2014}.

In stark contrast, empirical evidence shows that averaging sequences in PDFOM and ALM  exhibit the theoretical worst case rate $\BigO{\frac{1}{k}}$ in practice compared to the last iterate of the algorithm (see Subsection~\ref{subsec:example1} for a concrete example), which is superior and often locally linear in many examples.\footnote{There exist examples showing arbitrarily slow convergence rate of ADMM, see, e.g., \cite{Davis2014b}.} 
However, for these methods, last iterate generally has convergence guarantees but has no rate guarantees.

Recently, \cite{TranDinh2015b} proposed an \textit{accelerated, smoothed gap reduction} (ASGARD) framework to solve nonsmooth convex optimization problems.
ASGARD combines acceleration, smoothing, and homotopy techniques to handle both unconstrained and constrained nonsmooth problems, including \eqref{eq:ns_cvx_prob}.

One notable feature of ASGARD is a non-ergodic optimal $\BigO{\frac{1}{k}}$ rate on the objective residual, and feasibility violation in the constrained settings.
Moreover, this method only requires one proximal operator of $f$, one matrix-vector multiplication, and one adjoint operator per iteration. 
When $f$ is separable, the algorithm can be naturally parallelized. 
However, as also noted in~\cite{TranDinh2015b}, ASGARD needs restarting to be competitive with state-of-the-art methods such as ADMM and PDFOM in practice.
This is not surprising since empirical evidence \cite{fercoq2016restarting,Giselsson2014,Odonoghue2012,Su2014} has shown that restarting significantly improves the actual convergence rate in practice. While there exists theory to support the restarting strategies in accelerated gradient-type methods, supporting theory of these strategies are not yet investigated in primal-dual methods.

In this paper, we introduce an analysis framework for restarting ASGARD and prove the same worst-case $\BigO{\frac{1}{k}}$ rate in a non-ergodic sense.
To our knowledge, this is the first time restarting is studied for primal-dual methods with a rigorous convergence rate guarantee.
While doing so, we identify that restarting ASGARD corresponds to an inexact ALM algorithm in the constrained case.
In contrast to existing works on this front, our method has explicit inner-loop termination rules and does not need to set a horizon (i.e., the maximum number of inner iterations or a predefined inner loop accuracy) for the algorithm.

As a result, we present a method which has the guarantees on the last iterate compared to ALM/ADMM methods and extend the guarantees of ASGARD to the restarting case which significantly improves the practical performance.
In addition, we allow general Bregman distances to be used for smoothing and proximal operators in contrast to the original ASGARD scheme.
A more thorough discussion and comparison between our method and existing state-of-the-arts is deferred to Section \ref{sec:prev-work} for the sake of presentation.

\beforeparagraph
\paragraph{\bf Composite vs. constrained settings:}
An interesting instance of \eqref{eq:ns_cvx_prob} is the following constrained convex setting:
\begin{equation}\label{eq:constr_cvx}
f^{\star} := \min_{x\in\R^p}\Big\{ f(x) ~\mid~ Ax - b \in \Kc \Big\},
\end{equation}
where $\Kc$ is a nonempty, closed, and convex set in $\R^n$, and $b\in\R^n$.
If we define $g(u) := \delta_{b + \Kc}(u)$, the indicator of $b+\Kc$, then \eqref{eq:constr_cvx} can be converted into \eqref{eq:ns_cvx_prob}.

In the general setting \eqref{eq:ns_cvx_prob}, under different choices of $f$ and $g$,  \eqref{eq:ns_cvx_prob} covers a wide range of applications from different fields including compressive sensing, image/signal processing, machine learning, statistics, and optimal control.
Classical and well-known examples such as LASSO, square-root LASSO, support vector machines, image denoising and deblurring, and matrix completion can be cast  into \eqref{eq:ns_cvx_prob},
see, e.g., \cite{Boyd2011,Combettes2011,Parikh2013,wright2017optimization} for some concrete examples.

For the setting \eqref{eq:constr_cvx}, we do not impose any restriction on $\Kc$.
Hence, it covers a large class of constrained problems including equality and inequality constraints.
When $\Kc$ is a given cone (e.g., $\R^n_{+}$, second-order cone, or symmetric positive semidefinite cone), problem \eqref{eq:constr_cvx} covers also problems with cone constraints such as linear programming, second-order cone, and semidefinite programming.
Although the theory for \eqref{eq:ns_cvx_prob} as well as for \eqref{eq:constr_cvx} are well developed, various numerical methods for solving these problems rely on different structure assumptions and do not have a unified analysis: \textit{cf.}, Section \ref{sec:prev-work}.

\beforeparagraph
\paragraph{\bf Contributions:}
Our contributions can be summarized as follows.
\begin{itemize}
\vspace{-1ex}
\item[$\mathrm{(a)}$] 
We propose a new self-adaptive, double-loop smoothing algorithm to solve nonsmooth convex optimization problems of the form \eqref{eq:ns_cvx_prob}.
Our algorithm is based on Nesterov's smoothing technique via general Bregman distance functions.
It self-adaptively selects the number of iterations in the inner loop to achieve a desired complexity bound without requiring the accuracy a priori as in variants of ALM.
Compared to ASGARD \cite{TranDinh2015b}, {it incorporates restarts,} updates the dual center, and can work with general Bregman distances instead of only Lipschitz gradient distances.

\item[$\mathrm{(b)}$]
We prove $\BigO{\frac{1}{k}}$-convergence rate on the last iterate of the outer sequence for both unconstrained and constrained settings in contrast to ergodic rates which are common in ALM/ADMM literature.
This rate is known to be optimal \cite{Nemirovskii1983,Nesterov2005c,woodworth2016tight} under just convexity and strong duality assumptions. 
Compared to existing inexact ALM or quadratic penalty methods such as \cite{necoara2015complexity,xu2016homotopy}, our analysis does not rely on the worst-case bounds of the subproblem solved by the inner loop.
Therefore, our algorithm can be viewed as a restarting technique applied to ASGARD but with rigorous theoretical guarantees {or as an inexact ALM with explicit inner loop termination rules and adaptive parameters.}

\item[$\mathrm{(c)}$]
As an upshot, we customize our algorithm  to solve general constrained problems of the form \eqref{eq:constr_cvx}.
We prove the same $\BigO{\frac{1}{k}}$-convergence rate guarantee on both the objective residual $\vert f(x^k) - f(x^{\star})\vert$ and the feasibility $\kdist{\Kc}{Ax^k - b}$.
This rate is given on the last iterate of the outer sequence.
\end{itemize}

Our algorithm is a primal-dual method,  which can solve composite convex problem with linear operators as in Chambolle-Pock's method \cite{Chambolle2011}.
It only requires one proximal operator of $f$ and $g^{*}$, one matrix-vector multiplication and one adjoint for each iteration.
It is parallelizable when $f$ is separable, i.e., $f(x) = \sum_{i=1}^Nf_i(x_{[i]})$. 
Under this structure, our method has more advantages than ADMM and Chambolle-Pock's method.
In the algorithm, we provide explicit rules to update all  algorithmic parameters. 
We also note that these updates can be modified to trade-off  between the primal or the dual progress.

\vspace{-2ex}
\paragraph{\bf Paper organization:}
The rest of this paper is organized as follows.
Section \ref{sec:asgard} recalls some mathematical background and the ASGARD algorithm in \cite{TranDinh2015b}.
Section \ref{sec:restart_asgard} presents our main result with algorithm and its convergence guarantee.
We study  both unconstrained and constrained cases. 
Section \ref{sec:extension} shows an extension of our method to three composite objective functions with linearization on potentially smooth terms.
We provide overall complexity bounds in Section~\ref{sec:complexity} for linear programming and discuss their superiority compared to existing results.
In Section \ref{sec:num_exp}, we provide seven numerical examples to test our algorithm against state-of-the arts.
Section \ref{sec:prev-work} compares our method and existing algorithms in the literature.

\beforesec
\section{Mathematical tools and ASGARD}\label{sec:asgard}
\aftersec
We review some key ingredients for the design of our primal-dual methods.
We also recall the ASGARD algorithm in \cite{TranDinh2015b} and discuss its possible variants. 
\beforeparagraph
\paragraph{\bf Notation:} 
We denote the norm in primal space $\Xc$ as $\Vert \cdot \Vert _\Xc$ and the norm in dual space $\Yc$ as $\Vert \cdot \Vert _\Yc$. 
Their dual norms are denoted as $\Vert \cdot \Vert _{\Xc, \ast}$ and $\Vert \cdot \Vert _{\Yc, \ast}$, respectively.
Given a proper, closed, and convex function $f$, we use $\mathrm{dom}(f)$ to denote its domain and $\partial f(x)$ to denote its subdifferential at $x$. 
When the function is differentiable, we denote its gradient at $x$ as $\nabla f(x)$.
For a given nonempty, closed, and convex set $\Kc$, we denote its indicator function as $\delta_\Kc (x)=0$, if $x\in\Kc$, $\delta_\Kc(x)=+\infty$, otherwise; and its support function as $s_\Kc(y)=\sup_{x\in\Kc} \iprods{x, y}$.
 We define the normal cone of $\Kc$ as $\Nc_{\Kc}(x) := \set{ w\in\R^n \mid \iprods{w, y-x} \geq 0,~y\in\Kc}$ if $x\in\Kc$; $\Nc_{\Kc}(x) := \emptyset$, otherwise.
We also define $\Kc^{o} := \set{w \in\R^n \mid \iprods{w, x} \leq 1, ~x\in\Kc}$ as the polar set of $\Kc$.
If $\Kc$ is a convex cone, then $\Kc^o = -\Kc^{\ast}$, where $\Kc^{\ast} := \set{ w\in\R^n \mid \iprods{w, x} \geq 0, ~x\in\Kc}$ the dual cone of $\Kc$.
The Fenchel conjugate of a function $f$ is defined as $f^\ast(y) := \sup_{x}\set{ \iprods{x, y} - f(x) }$.
We say that  $f: \Xc \rightarrow \mathbb{R}$ has Lipschitz gradient if it satisfies $\norm{\nabla f(x) - \nabla f (y)}_{\Xc, \ast} \leq L_f \norm{ x - y}_\Xc$, for any $x, y\in\Xc$. 
This is equivalent to $f(x) \leq f(y) + \iprods{\nabla f(y), x-y} + \tfrac{L_f}{2} \norm{ x-y}_\Xc^2$, for all $x, y\in\Xc$.
Given a positive real number $a$, $\lfloor a\rfloor$ denotes the largest integer that is less than or equal to $a$.

Given a proper, closed, and convex function $f : \R^p\to\Rext$, $\prox_{f}(x) := \argmin_{u}\set{ f(u) + (1/2)\norms{u-x}_{\Xc}^2}$ is called the proximal operator of $f$.
We say that $f$ is ``proximally tractable'' if $\prox_f$ can be computed efficiently, e.g., in a closed form, or by a polynomial algorithm. 
By Moreau's identity, we have $\prox_{\gamma f}(x) + \gamma\prox_{f^{\ast}/\gamma}(\gamma^{-1}x) = x$ for any $x\in\dom{f}$.

\beforesubsec
\subsection{\bf Primal-dual formulation}
\aftersubsec
\noindent\textbf{\textit{Dual problem and min-max formulation:}}
Associated with the primal problem \eqref{eq:ns_cvx_prob}, we also consider the corresponding dual problem:
\begin{equation}\label{eq:ns_cvx_dual}
D^{\star} := \min_{y \in \R^n}\Big\{ D(y) := f^{\ast}(-A^{\top}y) + g^{\ast}(y) \Big\},
\end{equation}
where $f^{\ast}$ and $g^{\ast}$ are the Fenchel conjugates of $f$ and $g$, respectively.
Clearly, we can write the primal and dual pair \eqref{eq:ns_cvx_prob}-\eqref{eq:ns_cvx_dual} in the following min-max saddle point problem:
\begin{equation}\label{eq:ns_cvx_minmax}
\begin{array}{ll}
P^{\star} &= \displaystyle\min_{x\in\R^p}\displaystyle\max_{y\in\R^n}\Big\{ f(x) + \iprods{Ax, y} - g^{\ast}(y) \Big\} \vspace{1ex}\\
& = \displaystyle\max_{y\in\R^n}\displaystyle\min_{x\in\R^p}\Big\{ - g^{\ast}(y) - (\iprods{x, -A^{\top}y} - f(x)) \Big\} = -D^{\star}.
\end{array}
\end{equation}
Under mild and standard assumptions, this min-max problem is solvable and achieves  zero duality gap, i.e., $P^{\star} + D^{\star} = 0$.
In particular, the dual problem of \eqref{eq:constr_cvx} can be written as follows:
\begin{equation}\label{eq:constr_cvx_dual}
D^{\star} := \min_{y\in\R^n}\set{ D(y) := f^{\ast}(-A^Ty) + \iprods{b, y} + s_{\Kc}(y) },
\end{equation}
where $s_{\Kc}(y) = \sup_{x\in\Kc} \langle y, x \rangle$ is the support function of $\Kc$.
Compared to \eqref{eq:ns_cvx_dual}, we have $g^{\ast}(y) = \iprods{b, y} + s_{\Kc}(y) = s_{b + \Kc}(y)$.
Let $\Xopt$ and $\Yopt$ be the solution sets of the primal problem \eqref{eq:ns_cvx_prob} (or \eqref{eq:constr_cvx}) and dual problem \eqref{eq:ns_cvx_dual} (or \eqref{eq:constr_cvx_dual}), respectively.

\vspace{1ex}
\noindent\textbf{Fundamental assumptions:}
Throughout this paper, we will develop methods for solving~\eqref{eq:ns_cvx_prob} and~\eqref{eq:constr_cvx}. 
Note that we will use different assumptions for these two cases, which are given below for ~\eqref{eq:ns_cvx_prob} and~\eqref{eq:constr_cvx}, respectively.

\begin{assumption}\label{as:A1}
The solution set $\Xc^{\star}$ of \eqref{eq:ns_cvx_prob} is nonempty.
Both $f$ and $g$ are proper, closed, and convex.
Moreover, $\dom{g^{\ast}}$ is bounded, or equivalently, $g$ is Lipschitz continuous.
Note that this implies the Slater condition $\boldsymbol{0} \in \relint{\dom{g} - A(\dom{f})}$, where $\relint{\Xc}$ is the relative interior of $\Xc$.
\end{assumption}
\begin{assumption}\label{as:A1b}
The solution set $\Xc^{\star}$ of \eqref{eq:constr_cvx} is nonempty.
The function $f$ is proper, closed, and convex.
The constraint set $\Kc$ is nonempty, closed, and convex, and $\boldsymbol{0}^n\in\Kc$.
The Slater condition $\relint{\dom{f}} \cap \set{x\in\R^p \mid Ax - b\in \relint{\Kc}} \neq\emptyset$ holds, where $\relint{\Xc}$ is the relative interior of $\Xc$.
\end{assumption}

\noindent 
Except for the boundedness of $\dom{g^{\ast}}$, Assumptions~\ref{as:A1} and \ref{as:A1b} are very standard in convex optimization. 
It guarantees the strong duality of \eqref{eq:ns_cvx_prob} (respectively, \eqref{eq:constr_cvx}) and \eqref{eq:ns_cvx_dual} to hold.
The boundedness of $\dom{g^{\ast}}$ is guaranteed if and only if $g$ is Lipschitz continuous as we mentioned.
We emphasize that we need boundedness of $\dom{g^\ast}$ only for~\eqref{eq:ns_cvx_prob} and we do not require it for~\eqref{eq:constr_cvx}

We note that the assumption $\boldsymbol{0}^n\in\Kc$ is not restrictive, since if $\boldsymbol{0}^n\notin\Kc$, we can fix any point $\boldsymbol{e}\in\Kc$, and consider the set $\tilde{\Kc} = \Kc - \boldsymbol{e}$, then $\boldsymbol{0}^n\in\tilde{\Kc}$, and $Ax - b\in\Kc$ becomes $Ax - b + \boldsymbol{e}\in\tilde{\Kc}$.
Note that, in the sequel, we will refer to the setting of~\eqref{eq:ns_cvx_prob} with Assumption~\ref{as:A1} as the bounded dual domain case, and to the setting of~\eqref{eq:constr_cvx} with Assumptions~\ref{as:A1b} as the constrained case.

\vspace{1ex}
\noindent\textbf{Optimality conditions:}
Associated with the primal and dual problems \eqref{eq:ns_cvx_prob}-\eqref{eq:ns_cvx_dual}, we have the following optimality conditions:
\begin{equation}\label{eq:kkt_cond}
0 \in \partial{f}(x^{\star}) + A^{\top}\partial{g}(Ax^{\star})~~~\text{and}~~~0 \in -A\partial{f^{\ast}}(-A^{\top}y^{\star}) + \partial{g^{\ast}}(y^{\star}).
\end{equation}
From~\eqref{eq:ns_cvx_minmax}, it is straightforward to see the relation $y^{\star} \in \partial{g}(Ax^{\star}) \Leftrightarrow Ax^{\star} \in \partial{g^{\ast}}(y^{\star})$ and $x^{\star} \in \partial{f^{\ast}}(-A^{\top}y^{\star}) \Leftrightarrow -A^{\top}y^{\star} \in\partial{f}(x^{\star})$, we can write this optimality condition into the following KKT condition:
\begin{equation*}
0 \in \partial{f}(x^{\star}) + A^{\top}y^{\star}~~~\text{and}~~~0 \in -Ax^{\star} + \partial{g^{\ast}}(y^{\star}).
\end{equation*}
For the constrained problem \eqref{eq:constr_cvx} these conditions are written as 
\begin{equation*}
0 \in \partial{f}(x^{\star}) + A^{\top}y^{\star},~~~Ax^\star - b \in \Kc, ~~\text{and}~~y^{\star} \in \Nc_{\Kc}\left(Ax^{\star} -  b\right),
\end{equation*}
where $\Nc_{\Kc}(\cdot)$ is the normal cone of $\Kc$ defined above.
If $\Kc$ is a closed, pointed, and convex cone, then $\Nc_{\Kc} \equiv -\Kc^{\ast}$ the dual cone of $\Kc$.
In this case, $y^{\star} \in -\Kc^{\ast}$.

\beforesubsec
\subsection{Bregman Distances and Generalized Proximal Operators}\label{sec:bregman}
\aftersubsec
In the sequel, we will use Bregman distances for smoothing and computing proximal operators. 
Therefore, we give basic properties on Bregman distances.

Let $p_\Zc$ be $\mu_p$-strongly convex, continuous, and differentiable on $\Zc$ with the strong convexity $\mu_{p} = 1$, where $\Zc = \text{dom}(p_\Zc)$. 
We call $p_{\Zc}$ a proximity function (or prox-function).
We define the Bregman distance induced by  $p_\Zc$ as
\begin{equation*}
b_\Zc(x, y) := p_\Zc(x) - p_\Zc(y) - \langle \nabla p_\Zc(y), x - y \rangle, ~~ \forall x, y \in \Zc.
\end{equation*}
We assume that $b_\Zc$ is $1$-strongly convex with respect to the norm $\Vert \cdot \Vert _\Zc$. Then
\begin{equation}
b_\Zc(x, y) \geq \tfrac{1}{2} \Vert x - y \Vert_\Zc^2, ~~ \forall x, y\in\Zc.
\end{equation}
A special case of prox-functions is $p_\Zc(x) = \frac{1}{2} \Vert x \Vert_2 ^2$, which corresponds to the well-known Euclidean distance $b_\Zc(x, y) = \frac{1}{2} \Vert x - y \Vert _2 ^2$.
Another example is the entropy function $p_\Zc(x) := \sum_i x_i \ln(x_i)$, which corresponds to the so-called KL divergence $b_\Zc(x, y) := \sum_i x_i \ln\left(\frac{x_i}{y_i}\right) - x_i + y_i$.
When a Bregman distance $b_{\Zc}$ has Lipschitz continuous gradient, we denote its Lipschitz constant by $L_{b_\Zc}$.

We also define the strong convexity of a function $f$ with respect to a prox-function $p_\Zc$ which induces the Bregman distance $b_\Zc$ as follows:
\begin{equation}
f(x) \geq f(y) + \langle \nabla f(y), x - y \rangle + b_\Zc(x,y), ~~~ \forall x, y \in\Zc.
\end{equation}
We refer to \cite{chen1993convergence,eckstein1993nonlinear,kiwiel1997proximal}  for several concrete examples of Bregman divergences.

\beforesubsec
\subsection{Nesterov's smoothing technique}
\aftersubsec
We focus on Nesterov's smoothing technique with general Bregman distances \cite{Beck2012a,Nesterov2005c}.
Since $g$ is nonsmooth, assuming that it admits a max-form as 
\begin{equation*}
g(u) = \max_{y \in\Yc}\set{\iprods{u, y} - g^{\ast}(y)},~~~\text{where}~\Yc = \dom{g^{\ast}},
\end{equation*}
we smooth it by
\begin{equation}\label{eq:smooth_g}
\begin{array}{ll}
g_{\beta}(u; \dot{y}) &:= \displaystyle\max_{y\in\Yc}\set{ \iprods{u, y} - g^{\ast}(y) - \beta b_\Yc (y, \dot{y})},
\end{array}
\end{equation}
where $\dot{y}\in\R^n$ is a given center point, and $\beta > 0$ is a smoothness parameter.
The function $g_{\beta}(\cdot;\dot{y})$ is convex and smooth, its gradient is given by 
\begin{equation}\label{eq:grad_smooth_g_bregman}
\nabla{g_{\beta}}(u;\dot{y}) = y^{\ast}_{\beta}(u;\dot{y}) = \argmin_{y\in\Yc}\Big\{ g^{\ast}(y) - \iprods{u, y} + \beta b_\Yc(y, \dot{y}) \Big\}.
\end{equation}
Clearly, $\nabla{g_{\beta}}(\cdot;\dot{y})$ is Lipschitz continuous with the Lipschitz constant $L_{g_{\beta}} = \frac{1}{\beta}$.
Moreover, we have 
\begin{equation}\label{eq:diameter_bound}
g_{\beta}(u; \dot{y}) \leq g(u) \leq g_{\beta}(u; \dot{y}) + \beta D_{\Yc},
\end{equation}
where $D_{\Yc} := \sup\set{b_\Yc(y, \dot{y}) \mid y\in\dom{g^{\ast}}}$ is the prox-diameter of $g^{\ast}$.
Here, $D_{\Yc}$ is finite if and only if $g$ is Lipschitz continuous with the Lipschitz constant $L_g := \sqrt{2D_{\Yc}}$, i.e., $\vert g(u) - g(v)\vert \leq \sqrt{2D_{\Yc}}\Vert u - v\Vert$ for all $u, v\in\dom{g}$ due to \cite[Proposition 4.4.6]{borwein2010convex}.

If we choose $b_\Yc(y, \dot{y}) = \frac{1}{2}\Vert y - \dot{y} \Vert _2 ^2$, then we can write $y_\beta^\ast (u; \dot{y})$ as:
\begin{equation}\label{eq:grad_smooth_g_l2}
{\!\!\!\!}\nabla{g_{\beta}}(u;\dot{y}) = \argmin_{y\in\R^n}\set{ g^{\ast}(y) - \iprods{u, y} + \tfrac{\beta}{2}\Vert y - \dot{y}\Vert^2}  = \prox_{g^{\ast}/\beta}\left(\dot{y} + \tfrac{1}{\beta}u\right).{\!\!\!\!}
\end{equation}
Smoothing techniques are widely used in the literature, including \cite{Beck2012a,boct2012variable,Bot2013,Devolder2012,Necoara2008}.
The idea of smoothing is to approximate the original problem \eqref{eq:ns_cvx_prob} by a (partially) smoothed problem. 
For example, in our setting, we smooth $g$ and consider the following smoothed problem:
\begin{equation}\label{eq:smoothed_cvx_prob}
P_{\beta}^{\star} := \min_{x\in \mathbb{R}^p}\Big\{ P_{\beta}(x;\dot{y}) := f(x)  + g_{\beta}(Ax;\dot{y}) \Big\}.
\end{equation}
We define the following generalized proximal operator with Bregman distance $d_\Xc$ induced by a prox-function $q_\Xc$:
\begin{equation}\label{eq:general_prox}
\Pc^{d_\Xc}_{\theta f}(u, y) := \argmin_{v \in \Xc} \Big\{ f(v) + \iprods{y, v - u} + \tfrac{1}{\theta} d_\Xc(v, u) \Big\}.
\end{equation}
Note that the setup described in this and previous subsections will allow us to use different Bregman distances for smoothing as in~\eqref{eq:smooth_g} and computing the proximal operator as in~\eqref{eq:general_prox}, depending on the geometry of the problem.
Given that the Bregman distance $d_\Xc$ is defined in $\Xc$ and $b_\Yc$ is defined in $\Yc$, we define the following operator norm of $A$:
\begin{equation}
\norm{A} :=  \max_{x\in\R^p} \set{ \frac{\norm{Ax}_{\Yc, \ast}}{\norm{x}_\Xc} }.
\end{equation}
Different from \cite{Beck2012a,Bot2013,Devolder2012,Necoara2008,Nesterov2005c}, our strategy allows one to update the smoothness parameter $\beta$ gradually at each iteration.
Similar work can be found in \cite{boct2012variable,Nesterov2005d}, which are also essentially different from ours as discussed in Section~\ref{sec:prev-work}.

\beforesubsec
\subsection{ASGARD: A primal-dual gap reduction framework}
\aftersubsec
In \cite{TranDinh2015b}, the authors  proposed two primal-dual algorithms to solve \eqref{eq:ns_cvx_prob}.
The first algorithm, ASGARD (Accelerated Smoothed Gap Reduction),  can be viewed as a variant of FISTA \cite{Beck2009} applied to the smoothed problem of \eqref{eq:constr_cvx}.
The second one, ADSGARD (Accelerated Dual Smoothed Gap Reduction) is a Nesterov's accelerated variant  \cite{Nesterov2004} applied to the smoothed problem of  the dual \eqref{eq:ns_cvx_dual}.

\beforeparagraph
\paragraph{\bf ASGARD:}
Let us recall the first algorithm, ASGARD, from \cite{TranDinh2015b} as in Algorithm \ref{alg:A1} for our further reference.

\begin{algorithm}[hpt!]
\caption{($\mathrm{ASGARD}$ - \textbf{A}ccelerated \textbf{S}moothed \textbf{GA}p \textbf{R}e\textbf{D}uction)}\label{alg:A1}
\begin{normalsize}
\begin{algorithmic}[1]
\State {\hskip0ex}\textbf{Initialization:} 
	\vspace{0.75ex}
	\State Choose $x^0 \in \mathbb{R}^p$, $\dot{y}\in\R^n$, and $\beta_0 > 0$ (e.g., $\beta_0 := \Vert A\Vert$). 
		\vspace{0.75ex}	
	          Set $\tau_0 \leftarrow 1$ and $\bar{x}^0=\hat{x}^0 \leftarrow x^0$. 
	          \State  Choose $d_\Xc(\cdot, \dot{x}) = \tfrac{1}{2} \| \cdot - \dot{x} \|^2$, and a Bregman distance $b_\Yc$ as in Section~\ref{sec:bregman}.{\!\!\!}
	\State \textbf{For $k := 0$ to $k_{\max}$ perform}\vspace{1ex}
		\vspace{0.75ex}
		\State{\hskip4ex}\label{eq:ASGARD}Update 
		\vspace{0.75ex}
		$\left\{\begin{array}{ll}
		\tilde{x}^k & \leftarrow (1-\tau_k)\bar{x}^k + \tau_k \hat{x}^k \vspace{0.75ex}\\
		\tilde{y}^{k+1} & \leftarrow  \displaystyle\argmin_{y\in\Yc}\set{ g^{\ast}(y) - \iprods{A\tilde{x}^k, y} + \beta_k b_\Yc(y, \dot{y})}   \vspace{0.75ex}\\
		\bar{x}^{k+1} &\leftarrow \Pc^{d_\Xc} _{(\beta_k/\norm{A}^{2})f}\left(\tilde{x}^k, A^\top \tilde{y}^{k+1}\right) \vspace{0.75ex}\\
		\hat{x}^{k+1} &\leftarrow \hat{x}^{k} + \frac{1}{\tau_k}(\bar{x}^{k+1} - \tilde{x}^k).
		\end{array}\right.$
		\vspace{0.75ex}
	        \State{\hskip4ex}Compute $\tau_{k+1} \in (0, 1)~\text{by solving}~\frac{\tau^3}{L_{b_{\Yc}}} + \tau^2 + \tau_k^2\tau - \tau_k^2 = 0$ in $\tau$.
	        \vspace{0.75ex}
	        \State{\hskip4ex}Update $\beta_{k+1} \leftarrow \frac{\beta_k}{1 + L_{b_{\Yc}}^{-1}\tau_{k+1}}$.
	        \vspace{0.75ex}
	\State\textbf{End~for}
\end{algorithmic}
\end{normalsize}
\end{algorithm}
The main step of ASGARD, Algorithm~\ref{alg:A1}, is Step~\ref{eq:ASGARD}, which requires one subproblem in $\tilde{y}^{k+1}$, one $\prox_f$ of $f$, one matrix-vector multiplication $Ax$ and its adjoint $A^{\top}y$ at each iteration.
If $b_{\Yc}(\cdot, \dot{y}) := \tfrac{1}{2}\norms{\cdot - \dot{y}}_2 ^2$, then the computation of $\tilde{y}^{k+1}$ reduces to proximal operator $\prox_{g^{\ast}}$ of $g^{\ast}$.
Therefore, the per-iteration complexity of Algorithm \ref{alg:A1} is the same as in several primal-dual first-order algorithms \cite{Chambolle2011,Esser2010a}.
In contrast to existing works, the convergence guarantee is on the last primal iterate $\set{\bar{x}^k}$, instead of its weighted average $\set{\tilde{x}^k}$.

\beforeparagraph
\paragraph{\textbf{Variants:}}
Although ASGARD relies on FISTA \cite{Beck2009}, one can replace Step~\ref{eq:ASGARD} by any other accelerated proximal-gradient scheme such as Tseng's variant (APG) in \cite{tseng2008accelerated}.
One can also use two proximal operator schemes from \cite{lan2011primal,Nesterov2005c} to substitute Step~\ref{eq:ASGARD}.
To avoid the overload of this paper, we skip the analysis of Algorithm~\ref{alg:A1} and its variants, which can be found in~\cite{TranDinh2015b}.

\beforeparagraph
\paragraph{\textbf{Convergence:}}
As proved in \cite{TranDinh2015b}, if Algorithm~\ref{alg:A1} is applied to solve  \eqref{eq:ns_cvx_prob} with $b_{\Yc}(\cdot;\dot{y}) := \frac{1}{2}\norms{\cdot - \dot{y}}_{\Yc}^2$, then, under  Assumption~\ref{as:A1}, one has
\begin{equation*}
P(\bar{x}^k) - P^{\star} \leq \BigO{\frac{\norms{A} \norms{x^0 - x^{\star}}_{\Xc} D_{\Yc}}{k}}.
\end{equation*}
If we apply Algorithm~\ref{alg:A1} to solve the constrained problem \eqref{eq:constr_cvx}, then, under Assumption~\ref{as:A1b}, we obtain the following guarantee
\begin{equation*}
\left\{\begin{array}{lll}
&\vert f(\bar{x}^k) - f^{\star}\vert &\leq \BigO{\dfrac{\norms{A}\norms{x^0 - x^{\star}}_{\Xc}\norms{  y^{\star}}_{\Yc}}{k}} \vspace{1ex}\\
&\kdist{\Kc}{A\bar{x}^k - b} &\leq \BigO{\dfrac{\norms{A}\norms{x^0 - x^{\star}}_{\Xc} \norms{y^{\star}}_{\Yc}}{k}}.
\end{array}\right.
\end{equation*}
Hence, the convergence rate of Algorithm~\ref{alg:A1} under Assumption~\ref{as:A1} or Assumption~\ref{as:A1b} is $\BigO{\frac{1}{k}}$ and is in a non-ergodic sense.

\beforesec
\section{Main results: Self-Adaptive Double-Loop ASGARD}\label{sec:restart_asgard}
\aftersec
In this section, we develop a self-adaptive double-loop accelerated smoothed primal-dual gap reduction algorithm to solve \eqref{eq:ns_cvx_prob} and \eqref{eq:constr_cvx}.
We first present the complete algorithm.
Next, we provide its convergence analysis.
Then, we specify our algorithm to handle the constrained setting \eqref{eq:constr_cvx}.
Finally, we extend our method to handle \eqref{eq:ns_cvx_prob} with the sum of three objective functions where the third function has Lipschitz gradient.

\beforesubsec
\subsection{\bf The algorithm and its convergence guarantee}
\aftersubsec
\textbf{Main idea:}
The proposed algorithm consists of two loops:
\begin{itemize}
\vspace{-1ex}
\item The inner loop performs an accelerated proximal gradient scheme (APG)~\cite{tseng2008accelerated} to solve the smoothed problem \eqref{eq:smoothed_cvx_prob} for a fixed $\beta$, which is a different strategy from \cite{TranDinh2015b}, where $\beta$ is updated at each iteration.
We note that in the constrained case, the smoothed problem~\eqref{eq:smoothed_cvx_prob} is the augmented Lagrangian.
\item The outer loop can be considered as a restarting step and simultaneously decreases the smoothness parameter $\beta$.
\vspace{-1ex}
\end{itemize}
The intuition behind our new strategy lies on the fact that when applied to~\eqref{eq:smoothed_cvx_prob} with a fixed $\beta$, APG gets $\BigO{\frac{1}{k^2}}$ rate, whereas ASGARD as presented in~\cite{TranDinh2015b} controls the parameters in such a way that the algorithm gets $\BigO{\frac{1}{k}}$ rate throughout its execution.
The idea is to take the advantage of the faster rate of APG for the inner loop while carefully adjusting the number of inner iterations and the smoothness parameter to get the same overall $\BigO{\frac{1}{k}}$ rate with better practical performance.
Our analysis also gives insights on the heuristic restart strategy outlined in~\cite{TranDinh2015b}.
For the sake of presentation and its flexibility for using Bregman distances in proximal operators, we choose Tseng's variant of APG \cite{tseng2008accelerated}. However, we can replace by another scheme such as FISTA \cite{Beck2009}.
We adaptively determine the number of inner iterations at each outer iteration. 
Therefore, there is no need to tune this parameter.
The outer  loop gradually decreases the smoothness parameter $\beta$ such that the algorithm is still guaranteed to converge to the true solution of \eqref{eq:ns_cvx_prob} or \eqref{eq:constr_cvx}.

\vspace{1ex}
\noindent\textbf{The algorithm:}
The complete algorithm is presented in Algorithm \ref{alg:A2}.

\begin{algorithm}[H]\caption{(\textrm{Self-Adaptive Double Loop ASGARD Algorithm})}\label{alg:A2}
\begin{normalsize}
\begin{algorithmic}[1]
\State {\hskip0ex}\textbf{Initialization:} 
	\vspace{0.75ex}
	\State Choose $\beta_0 > 0$, $\omega > 1$, a positive integer $m_0 \geq 1$, $\bar{x}^0 \in \mathbb{R}^p$, and $\dot{y}^0 \in\R^n$. 
	\vspace{0.75ex}
	\State Choose a Bregman distance $b_{\Yc}$ for $y$ and $d_{\Xc}$ for $x$.
	\vspace{0.75ex}	
	\State Set $K_0 \leftarrow 0$, $\hat{x}^0 \leftarrow \bar{x}^0$, and $\tau_0 := 1$.
	\vspace{0.75ex}	
	\State{\hskip0ex}\textbf{For}~{$s = 0$ {\bfseries to} $S_{\max} - 1$, \textbf{perform:}}
	\vspace{0.75ex}	
		\State {\hskip4ex}\textbf{For $j := 0$ to $m_s - 1$ perform}
		\vspace{0.75ex}	
		\State {\hskip8ex} Set $k {~~~~}\leftarrow K_s + j$.
		\vspace{0.75ex}	
		\State {\hskip8ex}\label{step:A2_mainstep} Update 
		$\left\{\begin{array}{ll}
		\tilde{x}^k & \leftarrow (1-\tau_k)\bar{x}^k + \tau_k \hat{x}^k  \vspace{1ex}\\
		\tilde{y}^{k+1} & \leftarrow \displaystyle\argmin_{y\in\Yc}\set{ g^{\ast}(y) - \iprods{A\tilde{x}^k, y} + \beta_s b_\Yc(y, \dot{y}^s)}\vspace{1ex}\\
		\hat{x}^{k+1}  & \leftarrow \Pc^{d_\Xc} _{\gamma_k f}\left(\hat{x}^k, A^\top\tilde{y}^{k+1}\right)~~~\text{with}~~\gamma_k \leftarrow \frac{\beta_s}{\Vert A\Vert^2 \tau_k}.
		\end{array}\right.$
		\label{eq:FISTA_step}
		\vspace{0.75ex}			
		\State {\hskip8ex}\label{step:xbar_next} Update $\bar{x}^{k+1}$ using one of the following two options:
		\[\left[\begin{array}{lll}
		\bar{x}^{k+1} & \leftarrow \tilde{x}^k + \tau_k(\hat{x}^{k+1} - \hat{x}^k) &~~~~~ \text{\textbf{Option 1:} Averaging step} \vspace{1ex}\\
		\bar{x}^{k+1} & \leftarrow  \Pc^{d_\Xc} _{\beta_s f/\norms{A}^2}\left(\tilde{x}^k, A^\top\tilde{y}^{k+1}\right) &~~~~~ \text{\textbf{Option 2:} Proximal step}.
		\end{array}\right. {\hskip-8ex}\]
		\State {\hskip8ex} Update $\tau_k \leftarrow \frac{2}{k-K_s +2}$.\label{step:A2_update_tau}
	\vspace{0.75ex}			
	\State {\hskip4ex} \textbf{End~for}
	\vspace{0.75ex}		
	\State {\hskip4ex} Update $K_{s+1} {~~}\leftarrow K_{s} + m_s$.
	\vspace{0.75ex}	
	\State {\hskip4ex} Restart $\bar{x}^{K_{s+1}} \leftarrow \hat{x}^{K_{s+1}} \equiv \hat{x}^{K_s + m_s}$.
	\vspace{0.75ex}		
	\State {\hskip4ex} Restart $\dot{y}^{s+1} {~}\leftarrow \prox_{\frac{1}{\beta_{s}}g^{\ast}}\left( \dot{y}^s + \frac{1}{\beta_s} A\bar{x}^{K_{s+1}}\right)$.\label{step:A2_up_center} 
	\vspace{0.75ex}		
	\State {\hskip4ex} Restart $\tau_{K_{s+1}}\leftarrow 1$, ~and update $\beta_s$ and $m_s$  by \eqref{eq:update_param1a}.\label{step:A2_update_beta_s}
	\vspace{0.75ex}		
\State{\hskip0ex}\textbf{End~for}
\end{algorithmic}
\end{normalsize}
\end{algorithm}

Algorithm~\ref{alg:A2} uses APG with \textbf{Option 1} at Step \ref{eq:FISTA_step} and has the same per-iteration complexity as Algorithm~\ref{alg:A1} except for the extra step, Step~\ref{step:A2_up_center}, where we update the dual center $\dot{y}^s$ at each outer loop iteration.
In general, the number of outer iterations is small as it is the number of restarting steps. Hence, Step~\ref{step:A2_up_center} does not significantly increase the overall computational cost of the entire algorithm. 
Note that $\bar{x}^{k+1}$ computed at Step \ref{step:xbar_next} using \textbf{Option 1} is a weighted averaging step.
To avoid this averaging, we can choose \textbf{Option 2}, which requires an additional generalized proximal operator of $f$.
Alternatively, we can replace Step~\ref{step:A2_mainstep} of Algorithm \ref{alg:A2} by the following FISTA step:
\begin{equation}\label{eq:fista_step}
\left\{\begin{array}{ll}
\tilde{y}^{k+1} & \leftarrow \displaystyle\argmin_{y\in\Yc}\set{ g^{\ast}(y) - \iprods{A\tilde{x}^k, y} + \beta_s b_\Yc(y, \dot{y}^s)}\vspace{1ex}\\
\bar{x}^{k+1}  & \leftarrow \Pc^{d_\Xc} _{\gamma_k f}\left(\tilde{x}^k, A^\top\tilde{y}^{k+1}\right)~~~\text{with}~~\gamma_k \leftarrow \frac{\beta_s}{\Vert A\Vert^2} \vspace{1ex}\\
\tilde{x}^{k+1} & \leftarrow \bar{x}^{k+1} + \frac{(1-\tau_k)\tau_{k+1}}{\tau_k}(\bar{x}^{k+1} - \bar{x}^k).
\end{array}\right.
\end{equation}
However, we need to replace the general Bregman distance $d_{\Xc}$ by an Euclidean distance $d_{\Xc}(\cdot, \dot{x}) := \frac{1}{2}\norms{\cdot - \dot{x}}^2$.
The scheme \eqref{eq:fista_step} allows us to compute $\bar{x}^k$ through $\Pc^{d_\Xc} _{\gamma_k f}(\cdot)$ instead of a weighted averaging step as with \textbf{Option 1}.

\vspace{1ex}
\noindent\textbf{Convergence guarantee:}
Now, we analyze convergence of Algorithm \ref{alg:A2} for solving~\eqref{eq:ns_cvx_prob} under Assumption~\ref{as:A1}.
Due to technical details, we separate the main theorem, Theorem \ref{th:convergence_of_A2}, and its proof into different sections.

\begin{theorem}\label{th:convergence_of_A2}
Let $\set{\bar{x}^{K_s}}$ be the sequence generated by Algorithm~\ref{alg:A2} and $\omega > 1$ be a given constant. 
Assume that \eqref{eq:ns_cvx_prob} satisfies Assumption~\ref{as:A1} and the parameters $\beta_s$ and $m_s$ are updated as 
\begin{equation}\label{eq:update_param1a}
\beta_{s+1} \leftarrow \tfrac{\beta_s}{\omega}~~~\text{and}~~~m_{s+1} \leftarrow \lfloor\omega(m_s+1) +  1 \rfloor - 1.
\end{equation}
where $m_0 \geq 1$ and $\beta_0 > 0$.
Let $\kappa_0 := m_0 + \frac{\omega}{\omega-1} > 0$. Then, we have
\begin{equation} \label{eq:A2_key_est2a}
P(\bar{x}^{K_{s+1}}) - P^{\star} \leq  \frac{\omega \kappa_0}{\beta_0\left[(\omega-1)K_{s+1} + \kappa_0\right]}\left[R_0^2 + \frac{\beta_0^2\omega D_{\Yc}}{(\omega-1)m_0}\right],
\end{equation}
where $R_0^2 := \frac{4\Vert A\Vert^2}{(m_0+1)^2} d_\Xc(x^\star, \bar{x}^{0}) + \beta_0^2 b_\Yc (y^\star, \dot{y}^0) $ and $D_{\Yc} := \sup \big\{b_\Yc( y, \dot{y}^s) \mid y\in\dom{g^{\ast}}, ~\forall s \geq 0 \big\}$.

Consequently, if $\dom{g^{\ast}}$ is bounded $($equivalently, g is Lipschitz continuous$)$, then Algorithm~\ref{alg:A2} achieves $\mathcal{O}\left( \frac{1}{K_s} \right)$ convergence rate at the last iterate of the outer loop, i.e. $\vert P(\bar{x}^{K_{s+1}}) - P^\star \vert \leq \BigO{\frac{1}{K_{s+1}}}$.
\end{theorem}

\noindent We make a few remarks about Theorem \ref{th:convergence_of_A2}.
\begin{itemize}
\vspace{-1ex}
\item The smoothness parameter $\beta$ is only updated at the outer loop but with a geometric rate, depending on the parameter $\omega$.
We can select different $\omega$ to observe the performance in particular applications. 
\item The convergence rate depends on both the prox-distance between $\bar{x}^0$ to $x^{\star}$ and $\dot{y}^0$ to $y^{\star}$ as well as $D_{\Yc}$, the prox-diameter of $\dom{g^{*}}$.
\item The convergence rate is given at the last iterate instead of the averaged sequence as often seen in other primal-dual methods \cite{Chambolle2011,Ouyang2014,Shefi2014,xu2017accelerated}.
\end{itemize}
\begin{remark}\label{re:for_every_k}
By using \eqref{eq:thm3ns_proof1b} in our analysis, we can show that $\set{ P(\bar{x}^k)}$ converges to $P^{\star}$ at the rate of $\BigO{\frac{1}{k}}$ for any $k\geq 1$ instead of $k = K_s$ at the outer loop only.
\begin{itemize}
\item If we use the averaging step of APG, then $\bar{x}^k$ is computed via a weighted averaging step of the inner loop.
\item However, if we use the proximal step in APG or the FISTA scheme \eqref{eq:fista_step}, then $\bar{x}^k$ is computed through the generalized proximal operator $\Pc^{d_\Xc} _{\beta_s f/\norms{A}^2}$.
This rate is fully non-ergodic.
\end{itemize}
\end{remark}

\beforesubsec
\subsection{\textbf{Application to constrained convex optimization}}
\aftersubsec
One important application of nonsmooth optimization is linearly constrained convex optimization.
Most of the works on smoothing, including Nesterov's seminal work~\cite{Nesterov2005c}, need to know the diameter of the dual domain $D_\Yc$ to set the smoothness parameter. 
For the case of linear equality constraints, $D_\Yc$ is unbounded, therefore, these algorithms cannot be applied.

We will illustrate now how to apply our algorithm to constrained convex optimization problem~\eqref{eq:constr_cvx}, without any dependence on $D_\Yc$.
In this section, we require the Bregman distance used in smoothing for the dual variables to have Lipschitz gradient.
Under this condition, we have
\begin{equation}
b_\Yc(y, \dot{y}) \leq \frac{L_{b_\Yc}}{2} \norms{y - \dot{y}}_{\Yc}^2.
\end{equation}
Let us define $g(Ax) := \delta_{\Kc}(Ax - b)$ the indicator function of $\Kc$, where $\Kc$ satisfies Assumption~\ref{as:A1b}. 
Then, we can write $g$ as
\begin{equation}\label{eq:gAx3}
g(Ax) := \sup_{y\in\R^n}\set{ \iprods{Ax - b, y} - s_{\Kc}(y)},
\end{equation}
where $s_{\Kc}(y) := \sup_{u\in\Kc}\iprods{y, u}$ is the support function of $\Kc$.
In this case, the smooth function $g_{\beta}(Ax;\dot{y})$ becomes
\begin{equation}\label{eq:gAx3_beta}
g_{\beta}(Ax;\dot{y}) := \max_{y\in\R^n}\Big\{ \iprods{Ax - b, y} - s_{\Kc}(y) - \beta b_\Yc( y, \dot{y}) \Big\}.
\end{equation}

\begin{example}\label{ex:cone_constr_l2}
Suppose that $b_\Yc(x, \dot{x})=\frac{1}{2}\Vert x - \dot{x}\Vert_2^2$. Then, the function $g_{\beta}(\cdot;\dot{y})$ defined by \eqref{eq:gAx3_beta} can be written as
\begin{equation}\label{eq:gAx3_beta_b}
g_{\beta}(Ax;\dot{y}) = \frac{1}{2\beta}\kdist{\Kc}{Ax - b + \beta\dot{y}}^2 - \frac{\beta}{2}\Vert\dot{y}\Vert^2.
\end{equation}
Moreover, the solution $y^{\ast}_{\beta}(Ax;\dot{y})$ of the maximization problem in \eqref{eq:gAx3_beta} is given in a closed form as
\begin{equation}\label{eq:ystar_beta}
y^{\ast}_{\beta}(Ax;\dot{y}) = \dot{y} + \tfrac{1}{\beta}\left(Ax - b - \kproj{\Kc}{Ax - b + \beta\dot{y}}\right), 
\end{equation}
where $\kproj{\Kc}{\cdot}$ denotes the projection onto $\Kc$. 

In particular, if $\Kc$ is a cone, then $y^{\ast}_{\beta}(Ax;\dot{y}) = \kproj{-\Kc^{\ast}}{\dot{y} + \tfrac{1}{\beta}(Ax - b)}$, where $\Kc^{\ast}$ is the dual cone of $\Kc$.
The dual step for computing $\tilde{y}^k$ at the second line of Step~\ref{eq:FISTA_step} of Algorithm~\ref{alg:A2} becomes
\begin{align}\label{eq:y_tilde}
\tilde{y}^{k+1} &\leftarrow \dot{y}^s + \tfrac{1}{\beta_s}\left(A\tilde{x}^k - b - \kproj{\Kc}{A\tilde{x}^k - b + \beta_s\dot{y}^s}\right) \\ &= \tfrac{1}{\beta_s}\kproj{-\Kc^{\ast}}{A\tilde{x}^k - b+ \beta_s\dot{y}^s} \notag.
\end{align}

\end{example}

The following lemma provides a key estimate for the optimality condition of \eqref{eq:constr_cvx}, whose proof is given in Appendix~\ref{apdx:le:opt_cond_constr_cvx2}.

\begin{lemma}\label{le:opt_cond_constr_cvx2}
Let $S_{\beta}(\bar{x};\dot{y}) := f(\bar{x}) + g_{\beta}(A\bar{x};\dot{y}) - f(x^{\star})$ and $\beta_b:=\beta L_{b_\Yc}$. Then:
\begin{equation}\label{eq:approx_opt_cond}
\left\{\begin{array}{ll}
f(\bar{x}) - f^{\star}                        &\geq  \beta_b\iprods{\dot{y}, y^{\star}} -\Vert y^{\star}\Vert\kdist{\Kc}{A\bar{x} - b + \beta_b \dot{y}} \vspace{1ex}\\
f(\bar{x}) - f^{\star}                               &\leq S_{\beta}(\bar{x};\dot{y}) -\frac{1}{2\beta_b} \kdist{\Kc}{A\bar{x} - b + \beta_b \dot{y}} + \frac{\beta_b}{2}\Vert \dot{y}\Vert^2\vspace{1ex}\\
\kdist{\Kc}{A\bar{x} - b + \beta_b\dot{y}} &\leq \beta_b\Big[\Vert y^{\star}\Vert  + \big(\Vert \dot{y} - y^{\star}\Vert^2 + \frac{2}{\beta_b}S_{\beta}(\bar{x};\dot{y})\big)^{1/2} \Big].
\end{array}\right.
\end{equation}
Here, the term $\Vert \dot{y} - y^{\star}\Vert^2 + \frac{2}{\beta_b}S_{\beta}(\bar{x};\dot{y}) \geq \Vert y^{\ast}_{\beta}(A\bar{x};\dot{y}) - \dot{y}\Vert^2 \geq 0$. In addition, we have the following bound for any $\beta$:
\begin{equation}\label{eq:constr_san}
\begin{array}{ll}
\kdist{\Kc}{A\bar{x}-b}-\beta\big(\Vert \dot{y}-y^\star \Vert &+ \Vert y^\star \Vert\big) \leq \kdist{\Kc}{A\bar{x}-b-\beta\dot{y}} \vspace{1ex}\\
&\leq \kdist{\Kc}{A\bar{x}-b} + \beta\left(\Vert \dot{y}-y^\star \Vert + \Vert y^\star \Vert\right).
\end{array}
\end{equation}
\end{lemma}

Now, we apply Algorithm~\ref{alg:A1} to solve the constrained convex problem \eqref{eq:constr_cvx}.
Then, the following steps are changed:
\begin{itemize}
\item The APG scheme at Steps \ref{step:A2_mainstep} and \ref{step:xbar_next} is replaced by the FISTA scheme \eqref{eq:fista_step}.
\item The dual step for computing $\tilde{y}^{k+1}$ of Algorithm~\ref{alg:A2} becomes
\begin{equation}\label{eq:y_tilde}
\tilde{y}^{k+1} \leftarrow \arg\max_{y\in\mathbb{R}^n} \left\{ \langle A\tilde{x}^k -b, y \rangle - s_\Kc(y) - \beta_s b_\Yc(y, \dot{y}^s) \right\}.
\end{equation}
\item The update rule of $\tau_k$ and $\beta_k$ is changed as \eqref{eq:update_param2a}.
\end{itemize}
Combining the result of Theorem \ref{th:convergence_of_A2} and Lemma \ref{le:opt_cond_constr_cvx2}, we obtain the following convergence result of this new variant of Algorithm \ref{alg:A2}.

\begin{theorem}\label{th:convergence_A2b}
Assume that Assumption~\ref{as:A1b} holds.
Let $\set{ \bar{x}^{K_s} }$ be the sequence generated by Algorithm~\ref{alg:A1} for solving \eqref{eq:constr_cvx} using \eqref{eq:y_tilde} for $\tilde{y}^{k+1}$.
Let the parameters $\beta_s$ and $m_s$ be updated as
\begin{equation}\label{eq:update_param2a}
\left\{\begin{array}{ll}
m_{s+1} &\leftarrow \lfloor\omega(m_s+1) + 1 \rfloor - 1~~~\text{with}~~~m_0 > \frac{1}{\omega - 1}\vspace{1ex}\\
\beta_{s+1} &\leftarrow \frac{\beta_s(m_{s+1}+1)}{\omega\sqrt{m_{s+1}(m_{s+1}+3)}}.
\end{array}\right.
\end{equation}
Then, we have
\begin{equation}\label{eq:key_estimate3b}
{\!\!\!\!\!\!\!}\left\{\begin{array}{ll}
f(\bar{x}^{K_{s+1}}) - f^{\star} &{\!\!\!}\geq -\Vert y^{\star}\Vert\kdist{\Kc}{A\bar{x}^{K_{s+1}}  - b}  - \frac{2\sqrt{2} \omega \beta_0 L_{b_\Yc} \kappa_0 \Vert y^\star \Vert R_0 }{\rho_0\left[(\omega-1)K_{s+1}+\kappa_0\right]}\vspace{1ex}\\
f(\bar{x}^{K_{s+1}}) - f^{\star} &{\!\!\!}\leq \frac{\omega\kappa_0R_0^2}{\rho_0\left[(\omega-1)K_{s+1} + \kappa_0\right]} + \frac{\omega\beta_0 L_{b_\Yc} \kappa_0 }{2\left[(\omega-1)K_{s+1}+\kappa_0\right]} \left( \Vert y^\star \Vert ^2 + \frac{2R_0^2}{\rho^2} \right) \vspace{1ex}\\
\kdist{\Kc}{A\bar{x}^{K_{s+1}} {\!}- b} &\leq \frac{\omega\beta_0 L_{b_\Yc}\kappa_0}{\left[(\omega-1)K_{s+1} + \kappa_0\right]}\left[ 2\Vert y^{\star}\Vert + \left(2\sqrt{2} + \sqrt{\frac{2}{L_{b\Yc}}}\right)\frac{R_0}{\rho_0} \right],
\end{array}\right.{\!\!\!\!\!\!}
\end{equation}
where $y^{\star}$ is any dual solution of \eqref{eq:ns_cvx_dual}, and
\begin{equation*}
\left\{\begin{array}{ll}
\rho_0 &:= \beta_0\left(1 - \frac{1}{(\omega-1)m_0}\right) \vspace{1ex}\\
\kappa_0 &:= m_0 + \frac{\omega}{\omega-1}\vspace{1ex}\\
R_0 &:= \left[\frac{4\Vert A\Vert^2}{(m_0+1)^2} d_\Xc(x^\star, \bar{x}^{0})  + \frac{\beta_0^2m_0(m_0+3)}{(m_0+1)^2}b_\Yc(  y^\star, \dot{y}^0)\right]^{1/2}.
\end{array}\right.
\end{equation*}
Consequently, Algorithm~\ref{alg:A2} achieves an $\BigO{\frac{1}{{K_s}}}$ convergence rate in a \textbf{non-ergodic sense}, i.e., $\vert f(\bar{x}^{K_s}) - f^{\star} \vert \leq \BigO{\frac{1}{{K_s}}}$ and $\kdist{\Kc}{A\bar{x}^{K_s}-b} \leq \BigO{\frac{1}{{K_s}}}$.
\end{theorem}

%
\beforesubsec
\subsection{\textbf{Extension to Composite  Case with Three Objective Terms}}\label{sec:extension}
\aftersubsec
It is straightforward to apply Algorithm \ref{alg:A2} in the presence of a smooth term in the objective. The problem template we focus on in this section is
\begin{equation}\label{eq:three_comp}
F^{\star}  := \min_{x\in\mathbb{R}^p}\Big\{F(x) := f(x) + g(Ax) + h(x) \Big\},
\end{equation}
where $f$ and $g$ are as described in Assumption \ref{alg:A1} and $h$ is a differentiable function with $L_h$-Lipschitz gradient. In this case, only Step~\ref{eq:FISTA_step} in Algoritm~\ref{alg:A2} needs to be modified as follows (see also in \cite{van2017smoothing}):
\begin{equation*}
\hat{x}^{k+1} \leftarrow \Pc^{d_\Xc} _{\gamma_k f}\left(\hat{x}^k, \nabla h(\tilde{x}^k)+ A^\top\tilde{y}^{k+1}\right)~~~\text{with}~~\gamma_k \leftarrow \tfrac{\beta_s}{\tau_k\left( \Vert A\Vert^2 + \beta_s L_h\right)}.
\end{equation*}
Note that, this modification only changes the analysis of the inner loop as in~\cite{van2017smoothing} which does not affect our analysis of the outer loop.
In addition, using $L_h$ in the stepsize is not restrictive. 
When the Lipschitz constant is not known, line search strategies can be employed, see \cite{van2017smoothing} for more details.
The convergence of this variant is still guaranteed by Theorem \ref{th:convergence_of_A2} but the quantity $R_0^2$ will depend on $L_h$.
We omit the details of this result here for succinctness.

\beforesubsec
\subsection{\textbf{Better complexity bounds for Linear Programming}}\label{sec:complexity}
\aftersubsec
As an application of our theory, we  analyze the overall complexity of our algorithm for linear programming:
\begin{equation}\label{eq:lp}
\min_{x\in\mathbb{R}^p}\Big\{ c^{\top}x ~~ \text{ s.t. }~~ Ax=b, ~x\geq 0 \Big\},
\end{equation}
where we define $f(x):=  c^{\top}x  + \delta_{\set{x \geq 0}}(x)$.

First-order methods for linear programming have been widely studied in the literature. 
However, it is generally not preferred since nonsmooth optimization methods have $\BigO{\frac{1}{\varepsilon^2}}$ complexity~\cite{renegar2014efficient,wang2017new}, where $\varepsilon$ is a desired accuracy. 
\cite{renegar2014efficient} proposed applying Nesterov smoothing to linear programs.
Unfortunately, since Nesterov smoothing does not apply to indicator function of the linear constraints, their per-iteration complexity requires projection to the domain defined by the linear constraint which is computationally expensive. 
We propose applying smoothing to the linear constraint since our theory supports it to come up with an algorithm with much cheaper iterations.

Given $f(x) := c^{\top}x + \delta_{\{x \geq 0\}}(x)$, we define $\tilde{x}^{\ast}$ to be an $\varepsilon$-solution to \eqref{eq:lp} if
\begin{equation*}
\vert f(\tilde{x}^{\ast}) - f(x^{\star}) \vert \leq \varepsilon \text{ and } \norms{A\tilde{x}^{\ast} - b} \leq \varepsilon.
\end{equation*}
Using the bounds from Theorem~\ref{th:convergence_A2b}, we can derive the iteration complexity of our method for linear programming.
Let us use $b_{\Yc}(y_1, y_2) := \frac{1}{2}\norms{ y_1 - y_2 }_2^2 $ and $d_{\Xc}(x_1, x_2) = \frac{1}{2} \norms{x_1 - x_2}_2^2$ and let $\bar{x}^0$ and $\dot{y}^0$ to be all zero vectors.
Then, we can achieve an $\varepsilon$-solution $x^K$ of \eqref{eq:lp} for any $K \geq C K_{\varepsilon}$, where
\begin{align}
K_{\varepsilon} = \frac{\max\left( \Vert y^\star \Vert ^2, R_0^2 \right)}{\varepsilon} = \frac{\max\left( \Vert y^\star \Vert ^2, \Vert A \Vert ^2 \Vert x^\star \Vert ^2 \right)}{\varepsilon},
\end{align}
and $C$ hides dimension independent quantities. 
The per-iteration complexity of our method is dominated by applying $A$ matrix which has $\mathrm{nnz}(A)$ complexity. 
Combining this with the iteration complexity, we can estimate the overall complexity of our method for linear programming as
\begin{equation}
\mathcal{O} \left( \max\left(\Vert A \Vert ^2 \Vert x^\star \Vert ^2, \Vert y^\star \Vert \right)\mathrm{nnz}(A) \left( \tfrac{1}{\varepsilon} \right) \right).
\end{equation}

In~\cite{renegar2014efficient}, the author focuses on a specific case of linear programming where the primal domain diameter can be bounded.
Here, we focus on general linear programs, therefore, it is fairer to compare our complexity with~\cite{wang2017new}.
Compared with~\cite{wang2017new}, from the overall complexity, we remove the terms $a_m := \max _i \| A_i \|$ and $\theta_{S^\ast}^2$ where $\theta_{S^\ast}^2$ is the Hoffman bound and $S^\ast$ is the solution set.
One thing to note here is that~\cite{wang2017new} has a better dependence in terms of $\epsilon$ since they have linear convergence.
Our method, has a sublinear rate for accuracy, but a better dimension dependence which is the main bottleneck in large scale linear programming.

\beforesec
\section{Numerical experiments}\label{sec:num_exp}
\aftersec
We will test standard ASGARD~\cite{TranDinh2015b,van2017smoothing}, ASGARD with restart~\cite{TranDinh2015b,van2017smoothing} and standard Chambolle-Pock's algorithm~\cite{Chambolle2011} on the following problems. 
Note that when there is a smooth term in the objective, we use the version of Chambolle-Pock which linearizes the smooth term, which is also known in the literature as Vu-Condat's algorithm~\cite{vu2013splitting,condat2013primal}.
We omit HOPS~\cite{xu2016homotopy} from the comparisons because it does not apply to Basis pursuit and Markowitz's portfolio optimization problems due to the unboundedness of the dual domain. For the $\ell_1$-SVM example, we observed it to be extremely slow and difficult to tune for different datasets.
In all the experiments, we have used the standard $b_{\Yc}(y_1, y_2) = \frac{1}{2} \| y_1 - y_2 \|_2^2$ and $d_{\Xc}(x_1, x_2) = \frac{1}{2} \| x_1 - x_2 \|_2^2$ for smoothing and computing the proximal operators for fair comparison with other methods which do not allow Bregman distances.
In the sequel, we refer to our algorithm as ASGARD-DL.
In some cases, we also compare with ADMM and its variants.

The parameters are set as follows.
For Chambolle-Pock's method, we set its step-size $\sigma = \tau = \frac{1}{\norms{A}}$, where $A$ is the linear operator in \eqref{eq:ns_cvx_prob}.
For ASGARD-DL, we choose $\omega := 1.2$ and $m_s := 6$ which give us comparable performance.
For restarting ASGARD, we set the restarting frequency to be $s = 10$ in all experiments.

\beforesubsec
\subsection{\textbf{Convergence guarantees: Ergodic vs. Non-ergodic}}\label{subsec:example1}
\aftersubsec
ALM, ADMM and Chambolle-Pock methods have the convergence rate guarantees in an ergodic sense. 
That is, they have the rate guarantees only on the averaged iterate sequence. 
In contrast, our guarantees are for the last iterate of the algorithm.
To illustrate the importance between these two, we consider two synthetic problems in this section.
The first one is a square root LASSO problem widely studied in the literature, which is given by:
\begin{equation*}
F^{\star} := \min_{x\in\R^p} \Big\{ F(x) := \Vert Ax-b \Vert_2 + \lambda \Vert x \Vert_1 \Big\},
\end{equation*}
where $A\in\R^{n\times p}$ is generated using a Gaussian distribution and is normalized such that column norms are equal to 1.
Given a groundtruth vector $x^\natural$, we generate the observations as $b=Ax^\natural + \sigma \mathbf{n}$, where $\mathbf{n}$ is a noise vector generated by a standard Gaussian distribution and $\sigma = 0.01$.
We set $\lambda = 0.03$ which is tuned to get a good recovery of $x^\natural$.

\begin{figure}[hpt!]
\centering
\hspace{-2.5mm}\includegraphics[width=.5\columnwidth]{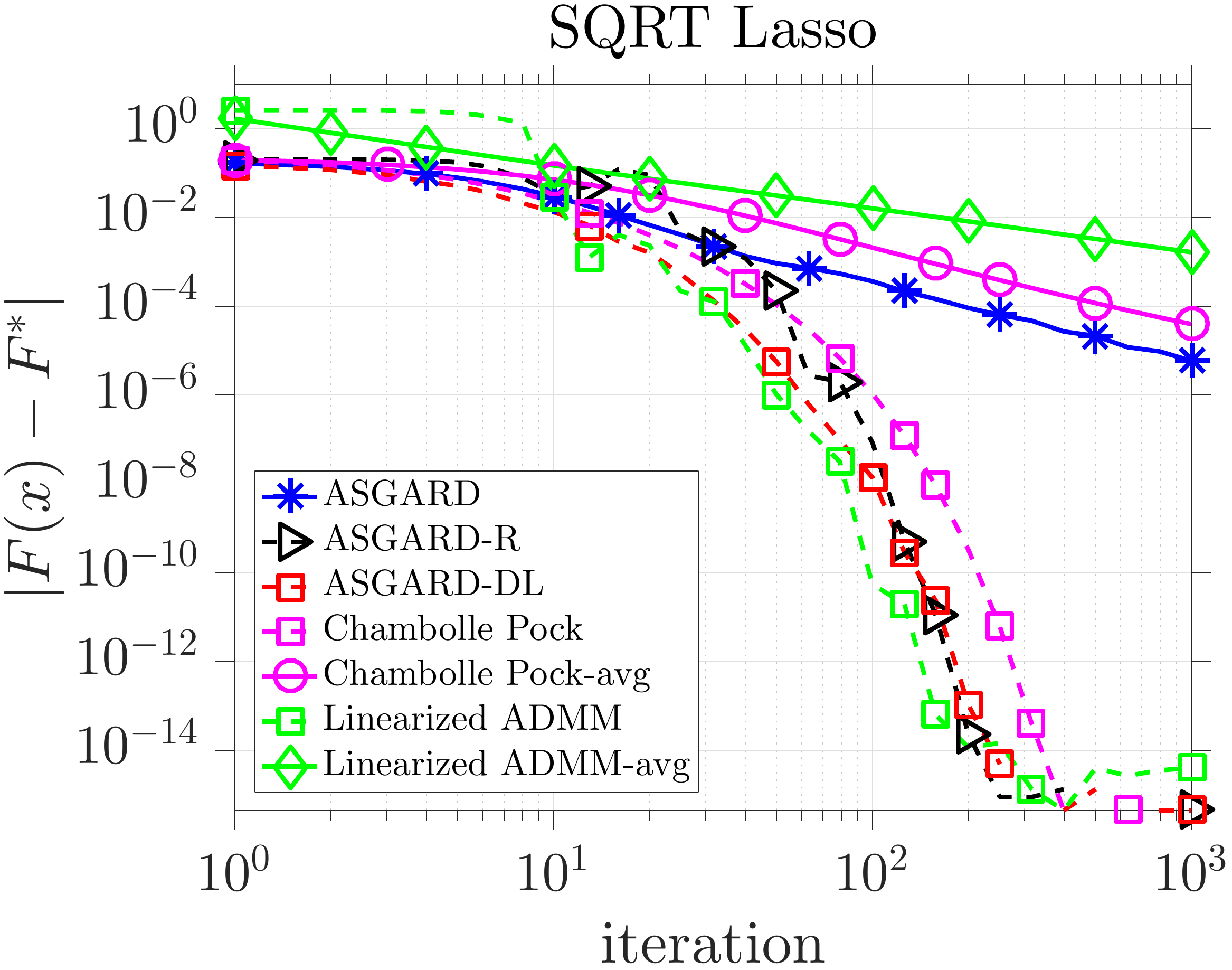} \hspace{0.1cm}
\hspace{-2.5mm}\includegraphics[width=.5\columnwidth]{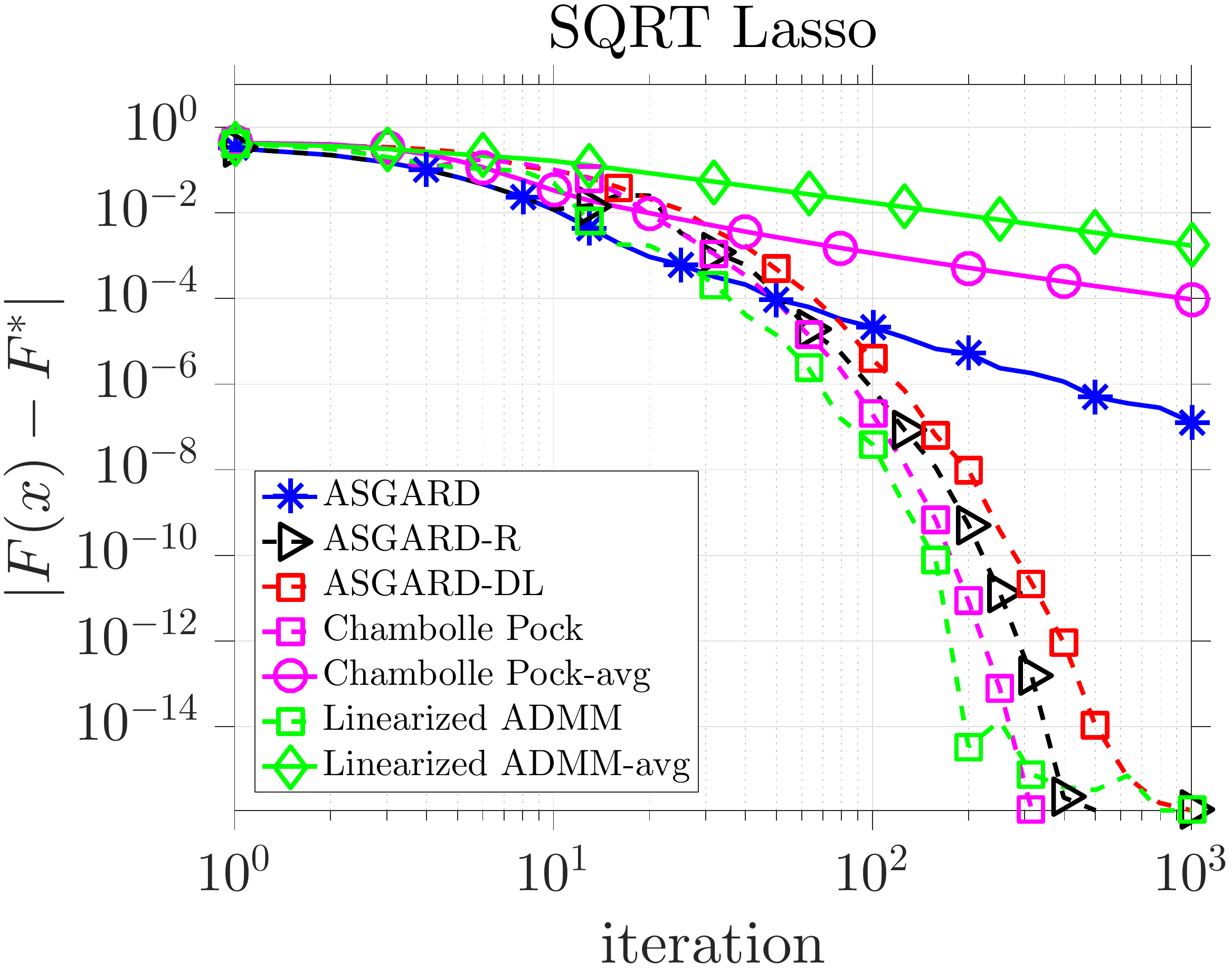}
\caption{Performance of 5 algorithms for solving square root LASSO problem. Left: $\sigma=0.1, \lambda=0.04$, Right: $\sigma=0.01, \lambda=0.03$}
\label{fig:sqrt_lasso}
\vspace{-4ex}
\end{figure}

In this experiment, we test the ergodic and non-ergodic variants of Linearized ADMM (in the sense that the augmented term in the Lagrangian is linearized)~\cite{gao2017first}, and Chambolle-Pock's algorithm~\cite{Chambolle2011}.
These methods have convergence guarantees for their last iterates, however, their rate guarantees only apply to the averaged sequence.
Moreover, they are very successful to solve this type of problems as can be seen from the litterature.
The behavior of the algorithms is given in Figure~\ref{fig:sqrt_lasso}.

As can be seen in Figure~\ref{fig:sqrt_lasso}, last iterates of Linearized ADMM and Chambolle Pock's algorithms seem to have the best performance.
However, the averaged iterates for which the methods have the rate guarantees shows the slowest convergence behavior.
Our method has the same rate as restarted ASGARD which does not have any convergence guarantees.

To illustrate the behavior of the last iterates of Linearized ADMM and Chambolle Pock's algorithm, we consider a degenerate linear program which is also studied in~\cite{TranDinh2015b}:
\begin{equation*}
\min_{x\in\R^p}\Big\{ h(x) := 2x_p ~\mid ~~ \sum_{k=1}^{p-1}x_k = 1, ~~ x_p - \sum_{k=1}^{p-1}x_k = 0 ~~(2\leq j \leq n), ~~ x_p \geq 0 \Big\}.
\end{equation*}
The second inequality is repeated $n-1$ times which causes the problem to be degenerate.
We define the linear constraint as
\begin{equation*}
Ax:= \Big[ \sum_{k=1}^{p-1}x_k ,~~ x_p - \sum_{k=1}^{p-1}x_k ,\cdots,~~x_p - \sum_{k=1}^{p-1}x_k \Big]^\top.
\end{equation*}
We have $b := \left( 1, 0, \cdots, 0 \right)^\top \in \R^n$.
We map the problem to our template in~\eqref{eq:three_comp} as $f(x) := \delta_{\set{x_p \geq 0}}(x_p)$, $g(x) := \delta_{\set{b}}(Ax)$, and $h(x) := 2x_p$.
For this problem, we pick $p=10$ and $n = 200$.

\begin{figure}[hpt!]
\vspace{0ex}
\centering
\hspace{-2.5mm}\includegraphics[width=1\columnwidth]{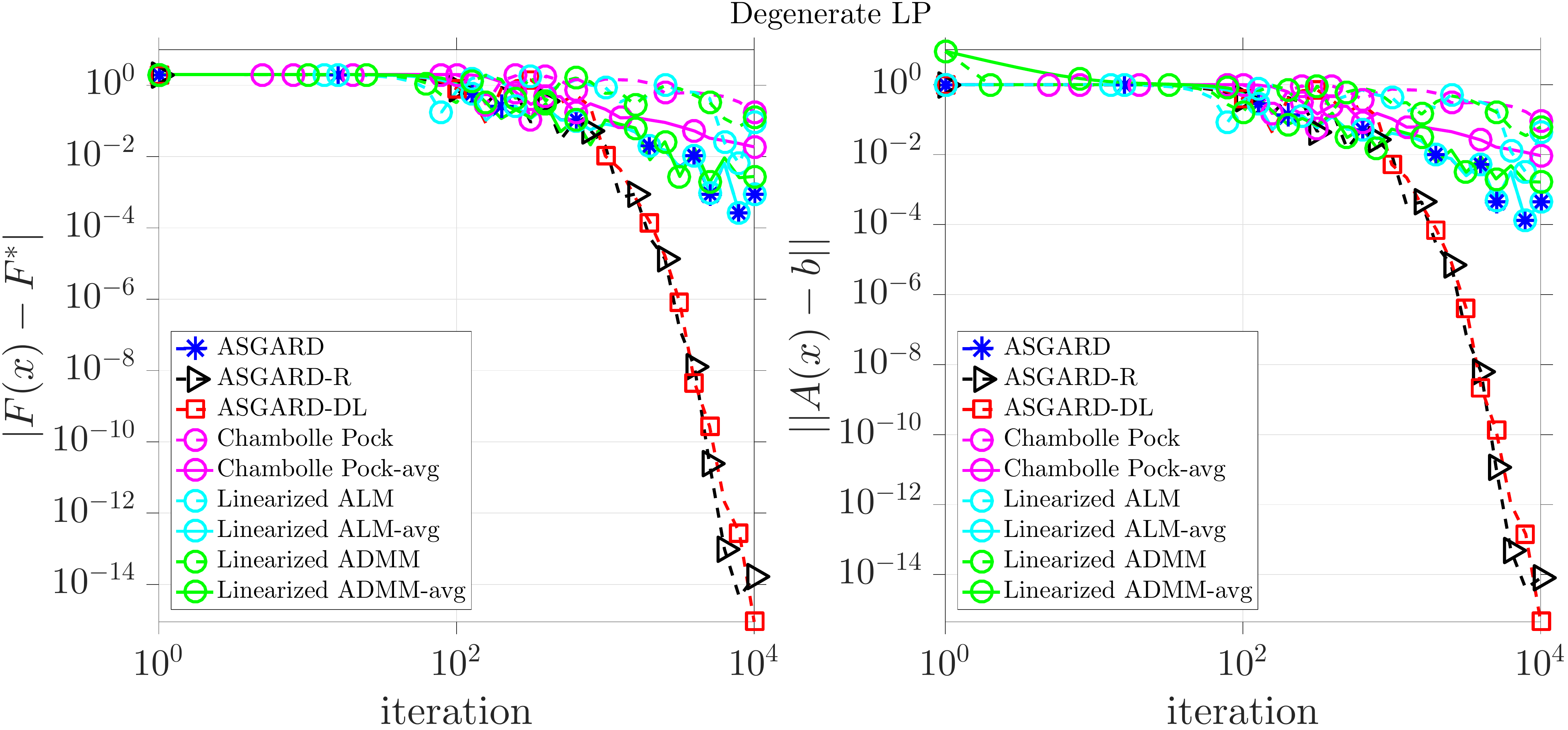}
\caption{Performance of 6 algorithms for solving the degenerate linear program.}
\label{fig:deg_lp}
\vspace{-4ex}
\end{figure}

In addition to Linearized ADMM and Chambolle-Pock's algorithm, we also include linearized ALM~\cite{gao2017first} to solve this example.
The result of this test is given in Figure \ref{fig:deg_lp}, where $F(x) = h(x)$.

As can be seen from Figure~\ref{fig:deg_lp}, Linearized ADMM, Linearized ALM and Chambolle-Pock's algorithm can get extremely slow where our algorithm and ASGARD with restart makes progress and converges to optimal value with a very high accuracy, and beyond the theoretical rate guarantee.

\beforesubsec
\subsection{\textbf{Basis Pursuit for recovering Bag-of-Words  of text documents}}
\aftersubsec
We first consider a basis pursuit problem which is used in signal/image processing, statistics, and machine learning~\cite{chen2001atomic,donoho2006compressed,candes2006robust}:
\begin{equation}\label{eq:bp}
\min_{x\in\R^{ p}}\big\{ F(x) := \Vert x \Vert _1 ~\mid~ Ax=b \big\},
\end{equation}
where $A\in\R^{n\times p}$ and $b\in\R^n$. 
This problem clearly fits into our template~\eqref{eq:ns_cvx_prob} by mapping $f(\cdot) = \norms{\cdot}_1$ and $g(\cdot) = \delta_{\set{ b }}(\cdot)$. 
It is also a special case of \eqref{eq:constr_cvx} with $\Kc = \set{\boldsymbol{0}}$.
Proximal operators of both terms are given in a closed form.

We apply this model to  text processing.
In~\cite{arora2018a}, the authors proposed using basis pursuit formulation to obtain bag-of-words representation from the unigram embedding representation of a text.
The setting can be briefly described as the following: For any word $w$, there exists a word vector $v_w \in \R^n$. 
For a given text document $\set{ w_1, \cdots, w_T }$, one defines the unigram embedding as $\sum_{i=1}^T v_{w_i}$. 
It is easy to see that unigram embeddings can be written as a linear system $Ax$ where $A\in\R^{n\times p}$ contains $v_{w_i}$ in the $i^{th}$ column and $x\in\R^p$ is the bag-of-words vector which counts the number of occurances of words in a text. 
This application is considered in text processing applications to obtain the original text document given the unigram embeddings~\cite{white2016generating}.

For this experiment, we have used the movie review dataset of~\cite{maas-EtAl:2011:ACL-HLT2011}.
We have selected $4$ different documents and computed the unigram embeddings using pre-trained word embeddings from \texttt{GloVe} \cite{pennington2014glove} with $n=50$ as the dimension of the word vectors and restricted the vocabulary size to $p=10,000$ for getting faster results with all algorithms.

\begin{figure}[hpt!]
\centering
\hspace{-2.5mm}\includegraphics[width=.95\columnwidth, height=.4\columnwidth]{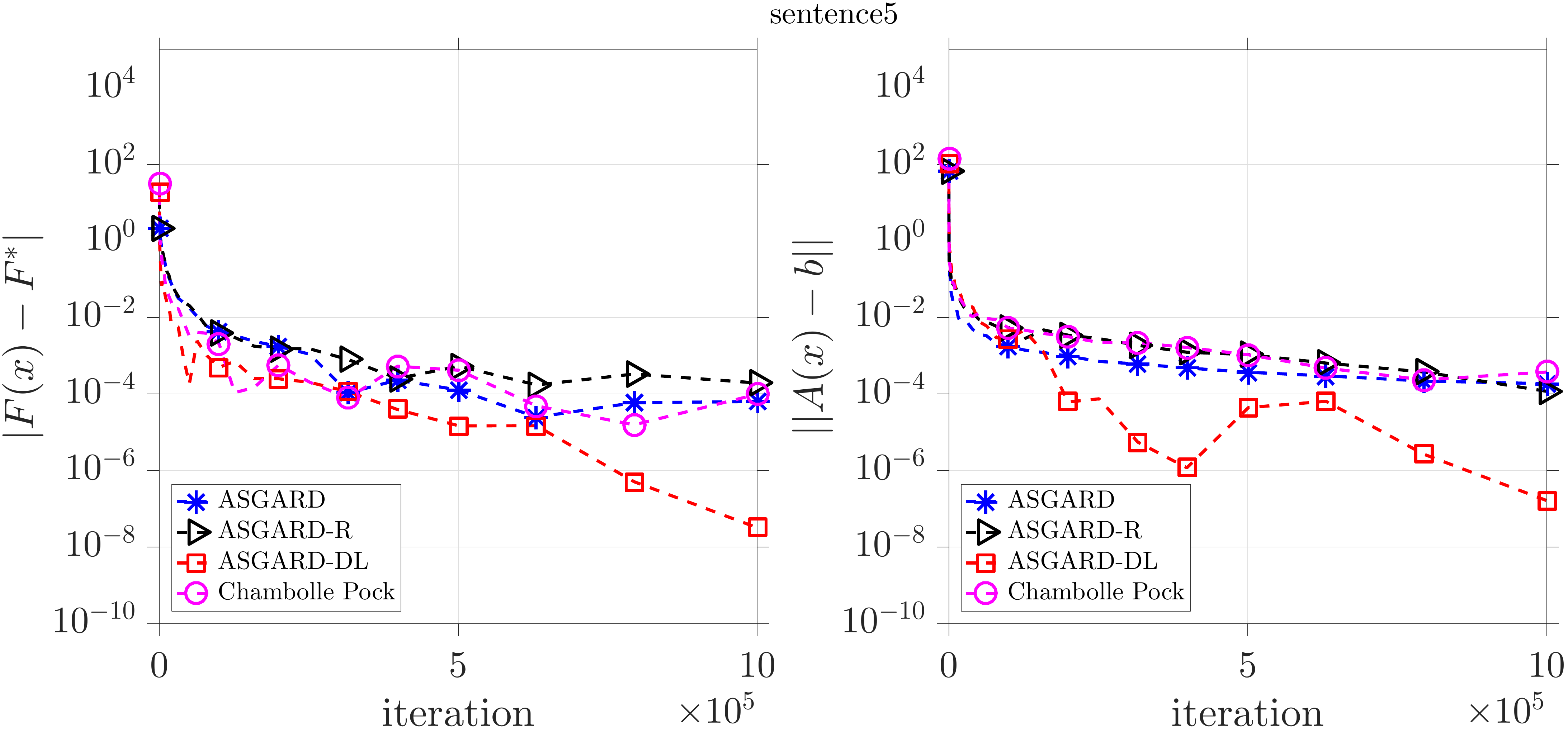} \\
\hspace{-2.5mm}\includegraphics[width=.95\columnwidth, height=.4\columnwidth]{figures/sentence5} \\
\hspace{-2.5mm}\includegraphics[width=.95\columnwidth, height=.4\columnwidth]{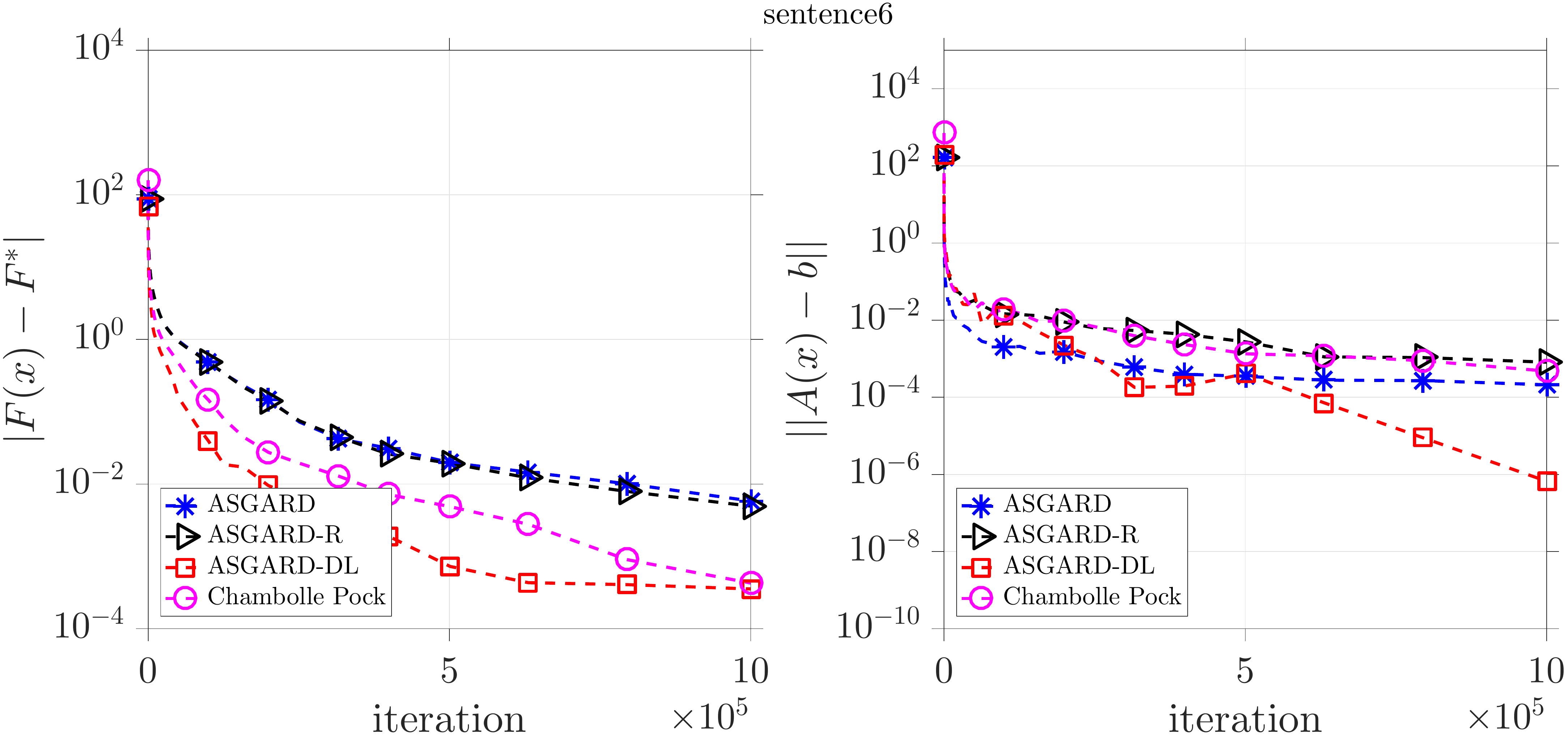} \\
\hspace{-2.5mm}\includegraphics[width=.95\columnwidth, height=.4\columnwidth]{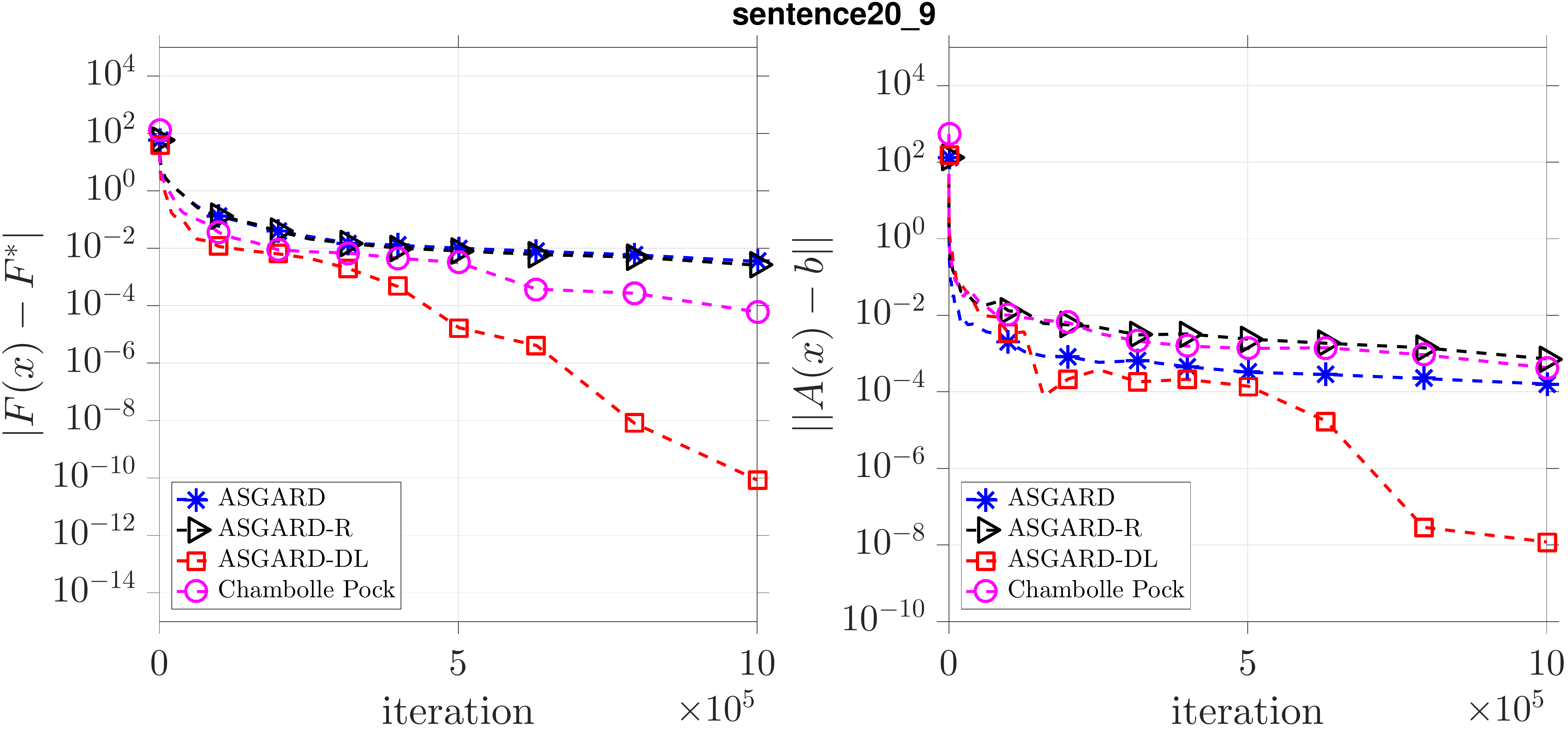}
\caption{Performance of 4 algorithms for solving basis pursuit for 4 text documents.}
\label{fig:bp}
\vspace{-4ex}
\end{figure}

We have applied $4$ methods to solve~\eqref{eq:bp} for $4$ different documents. 
Here, the parameter $\beta_0$ in ASGARD, ASGARD-restart, and ASGARD-DL is set to $\beta_0 := 10\norms{A}$.
Note that this choice is not optimal, but give us reasonable results in all test.
The results are compiled in Figure~\ref{fig:bp}.

As we can observe from Figure \ref{fig:bp}, our new algorithm works quite well and is comparable with state-of-the-art methods for low accuracy.
It outperforms them if we run the algorithms long enough to get more accurate solutions than $\varepsilon = 10^{-5}$ both in objective residual and feasibility.
Note that \eqref{eq:bp} is fully nonsmooth, and $A$ is non-orthogonal. 
If we apply ADMM to solve \eqref{eq:bp}, then it requires to solve a general convex subproblem, or a linear system, which has higher per-iteration complexity than four methods we used in this example.

\beforesubsec
\subsection{\textbf{The $\ell_1$-Regularized Least Absolute Deviation Problem (LAD)}}
\aftersubsec
Our second example is the  $\ell_1$-regularized least absolute deviation regression problem, also known as LAD-Lasso in the literature. 
It is known that when the noise has a heavy tailed distribution such as Laplace distribution, LAD-Lasso is more robust to the outliers \cite{wang2007robust}. 
The optimization model of this problem is
\begin{equation}
\min_{x\in\R^p}\Big\{ F(x) := \norm{Ax - b}_1 + \lambda \norm{x}_1 \Big\},
\end{equation}
where $A\in\R^{n \times p}$ is generated according to a normal distribution and the noise $\mathbf{n}\in\R^n$ is generated by Laplace$(0,1)$ distribution.
We generate an observed vector $b := Ax^{\natural} + \sigma \mathbf{n}$, where $\sigma := 0.1$ and $x^{\natural}$ is a $s$-sparse vector of ground-truth.
We choose $\lambda := 1/n$ for the regularization parameter, which gives us a good recovery of $x^{\natural}$.

\begin{figure}[hpt!]
\centering
\begin{tabular}{ccccc}
\hspace{-2.5mm}\includegraphics[width=0.33\columnwidth]{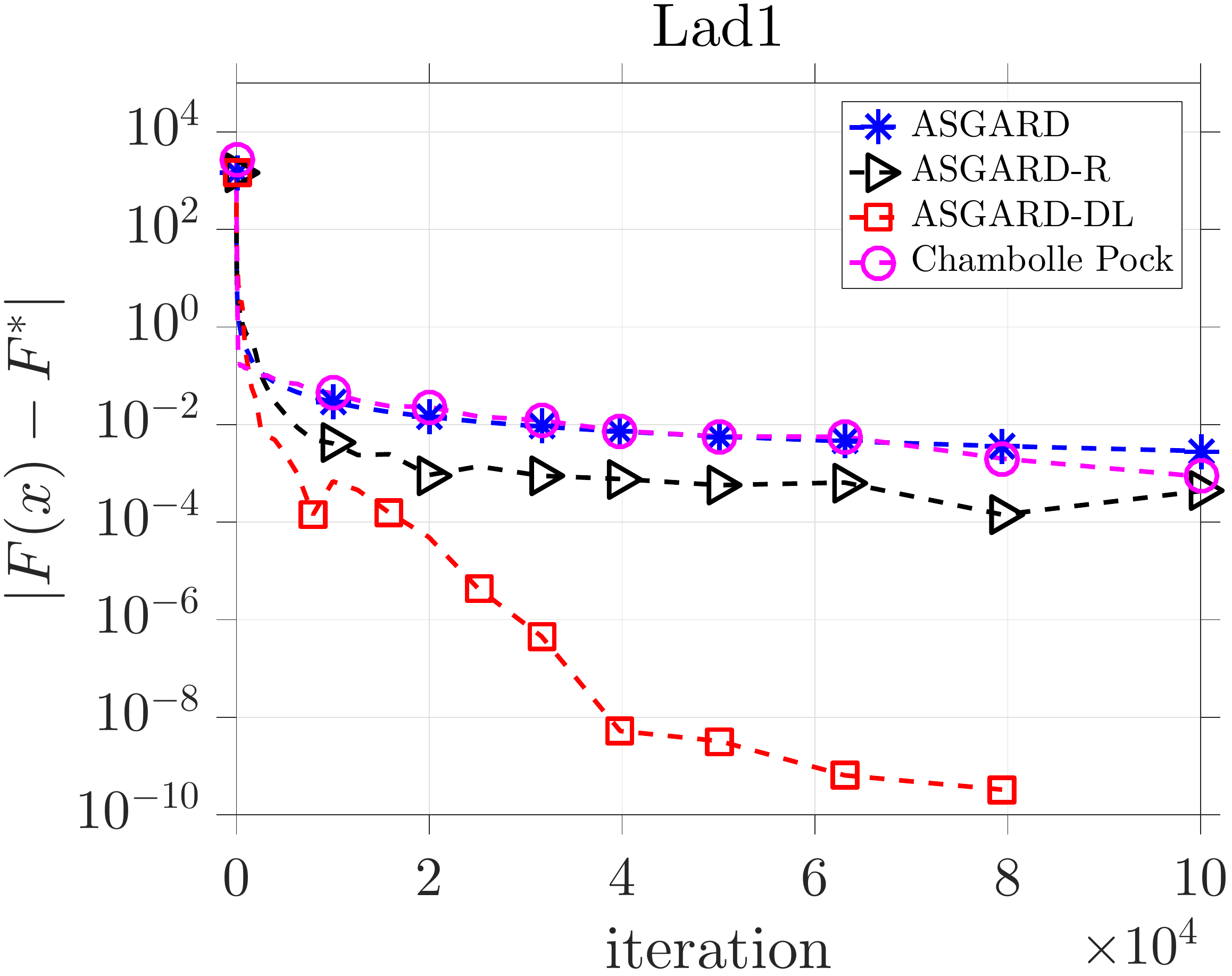} &
\hspace{-2.5mm}\includegraphics[width=0.33\columnwidth]{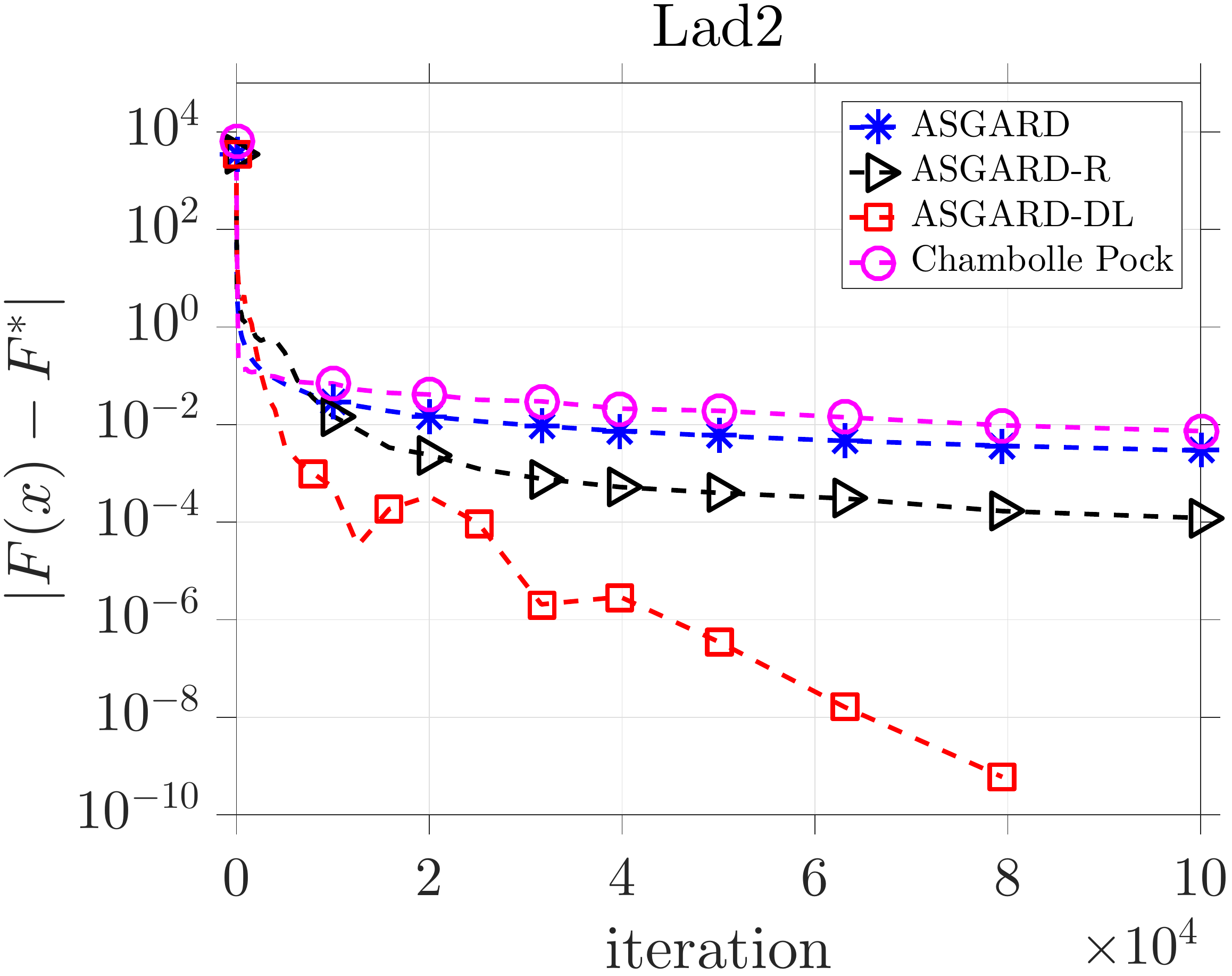} &
\hspace{-2.5mm}\includegraphics[width=0.33\columnwidth]{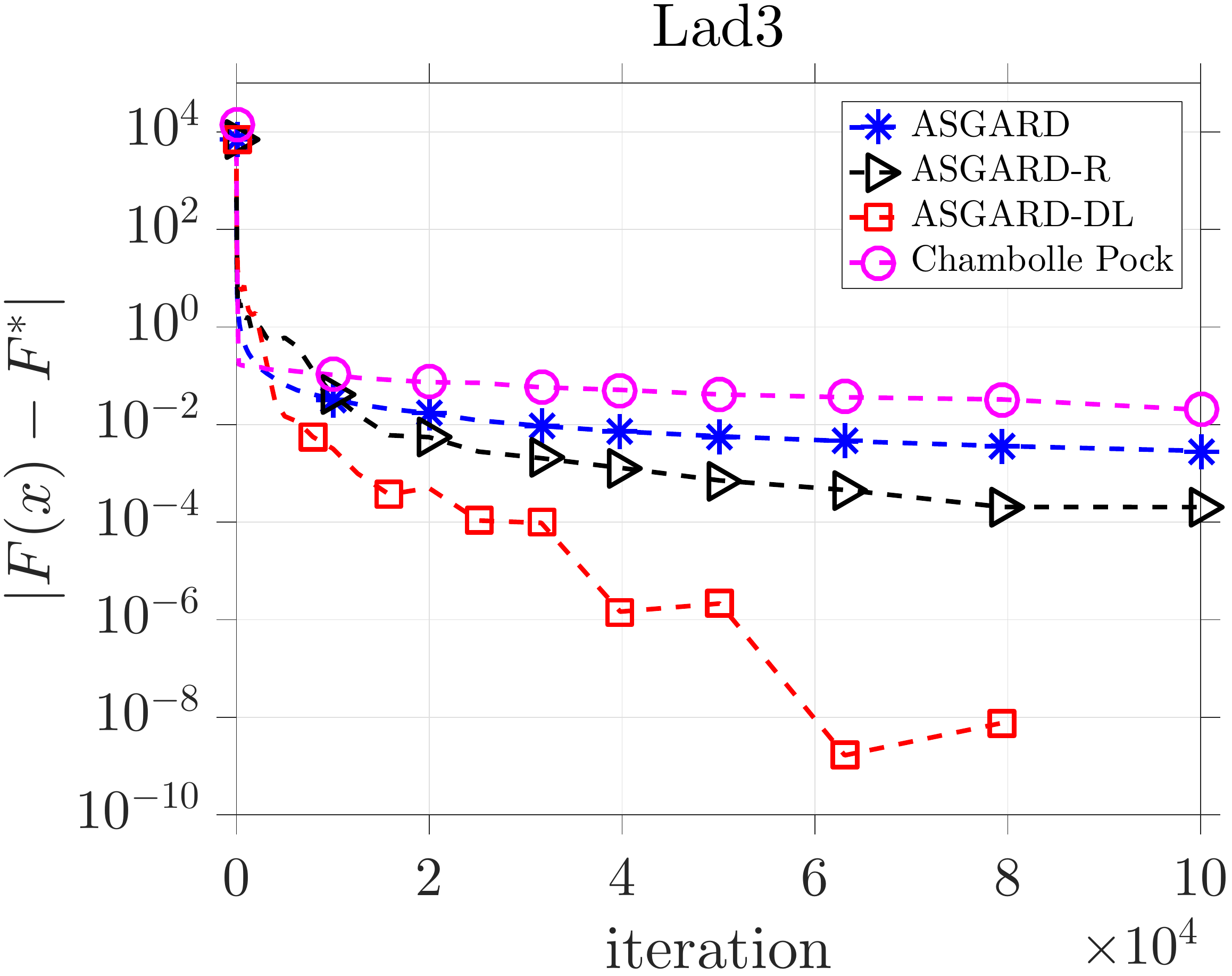} &\end{tabular}
\caption{Performance of 4 algorithms for LAD-Lasso problem in 3 different realization of varying problem size.}
\label{fig:fig_lad}
\vspace{-4ex}
\end{figure}

This problem fits to our template by setting $f(\cdot) := \lambda \norms{ \cdot }_1 $ and $g(\cdot) = \norms{ \cdot - b}_1$. 
We set $\beta_0$ in ASGARD, ASGARD-restart, and ASGARD-DL as $\beta_0 := 100\norms{ A}$.
We generated three problem instances of the size $n := 340r$, $p :=1000r$, $s :=100r$, where $s$ is the sparsity level, and $r = 1,2,3$ for the first, second and third instances, respectively. 
We present the results of this example in Figure~\ref{fig:fig_lad}. 

As we can see from Figure~\ref{fig:fig_lad} that, with the same per-iteration complexity, our method significantly outperforms the other algorithms after accuracy $10^{-4}$.
It beats other algorithms after a couple of hundred iterations and continues to decrease the objective values.
Although this problem is fully nonsmooth, heuristic restart such as in ASGARD still improves the performance of the non-restart one, but does not significant outperform.

\beforesubsec
\subsection{\textbf{Support Vector Machines}}
\aftersubsec
Our next example is the following primal support vector machines (SVM) problem in binary classification. 
Instead of classical models, we consider the following $\ell_1$-regularized nonsmooth hinge loss as proposed in~\cite{zhu20041}:
\begin{equation}\label{eq:PrimalSVM}
\min_{x\in\R^p}\Big\{ F(x) := \frac{1}{n} \sum_{i=1}^n \max\set{0, 1-b_i\iprods{a_i, x}} + \lambda\norms{x}_1 \Big\},
\end{equation} 
where $a_i \in \R^p$ are the feature vectors and $b_i \in \{ -1, +1 \} $ are the labels for $i=1,\cdots, n$.
We can cast \eqref{eq:PrimalSVM} into our template by setting $f(\cdot) := \lambda\norms{\cdot}_1$ and 
\begin{equation*}
g(Ax) = \frac{1}{n} \sum_{i=1}^n \max \{ 0, 1-b_i \langle a_i, x \rangle \}  
= \max_{y\in[0,1]^n}\iprods{y, Ax+\frac{1}{n}\mathbbm{1}},
\end{equation*}
where $A :=-\frac{1}{n} \begin{bmatrix}  b_1 a_1, b_2 a_2, \cdots, b_n a_n \end{bmatrix}^\top$ and $\mathbbm{1}$ is a vector of all ones. 
Clearly, the proximal operator of $g$ is simply a projection onto $[0,1]^n$.

\begin{figure}[hpt!]
\vspace{0ex}
\centering
\begin{tabular}{ccc}
\hspace{-2.5mm}\includegraphics[width=0.33\columnwidth]{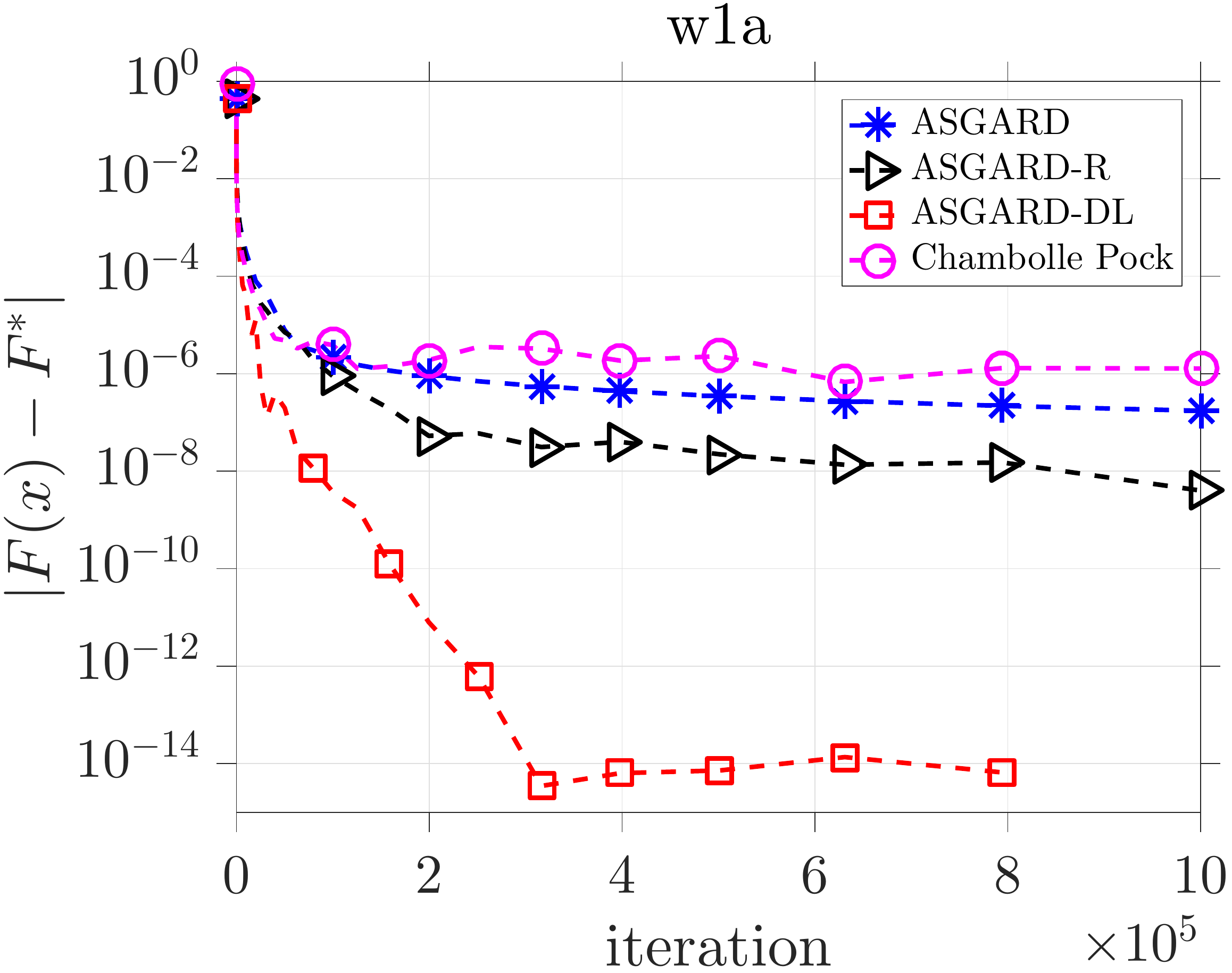} &
\hspace{-2.5mm}\includegraphics[width=0.33\columnwidth]{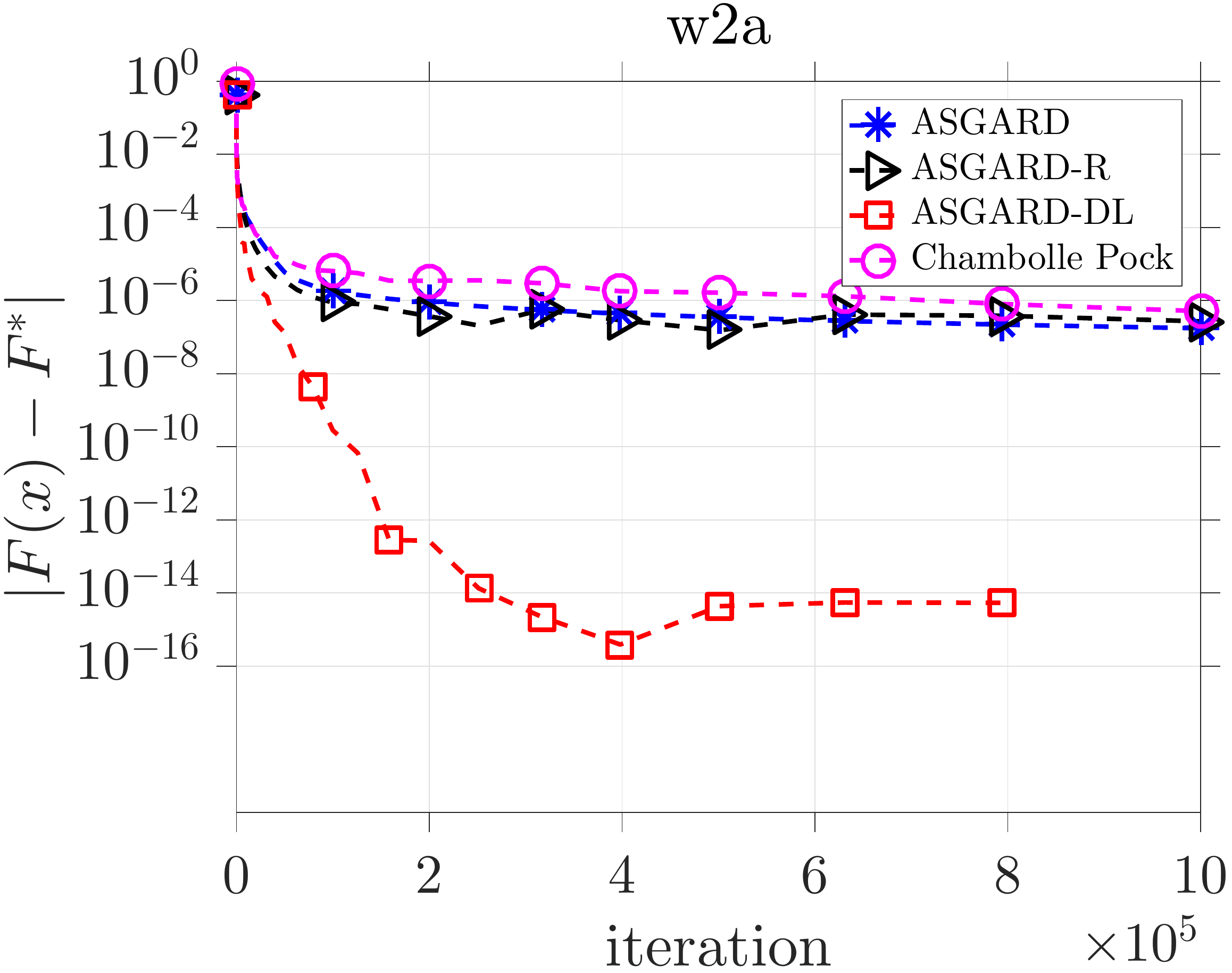} &
\hspace{-2.5mm}\includegraphics[width=0.33\columnwidth]{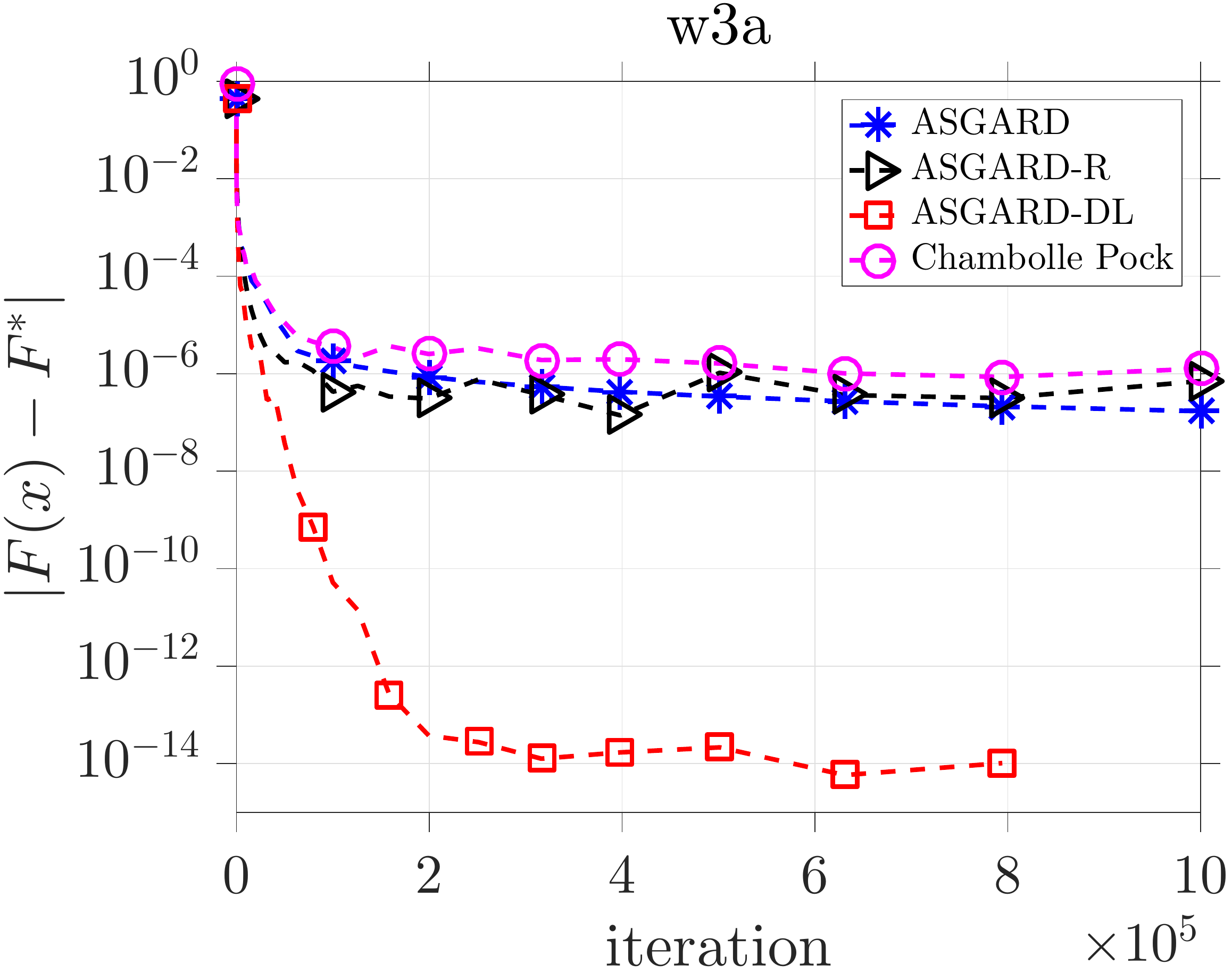} \\
\hspace{-2.5mm}\includegraphics[width=0.33\columnwidth]{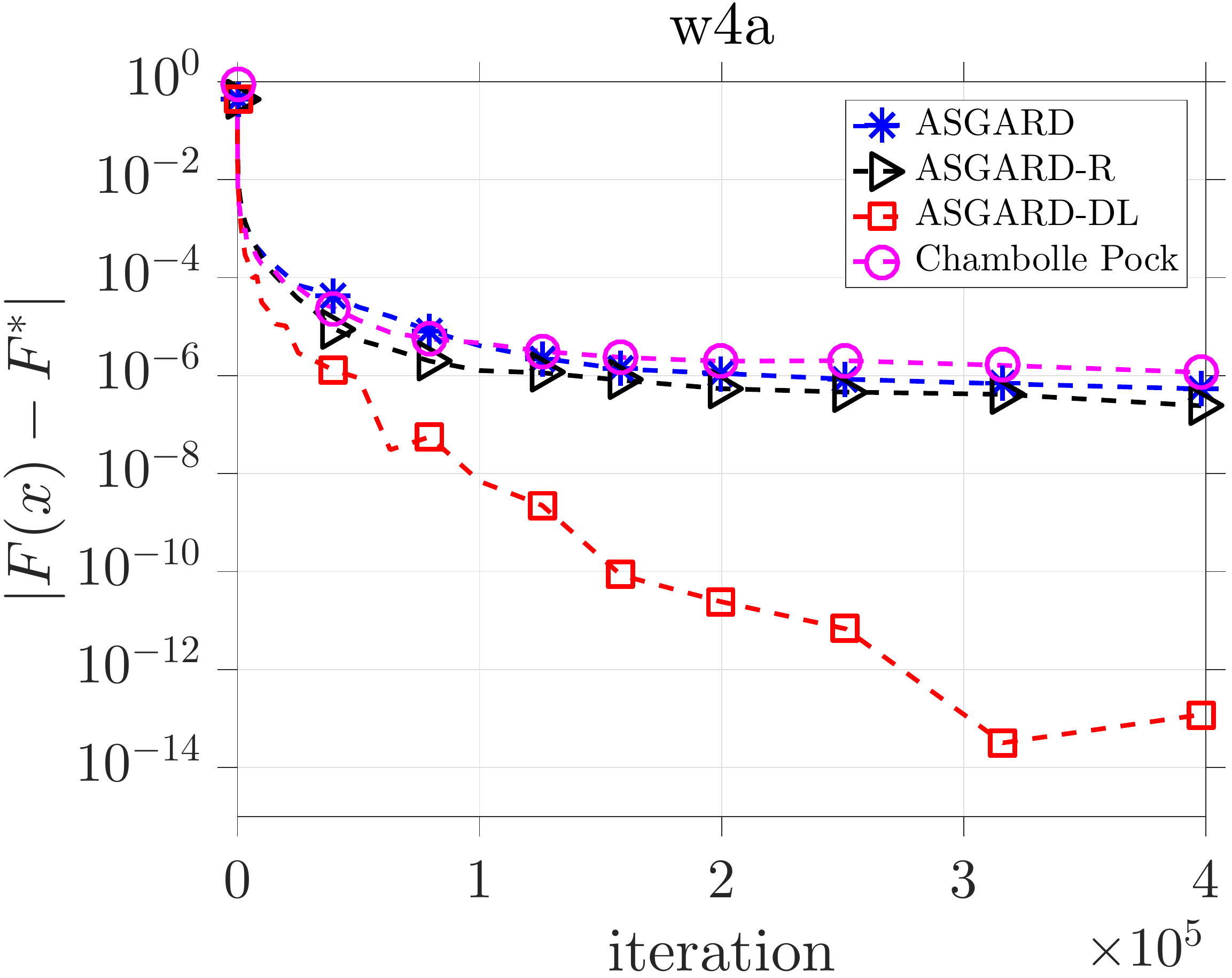} &
\hspace{-2.5mm}\includegraphics[width=0.33\columnwidth]{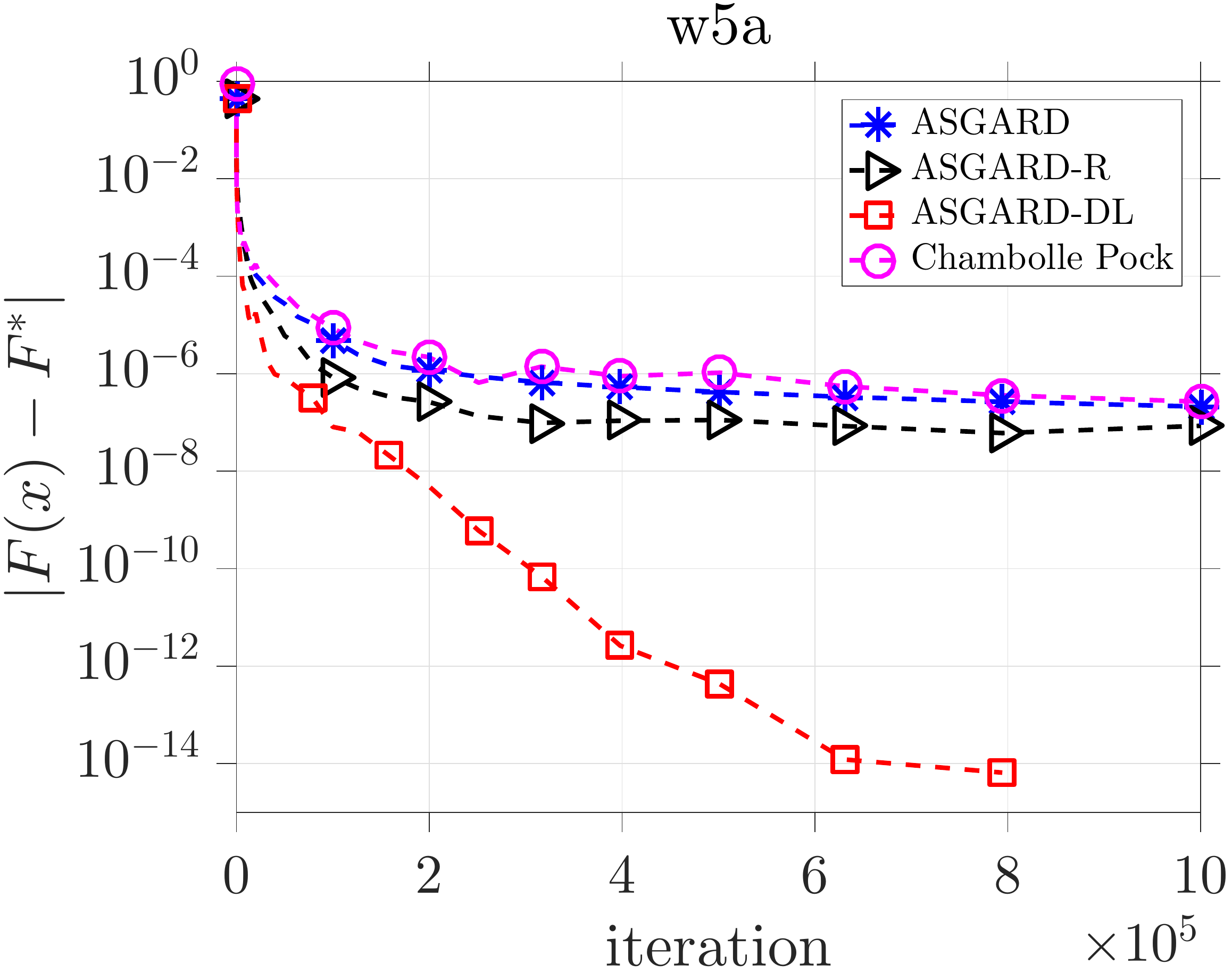} &
\hspace{-2.5mm}\includegraphics[width=0.33\columnwidth]{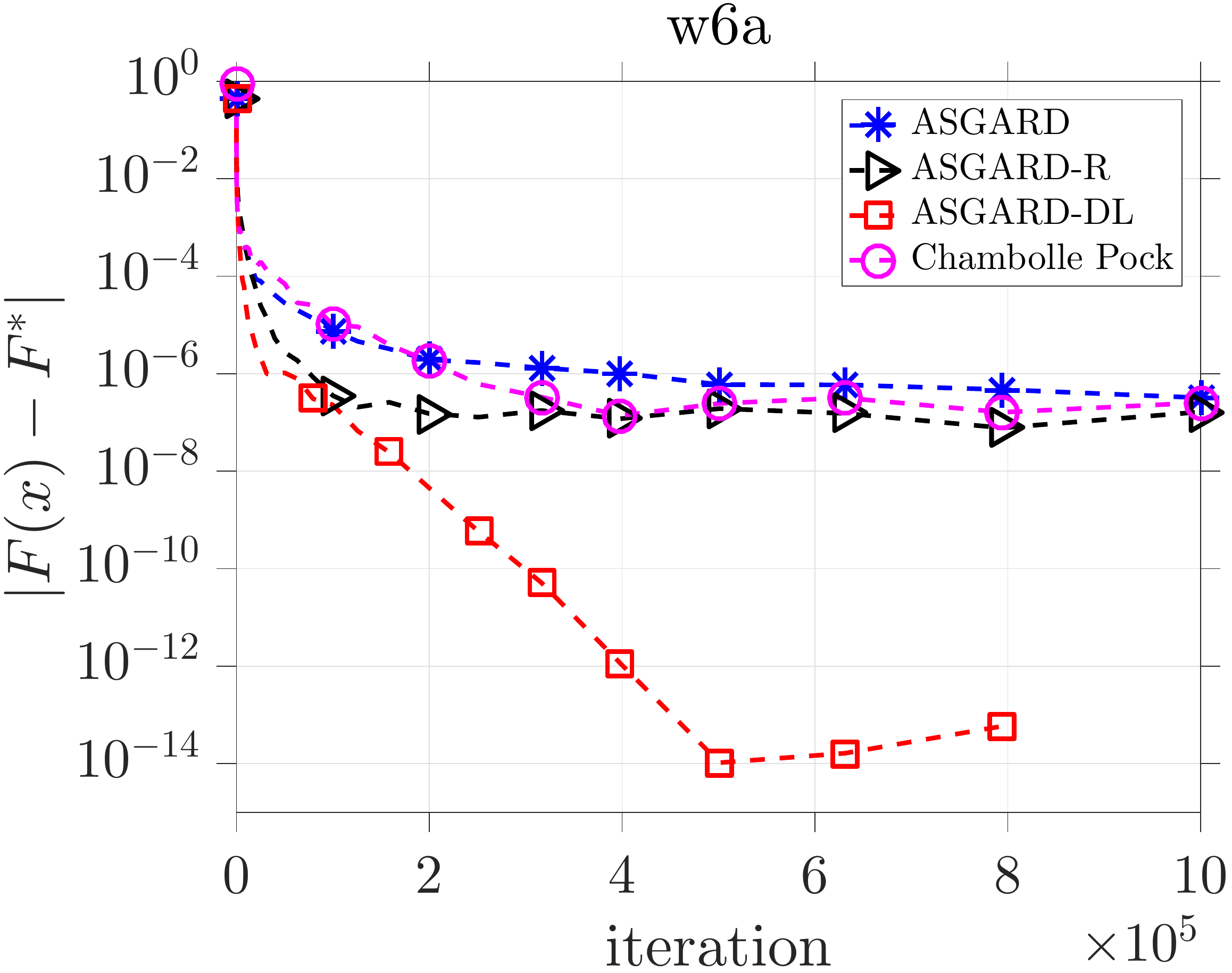} \\
\hspace{-2.5mm}\includegraphics[width=0.33\columnwidth]{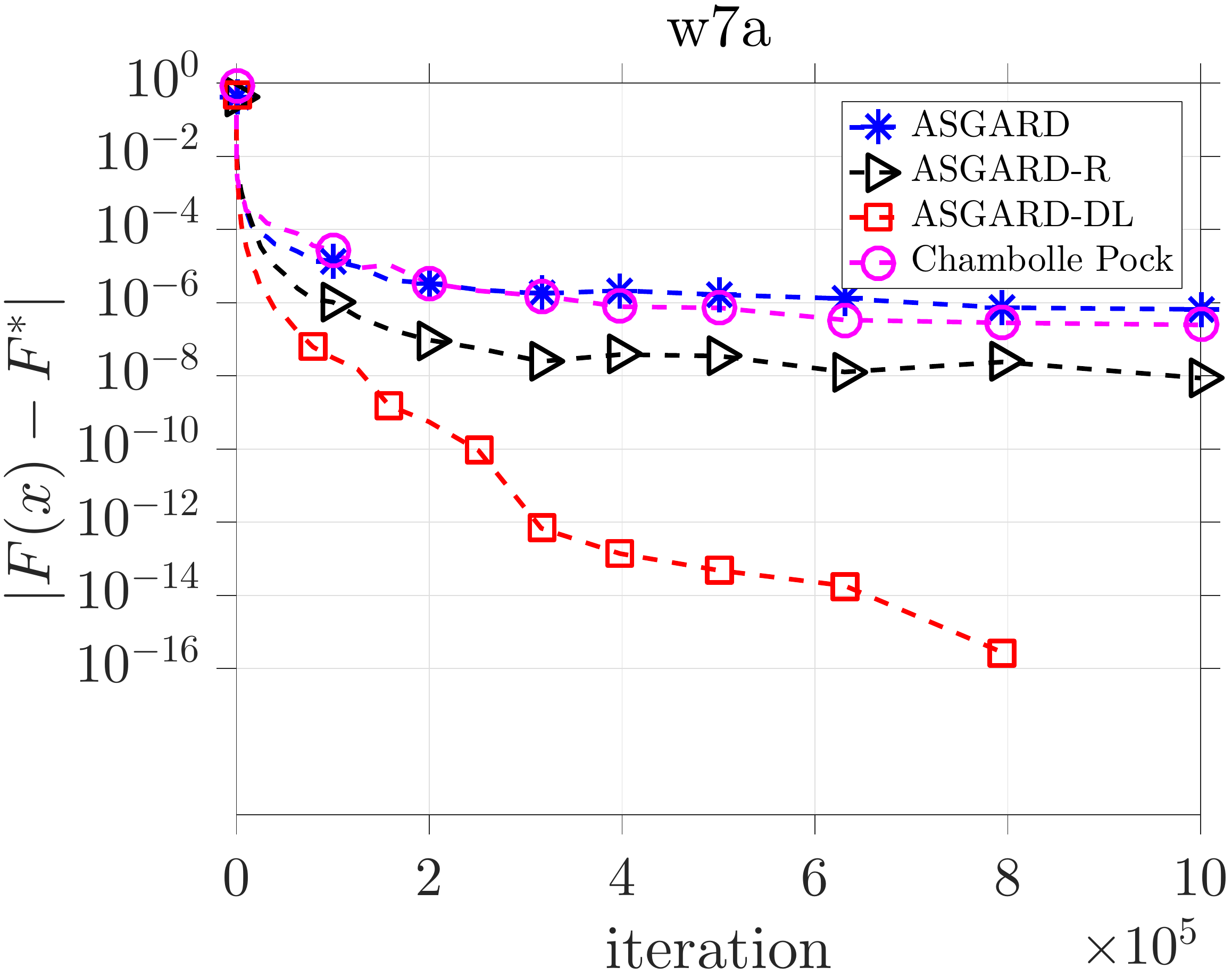} &
\hspace{-2.5mm}\includegraphics[width=0.33\columnwidth]{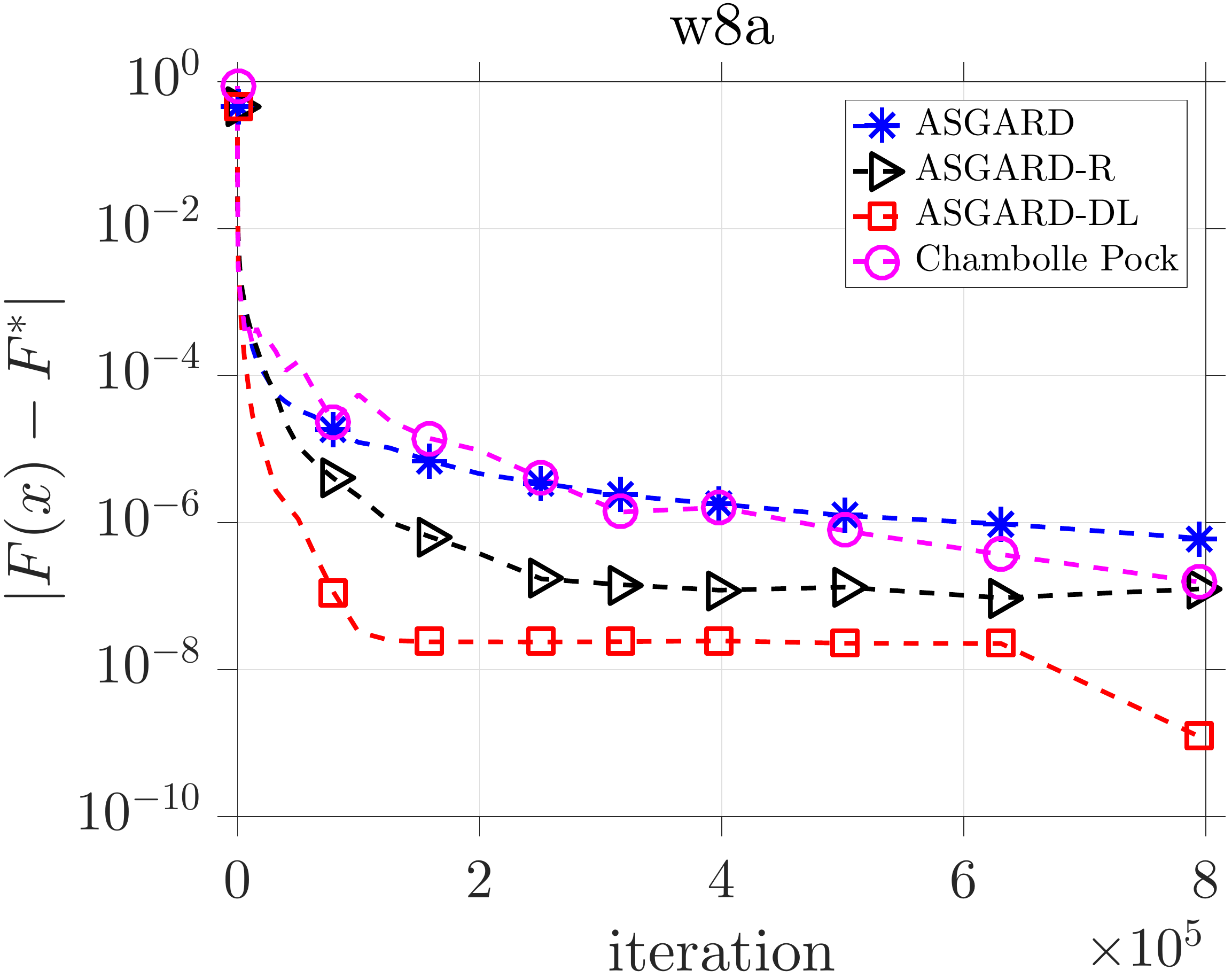} &
\end{tabular}
\caption{Performance of $4$ algorithms for the $\ell_1$-regularized SVM problem on $\set{\texttt{w1a}, \cdots, \texttt{w8a}}$.}
\label{fig:fig1_svm}
\vspace{-4ex}
\end{figure}

\begin{figure}[hpt!]
\centering
\begin{tabular}{ccc}
\hspace{-2.5mm}\includegraphics[width=0.475\columnwidth]{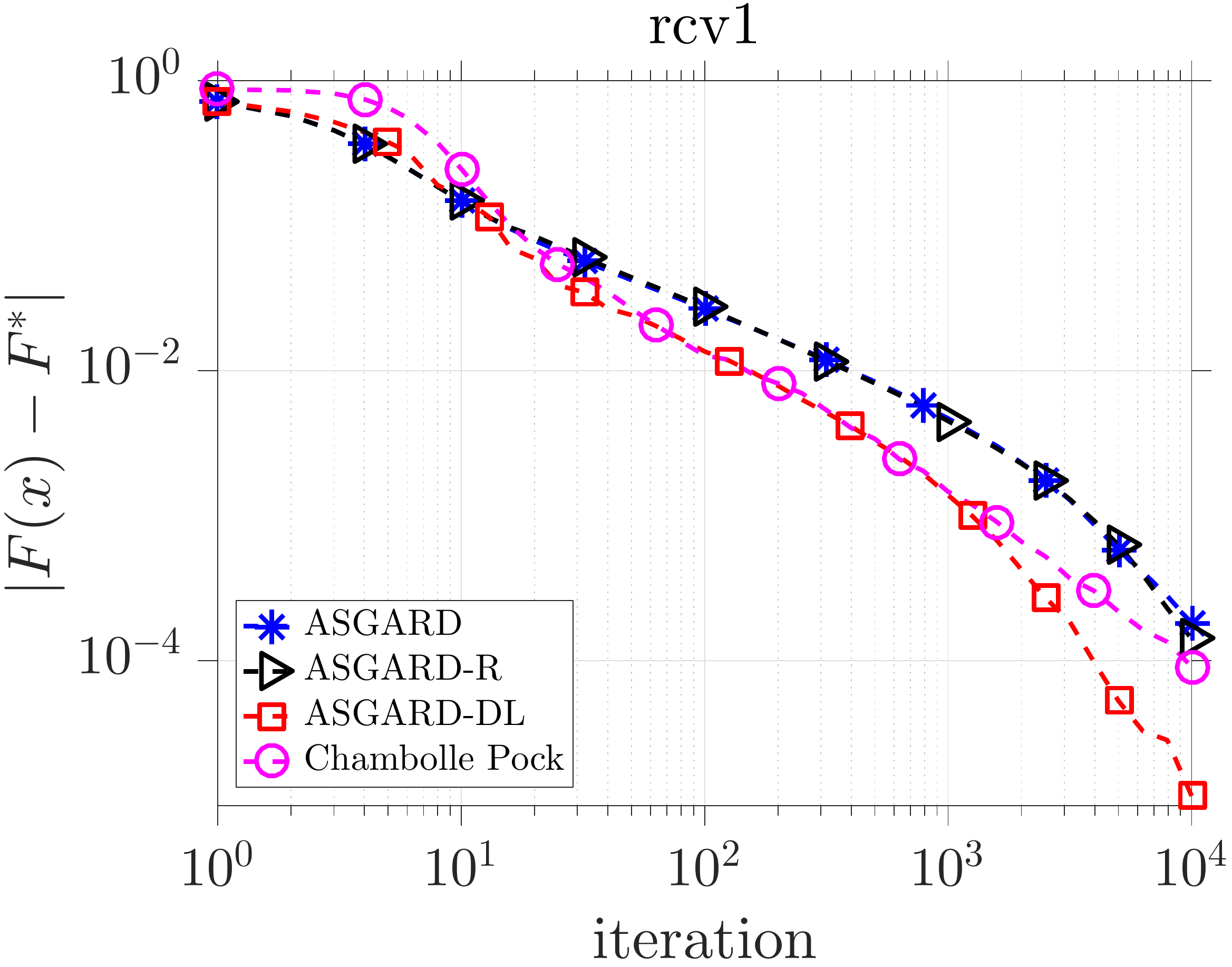} &
\hspace{-2.5mm}\includegraphics[width=0.475\columnwidth]{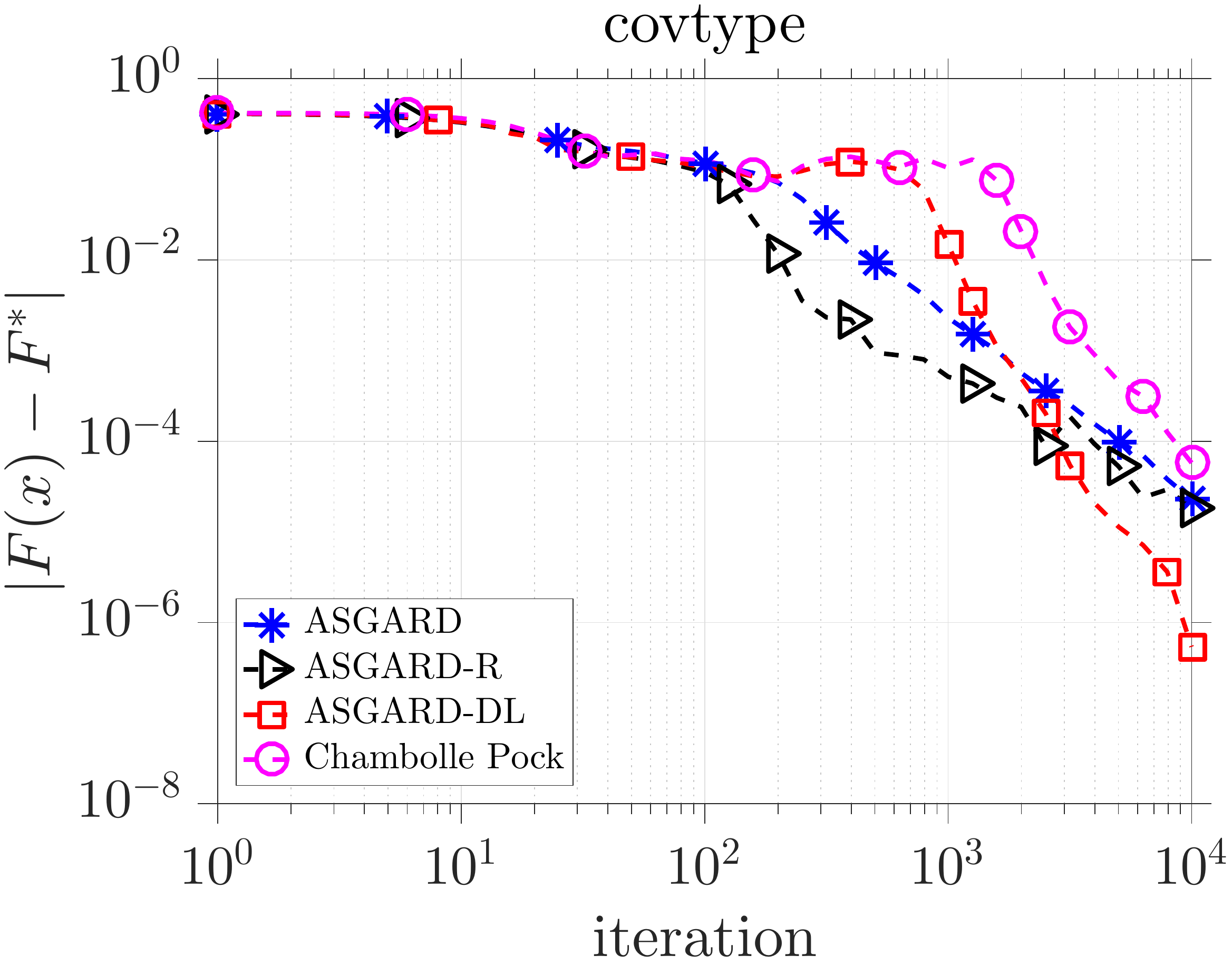} &
\end{tabular}
\caption{Performance of $4$ algorithms for the $\ell_1$-regularized SVM problem on $\set{\texttt{rcv1}, \texttt{covtype}}$.}
\label{fig:fig2_svm}
\vspace{-4ex}
\end{figure}

We use 10 different datasets from \texttt{libsvm}~\cite{CC01a} to test four different algorithms. 
The initial value $\beta_0$ in ASGARD, ASGARD-restart, and ASGARD-DL is set to $\beta_0  := 0.1\norms{A}$.
But for \texttt{covtype} dataset, we used $\beta_0 := 0.01 \norms{A}$.
The details about the datasets are given in Table~\ref{table:table_svm}.
\begin{table}[hpt!]
\begin{center}
\caption{Datasets used for classification.}\label{table:table_svm}
\begin{tabular}{l | r  | r}\toprule
    Data set & Training size & Number of features \\ \midrule
    w1a & 2,477 & 300 \\ 
    w2a & 3,470 & 300 \\  
    w3a & 4,912 & 300 \\  
    w4a & 7,366 & 300 \\  
    w5a & 9,888 & 300 \\  
    w6a & 17,188 & 300 \\ 
    w7a & 24,692 & 300 \\  
    w8a & 49,749 & 300 \\  
    rcv1 & 20,242 & 47,236 \\ 
    covtype & 581,012 & 54 \\    
\bottomrule
\end{tabular}
\vspace{-2ex}
\end{center}
\end{table}
We test $4$ algorithms on these ten datasets. 
The results of the first $8$ problems are given in Figure~\ref{fig:fig1_svm}, and the results of the two last problems are in Figure~\ref{fig:fig2_svm}. 

We again observe that Algorithm \ref{alg:A2} significantly outperform the other methods.
Since these algorithms have the same per-iteration complexity, it is sufficient to compare them in terms of iteration numbers.
Although all the algorithms have $\BigO{\frac{1}{k}}$-worst-case convergence rate, due to its double-loop, Algorithm \ref{alg:A2} performs much better than the others, especially for high accurate solutions.
This is not surprise.
The double-loop allows Algorithm \ref{alg:A1} to use large stepsize by frequently restarting $\tau_k$ and $\beta_k$, while ASGARD gradually decreases these parameters to zero, and Chambolle-Pock's method fixes the step-size.
Note that the $\BigO{\frac{1}{k}}$ rate of Chambolle-Pock's method is achieved via the averaging sequence, which is often much slower than the last iteration as we showed in Figures \ref{fig:fig1_svm} and \ref{fig:fig2_svm}. 

\beforesubsec
\subsection{\textbf{Markowitz Portfolio Optimization}}
\aftersubsec
We consider a classical example from Markowitz portfolio optimization~\cite{brodie2009sparse}.
The setting we consider here aims at maximizing the expected return for a given risk level. 
Assume that we are given a vector $\rho \in \R^n$, where $\rho$ is composed of expected returns from $n$ assets. 
This problem can be formulated as
\begin{equation}\label{eq:portfolio1}
\max_{x\in\R^{p}}\Big\{  \rho^{\top}x ~\mid~ x\in \triangle, ~\mathbb{E} \left[ \vert (a_i - \rho)^{\top}x  \vert^2 \right]  \leq \epsilon \Big\},
\end{equation} 
For our setting, we use empirical sample average instead of the expectation and convert the problem to a minimization problem by negating the objective:
\begin{equation}\label{eq:portfolio}
\min_{x\in\R^{ p}}\set{- \iprods{\rho, x} ~\mid~ x\in\triangle, ~\tfrac{1}{p} \Vert Ax \Vert_2 ^2 \leq \epsilon },
\end{equation} 
where $A = [(a_1 - \rho), (a_2 - \rho), \dots, (a_n - \rho)]^\top$. 
We map this problem to our template~\eqref{eq:three_comp} by mapping $f(\cdot) := \delta_{\triangle}(\cdot)$, $g(\cdot) := \delta_{\set{ \norm{\cdot}_2 \leq \sqrt{p\epsilon} }}(\cdot)$, and $h(x) := - \langle \rho, x \rangle$.
One key step of primal-dual algorithms is computing the projection onto an $\ell_2$-norm ball and on a simplex. 
Here, the complexity of simplex projection is $\mathcal{O}(p \log p)$.

As before, we apply $4$ algorithms to solve \eqref{eq:portfolio1}.
We use $4$ datasets that are also considered in~\cite{borodin2004can}. 
The details about the datasets are given in Table~\ref{table:table_portfolio}.
\begin{table}[hpt!]
\begin{center}
\vspace{-3ex}
\caption{Portfolio optimization datasets and parameters of algorithms.}\label{table:table_portfolio}
    \begin{tabular}{l r r r | r r r r r r}\toprule
    \multicolumn{4}{c|}{The size of datasets} &     \multicolumn{6}{c}{Parameters used in $4$ algorithms.} \\ \midrule
    Datasets & $n$ & $p$ & $\epsilon$ in~\eqref{eq:portfolio} & $\beta_0$ & RF & $\omega$ & $m_s$ & $\tau$ & $\sigma$ \\ 
    \midrule
    DJIA & 507 & 30 & 0.002 & $\norms{A}$ & 10 & 1.1 & 11 & $\frac{1}{\norms{A}}$ & $\frac{1}{\norms{A}}$ \\  
    NYSE & 5651 & 36 & 0.02 & $100\norms{A}$ & 10 & 1.1 & 11 & $\frac{1}{\norms{A}}$ & $\frac{1}{\norms{A}}$ \\ 
    SP500 & 1276 & 25 & 0.02 & $100\norms{A}$ & 10 & 1.2 & 6 & $\frac{1}{\norms{A}}$ & $\frac{1}{\norms{A}}$ \\ 
    TSE & 1258 & 88 & 0.002 & $100\norms{A}$ & 10 & 1.1 & 11 & $\frac{1}{\norms{A}}$ & $\frac{1}{\norms{A}}$ \\ 
\bottomrule
\end{tabular}
\end{center}
\vspace{-4ex}
\end{table}

We summarized the parameters that we used for these algorithms in Table~\ref{table:table_portfolio}, where $\beta_0$ is common to ASGARD, ASGARD-restart, our algorithm, restart frequency (RF) is specific to ASGARD-restart, $\omega$ and $m_s$ are specific to our algorithm and $\tau$ and $\sigma$ are specific to Chambolle-Pock's algorithm.

{\begin{figure}[hpt!]
\centering
\hspace{-2.5mm}\includegraphics[width=.95\columnwidth, height=.4\columnwidth]{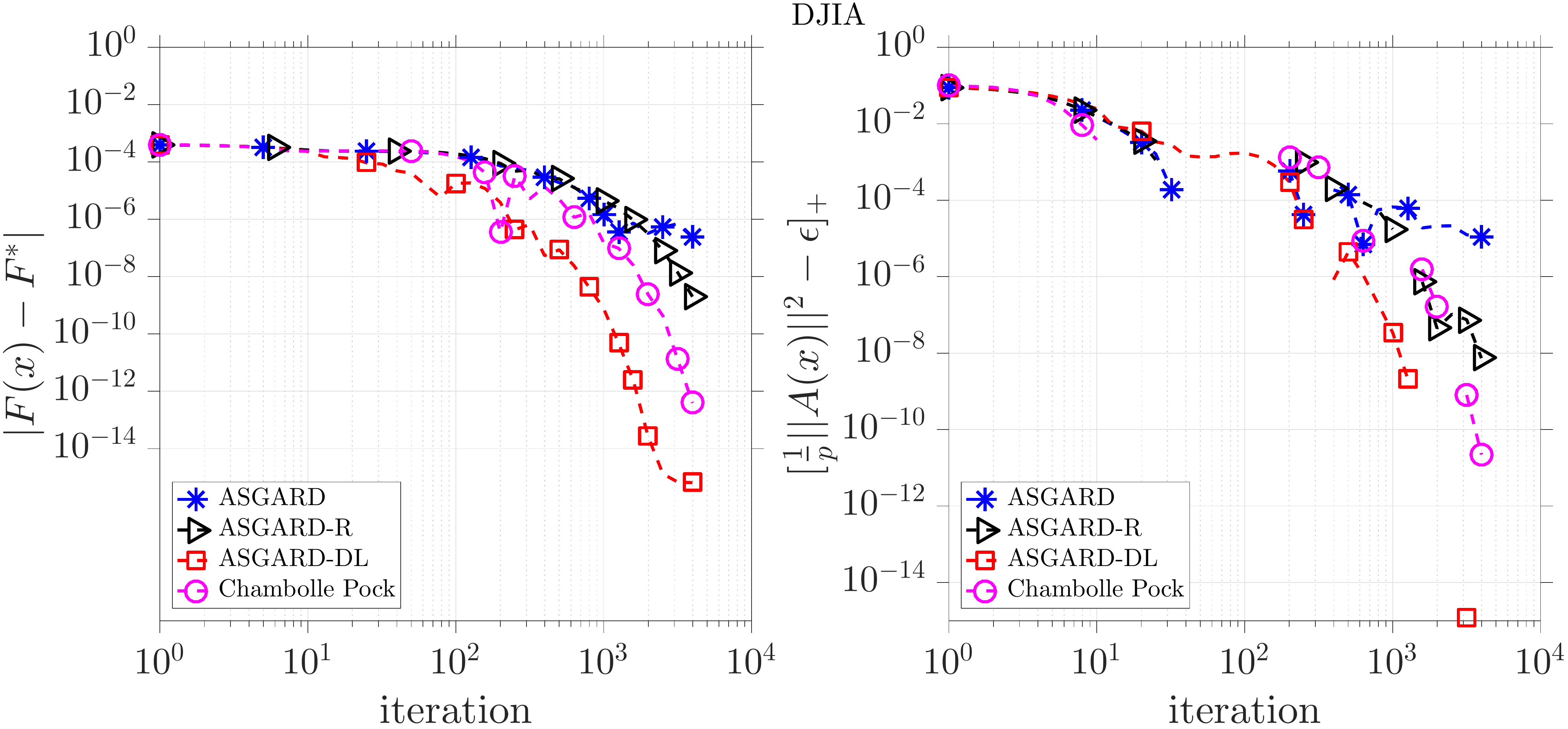} \\
\hspace{-2.5mm}\includegraphics[width=.95\columnwidth, height=.4\columnwidth]{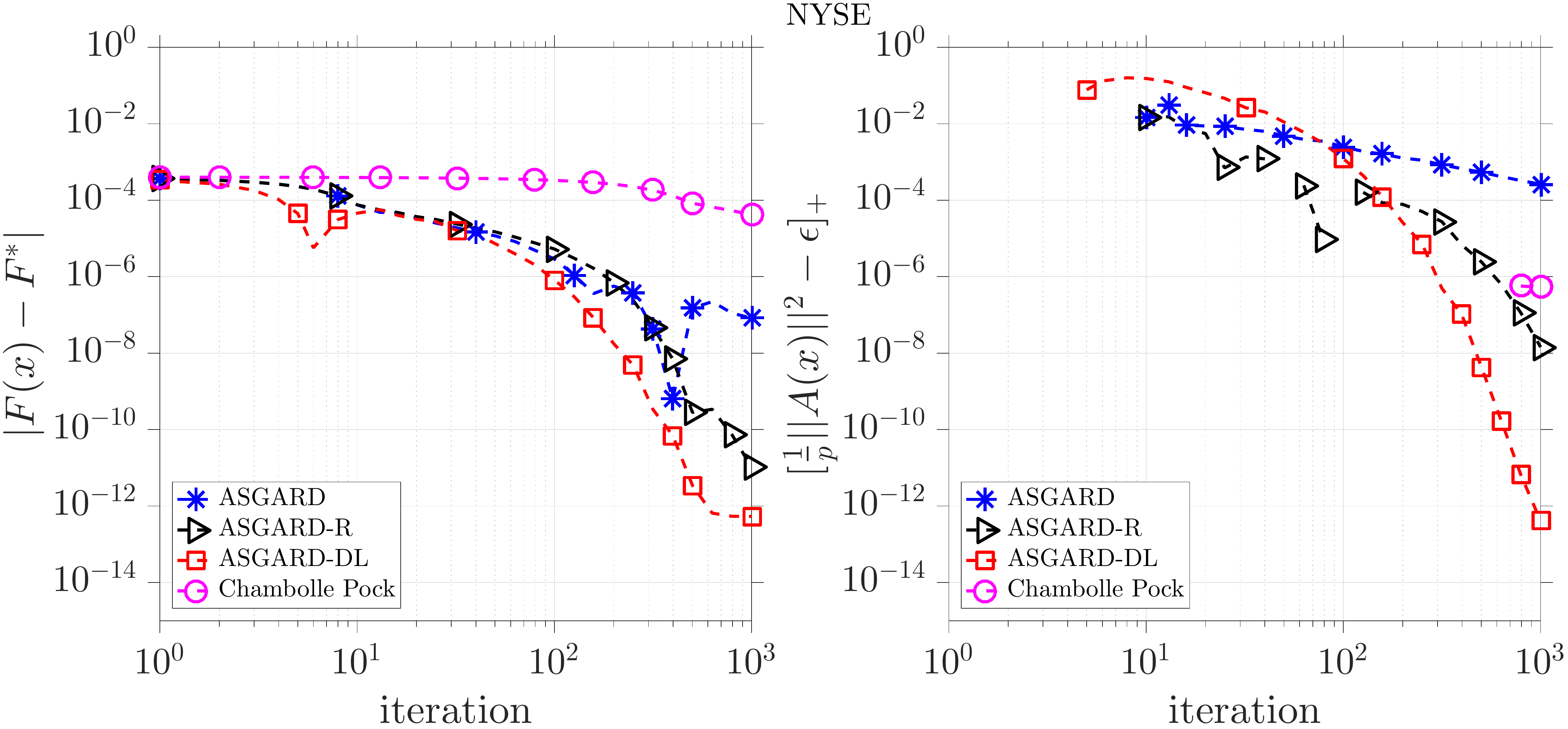} \\
\hspace{-2.5mm}\includegraphics[width=.95\columnwidth, height=.4\columnwidth]{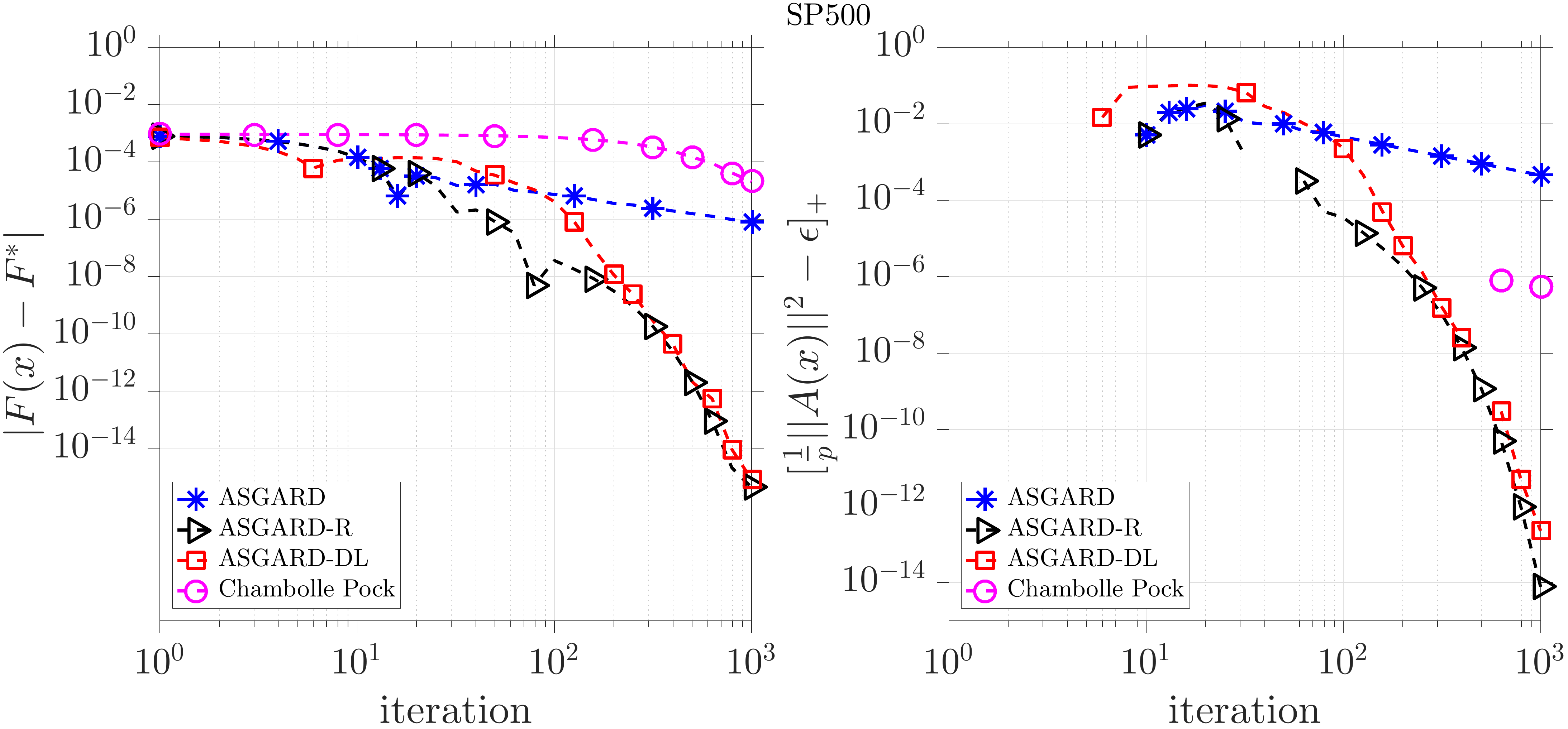} \\
\hspace{-2.5mm}\includegraphics[width=.95\columnwidth, height=.4\columnwidth]{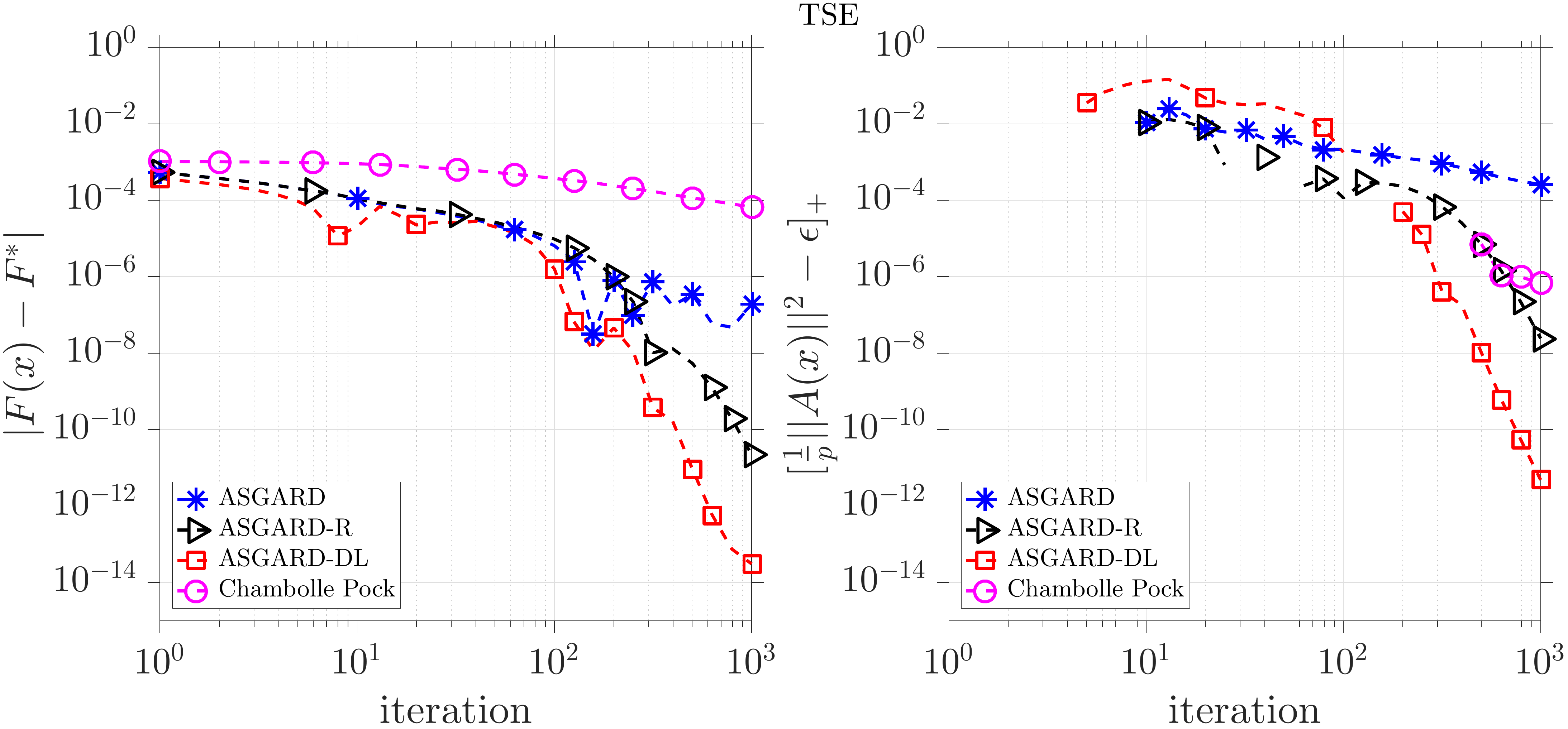}
\caption{Performance of 4 algorithms for Markowitz portfolio optimization problem on 4 real datasets.}
\label{fig:portfolio}
\vspace{-4ex}
\end{figure}}

We have tested 4 algorithms on 4 real datasets and the results are compiled in Figure~\ref{fig:portfolio}. 
As can be seen, except for SP500 dataset, Algorithm \ref{alg:A2} significantly outperforms the other methods and shows a much faster practical performance than $\BigO{\frac{1}{k}}$ guarantee.
For SP500 dataset, ASGARD-restart algorithm shows a comparable performance to our method.
However, as discussed in~\cite{TranDinh2015b}, the effect of restarting to ASGARD method is not understood theoretically.
Our algorithm theoretically preserves the best-known $\BigO{\frac{1}{k}}$ guarantee while performing as fast as, and most of the times faster than the heuristic restarting ASGARD method.

\beforesubsec
\subsection{\textbf{Sparse Subspace Clustering}}
\aftersubsec
In the last example, we consider the following sparse subspace clustering problem which has broad applications in machine learning, computer vision and image processing. 
This problem is studied extensively in the literature~\cite{elhamifar2013sparse,elhamifar2009sparse,pourkamali2018efficient}. 
In this problem setting, we assume that there exist $n$ points $\set{ x_1, x_2, \cdots, x_n } \in \R^p$ lying in the union of subspaces in $\R^p$. 
We form a matrix $X\in\R^{p \times n}$ by stacking $\set{ x_1, x_2, \cdots, x_n }$ as the columns. 
With this notation, each point can be represented as
\begin{equation*}
x_j = Xc_j + e_j, ~~\text{s.t.}~~ [c_j]_j = 0 ~~\text{and}~~ \mathbbm{1}^{\top}c_j = 1.
\end{equation*}
where $c_j\in\R^n$ represents the coefficients to represent point $x_j\in\R^p$ as an affine combination of other points, $e_j\in\R^p$ is the representation error and $\mathbbm{1}\in\R^n$ is a vector of $1$'s.

This formulation can be represented compactly by stacking $c_j$ to the $j^{\text{th}}$ column of matrix $C$ as follows:
\begin{equation}
X = CX ~~ \text{s.t.}~~ \mathrm{diag}(C)=0, ~C^{\top} \mathbbm{1} = \mathbbm{1}.
\end{equation}
The optimization problem that we will tackle in this subsection is referred to as an SSC-Lasso problem in the literature, and is written as
\begin{equation}\label{prob:ssc}
\min_{C\in\R^{n\times n}}\Big\{ \norm{C}_1 + \tfrac{\lambda}{2} \norm{X-CX}^2~\mid~\mathrm{diag}(C)=0, ~C^{\top}\mathbbm{1} = \mathbbm{1} \Big\}.
\end{equation}
In~\cite{elhamifar2009sparse} and~\cite{elhamifar2013sparse}, ADMM is used to solve~\eqref{prob:ssc} and recently, \cite{pourkamali2018efficient} proposed an efficient implementation of ADMM and an application of standard accelerated proximal scheme to this setting. 
One drawback of applying accelerated proximal schemes to~\eqref{prob:ssc} is the evaluation of the proximal operator of an $\ell_1$-norm over the linear constraint $C^\top \mathbbm{1} = \mathbbm{1}$.
This requires additional computation cost of $\log(n)pn^2$. 
We fit \eqref{prob:ssc} into our template \eqref{eq:three_comp} by defining $f(\cdot) := \norms{ \cdot }_1 + \delta_{\set{\mathrm{diag}(\cdot)=0}}(\cdot)$, $g(\cdot) := \delta_{\set{\langle \cdot, \mathbbm{1} \rangle=\mathbbm{1}}}(\cdot)$, and $h(\cdot) = \frac{\lambda}{2} \norms{X - X(\cdot)}^2$. 
If we apply Algorithm \ref{alg:A2} to solve this reformulation, then no  extra computation cost is incurred as in accelerated proximal gradient methods.

We use a classic benchmark Extended Yale B dataset~\cite{georghiades2001few} to test the sparse subspace clustering problem \eqref{prob:ssc}. 
This dataset contains face pictures of $38$ individuals taken under $64$ different environmental conditions. 
As previous works, we use downsampled images of size $48\times 42$ pixels which correspond to $p=2016$. 
We ran experiments with ADMM, TFOCS, and our method ASGARD-DL. 
We note that our method includes tuning parameters similar to ADMM. 
We use  $\beta_0 := \sqrt{\norms{M}}$, where $M(C) = C^{\top}\mathbbm{1}$.
We randomly selected $m=2, 3, 5$ clusters and ran $3$ trials for each case. 
We have used the implementation of ADMM~\cite{elhamifar2013sparse,elhamifar2009sparse} and TFOCS~\cite{pourkamali2018efficient} provided by the authors of~\cite{pourkamali2018efficient,elhamifar2013sparse,elhamifar2009sparse}. 
For fair comparison, we ran the algorithms for the same duration of time and reported the results accordingly. 
We used objective value and clustering error as comparison measures as~\cite{pourkamali2018efficient}.

\begin{table}[hpt!]
\vspace{-3ex}
\begin{center}
\caption{Comparison of 3 methods on the SSC-Lasso problem with $m=2, 3, 5$ clusters and $3$ independent trials of each.}\label{table:table_ssc}                                    
\begin{tabular}{ c | r | r | r } \toprule
Problem & ADMM & TFOCS & ASGARD-DL \\                                
\midrule
($n=2$)-objective-trial 1 & 236.5653 & 226.4371 & 225.7578 \\   
($n=2$)-Clustering error-trial 1 & 0.0312 & 0.0391 & 0.0391 \\  
\midrule                                                    
($n=2$)-objective-trial 2 & 200.3710 & 192.2985 & 191.5177 \\   
($n=2$)-Clustering error-trial 2 & 0.0234 & 0.0469 & 0.0469 \\  
\midrule                                                    
($n=2$)-objective-trial 3 & 197.2510 & 188.7655 & 188.1555 \\   
($n=2$)-Clustering error-trial 3 & 0.0703 & 0.0938 & 0.0938 \\  
\midrule                                                    
($n=3$)-objective-trial 1 & 329.9188 & 320.9690 & 319.5887 \\   
($n=3$)-Clustering error-trial 1 & 0.0156 & 0.0156 & 0.0312 \\  
\midrule                                                    
($n=3$)-objective-trial 2 & 341.1980 & 330.6395 & 329.5704 \\   
($n=3$)-Clustering error-trial 2 & 0.0729 & 0.0677 & 0.0677 \\  
\midrule                                                    
($n=3$)-objective-trial 3 & 398.8778 & 389.3963 & 388.0739 \\   
($n=3$)-Clustering error-trial 3 & 0.4375 & 0.3594 & 0.3646 \\  
\midrule                                                    
($n=5$)-objective-trial 1 & 549.8250 & 530.0340 & 526.1905 \\   
($n=5$)-Clustering error-trial 1 & 0.1625 & 0.1156 & 0.0906 \\  
\midrule                                                    
($n=5$)-objective-trial 2 & 482.8483 & 467.0535 & 461.6563 \\   
($n=5$)-Clustering error-trial 2 & 0.2188 & 0.1125 & 0.1562 \\  
\midrule                                                    
($n=5$)-objective-trial 3 & 1029.5459 & 1017.7089 & 1025.6752 \\
($n=5$)-Clustering error-trial 3 & 0.3156 & 0.3469 & 0.3156 \\  
\bottomrule                                                    
\end{tabular}                                                 
\end{center}
\vspace{-5ex}
\end{table}                      

We can see  from Table~\ref{table:table_ssc} that our method consistently outperforms other methods in terms of objective values, and has similar performance in terms of the clustering error.
We present our algorithm as another candidate for solving the classical sparse subspace clustering problem with a lower per iteration cost than previous approaches ADMM and TFOCS and similar performance.

\beforesec
\section{Further discussion and comparison with previous work}\label{sec:prev-work}
\aftersec
Theory and numerical methods for solving \eqref{eq:ns_cvx_prob} and \eqref{eq:constr_cvx} are well-studied in the literature.
Due to such a large proportion of solution methods, we only focus on some recent works that are the most related to our method developed in this paper.
We briefly survey these results to highlight the similarities and differences with our work.

In~\cite{Nesterov2005c}, Nesterov proposed combining smoothing technique and accelerated gradient methods to obtain $\BigO{\frac{1}{\varepsilon}}$-iteration complexity to obtain an $\varepsilon$-approximate solution to~\eqref{eq:ns_cvx_prob}. 
However, this method requires $\varepsilon$ to be predefined, and both primal and dual domains are bounded.
In addition, the step-size of the underlying gradient-type scheme is proportional to $\varepsilon$, which is often small. 
This leads to a poor performance in early iterations.
In  \cite{Nesterov2005d}, Nesterov introduced an excessive gap technique to develop new algorithms that allow the smoothness parameter to be adaptively updated.
Nevertheless, these methods still require both primal and dual domains to be bounded, and one additional proximal operator for every two iterations.

In ~\cite{Chambolle2011,chambolle2016ergodic}, A. Chambolle and T. Pock proposed a primal-dual algorithm to solve~\eqref{eq:ns_cvx_prob} that achieves $\BigO{\frac{1}{k}}$-convergence rate. 
This rate is guaranteed on a gap function and also requires both primal and dual domains to be bounded, which is unfortunately not applicable to \eqref{eq:constr_cvx}.
In addition, the guarantee of their methods relies on ergodic or weighted averaging sequences. 
Note that, in sparse and low-rank optimization and image processing, taking averaging sequence unfortunately destroys desired structures of approximate solutions.
In addition, as also presented with numerical evidence, averaging sequences perform poorly in practice.

In \cite{xu2016homotopy}, the authors proposed a homotopy algorithm called Homotopy Smoothing algorithm (HOPS) which also essentially relies on Nesterov's smoothing technique \cite{Nesterov2005c}.
HOPS employs a similar strategy to ours in the sense of having a double loop structure. 
However, this method suffers from several drawbacks. 
First, it only applies to unconstrained problems as in \eqref{eq:ns_cvx_prob}, but not to \eqref{eq:constr_cvx} due to the unboundedness of  the dual domain.
Second, it requires knowing $\varepsilon_0 = P(x^0) - P(x^\star)$ to be able to set the initial smoothness parameter.
Third, HOPS requires tuning the number of inner iterations and the rate at which the smoothness parameter is going to be reduced.
The alternative of HOPS to alleviate this issue requires a bounded primal domain which further restricts the usage of their method.

\begin{sidewaystable}
    \centering
    \caption{Summary of algorithms that require two proximal operators each iteration.
    Note that $z^k := \sum_{k=1}^K w^k x^k$ where K is the maximum number of iterations and $w^k$ are the weights.
    For solving the unconstrained problem with ALM/ADMM methods, we split the problem.
    Practicality column refers to whether using the iterate in the convergence rate gives a fast practical performance or not. }
   \begin{tabular}{| l | c | c | c | c | c |}
   \hline
    Algorithm & $g$ is Lipschitz & $g=\delta_{\{b\}}$ & Type of rate & Set $\epsilon$ & Practicality \\
   \hline
    Nesterov~\cite{Nesterov2005c} & $P(x^k) - P^\star \leq \mathcal{O}\left( \max \left( \epsilon, \frac{1}{\epsilon k^2} \right) \right) $ & Not applicable & \textbf{Non-ergodic} & Yes & No \\
    Chambolle-Pock~\cite{Chambolle2011} & $G(z^k) \leq \mathcal{O} \left( \frac{1}{k} \right)$ & Convergence & Ergodic & \textbf{No} & No\\
    Linearized ALM~\cite{xu2017accelerated} & $ P(z^k) - P^\star \leq \mathcal{O} \left( \frac{1}{k} \right)$  & 
    \begin{minipage}{0.1\textwidth}
    \begin{align*}
    \vert f(z^k) - f^\star \vert \leq \mathcal{O} \left( \frac{1}{k}  \right)\\
    \Vert Ax-b \Vert \leq \mathcal{O} \left( \frac{1}{k}  \right)
    \end{align*}
    \end{minipage}
     & Ergodic & \textbf{No} & No \\
    Inexact ALM~\cite{xu2017iteration} & $ P(z^k) - P^\star \leq \mathcal{O} \left( \max\left(\epsilon_k, \beta_k \right)  \right)$  & 
    \begin{minipage}{0.1\textwidth}
    \begin{align*}
    \vert f(z^k) - f^\star \vert \leq \mathcal{O} \left( \max\left( \epsilon_k + \beta_k \right)  \right)\\
    \Vert Ax-b \Vert \leq \mathcal{O} \left( \beta_k \right) ~~~~~~~~~~
    \end{align*}
    \end{minipage} & \textbf{Non-ergodic} & Yes & No\\
    Linearized ADMM~\cite{xu2017accelerated} &$ P(z_1^k, z_2^k) - P^\star = \mathcal{O} \left( \frac{1}{k} \right)$  & 
    \begin{minipage}{0.1\textwidth}
    \begin{align*}
    \vert f_1(z_1^k)+f_2(z_2^k) - f_1^\star - f_2^\star \vert \leq \mathcal{O} \left(  \frac{1}{k} \right)\\
    \Vert A_1z_1^k + A_2z_2^k-b \Vert  \leq \mathcal{O} \left( \frac{1}{k} \right) 
    \end{align*}
    \end{minipage}& Ergodic & \textbf{No} & No\\
    ASGARD~\cite{TranDinh2015b} &$P(x^k) - P^\star \leq \mathcal{O} \left( \frac{1}{k} \right)$  & 
    \begin{minipage}{0.1\textwidth}
    \begin{align*}
    \vert f(x^k) - f^\star \vert \leq \mathcal{O} \left( \frac{1}{k}  \right)\\
    \Vert Ax^k-b \Vert \leq \mathcal{O} \left( \frac{1}{k}  \right)
    \end{align*}
    \end{minipage} & \textbf{Non-ergodic} & \textbf{No}& No\\
    This paper (Algorithm~\ref{alg:A2}) & $P(x^k) - P^\star \leq \mathcal{O} \left( \frac{1}{k} \right)$  & 
    \begin{minipage}{0.1\textwidth}
    \begin{align*}
    \vert f(x^k) - f^\star \vert \leq \mathcal{O} \left( \frac{1}{k}  \right)\\
    \Vert Ax^k-b \Vert \leq \mathcal{O} \left( \frac{1}{k}  \right)
    \end{align*}
    \end{minipage} & \textbf{Non-ergodic} & \textbf{No} & \textbf{Yes}\\
    \hline
    \end{tabular}
\end{sidewaystable}

For constrained problem \eqref{eq:constr_cvx},
among different methods, augmented Lagrangian (ALM), alternating direction method of multipliers (ADMM), alternating minimization algorithms (AMA), and penalty methods are the most popular.
Inexact augmented Lagrangian methods (iALM)~\cite{lan2016iteration,nedelcu2014computational,xu2017iteration} relies on a double loop structure similar to our method.
However, termination rules for these methods require the desired accuracy $\varepsilon$ to be set a priori.
In addition, in practice, it is not easy to check when the inner problem is solved to an $\varepsilon_k$-accuracy in the $k$-th iteration. 
Such an estimate is often derived from the worst-case complexity bound of the underlying solution method, and therefore, the corresponding algorithm is not efficient in practice.

While ADMM works really well and is widely used in practice, AMA is rarely used and requires additional conditions to converge. 
The best-known convergence rate of ADMM and its variants such as linearized ADMM and preconditioned ADMM is $\BigO{\frac{1}{k}}$ under standard assumptions \cite{He2012,Monteiro2012b,Monteiro2010,ouyang2015accelerated,xu2017accelerated}.
Moreover, this rate is given in an ergodic sense, and examples show that such a rate is optimal.
See \cite{Boyd2011} for more information about the behavior of ADMM.
In practice, however, the ergodic rate is rather pessimistic, which is much slower than the last iterate sequence (see Subsection \ref{subsec:example1} as an example).
So far, we are not aware of any work showing an $\BigO{\frac{1}{k}}$-rate of the standard ADMM or its linearized and preconditioned ADMM in the last iterate.
A recent work \cite{li2016accelerated} combined preconditioned/linearized ADMM and Nesterov's accelerated schemes to achieve an $\BigO{\frac{1}{k}}$-non-ergodic convergence rate.

Penalty methods use a quadratic penalty term to move the constraints to the objective and solve the subproblems by changing the penalty parameter~\cite{lan2013iteration,necoara2015complexity}.
Similar to iALM, these methods also do not have clear implementable termination rules for the inner loop.
In addition, they do not involve dual variables.
Therefore, they are often less competitive with primal-dual methods.
A recent work \cite{tran2017proximal} proposed a new alternating quadratic penalty algorithm to solve \eqref{eq:constr_cvx} that has the same $\BigO{\frac{1}{k}}$-non-ergodic convergence rate as in this paper.
Nevertheless, this method is completely different from this paper and does not have an update on the dual center.

Compared to our previous work  \cite{TranDinh2015b}, ASGARD, our new algorithm shares some similarities but also has several differences.
First, it has inner and outer loops but the guarantee is on the overall iterations.
Second, it works with any Bregman divergence induced by a general prox-function when solving \eqref{eq:ns_cvx_prob}, while ASGARD only works with the Bregman distances induced by a strongly convex and Lipschitz gradient prox-function.
This excludes some important Bregman divergences such as the Kullback-Leibler (KL)  divergence.
Third, our algorithm allows us to use different norms while computing proximal operators, compared to ASGARD which works with only Euclidean norms.
Fourth, it automatically restarts both the primal and dual variables as well as the parameters.
It also has a rigorous convergence guarantee, while the practical restarting variant of  ASGARD does not have convergence guarantee.

We developed a novel analysis for our double loop structured smoothing algorithm which allowed us to derive flexible rules for parameters in both unconstrained and constrained problems, in contrast to~\cite{xu2016homotopy}.
Our analysis gives insights on the heuristic restarting strategies in~\cite{TranDinh2015b} as well as on the number of inner iterations in the algorithm.
It also gives explicit number of iterations for the inner subproblems and does not require to predefine the horizon as opposed to iALM.
Table \ref{tbl:comparison} summarizes the key differences between different methods we have discussed in this paper.
\begin{table}[htp!]
\vspace{-4ex}
\begin{center}
\caption{A comparison with previous work $($$\beta$ is a smoothness parameter defined in \eqref{eq:smooth_g}$)$. }\label{tbl:comparison}
\begin{tabular}{ p{3.3cm} | p{3.3cm} | p{3.3cm} }\toprule
\multicolumn{1}{c|}{ADMM/iALM}  & \multicolumn{1}{c|}{Penalty / HOPS / ASGARD} & \multicolumn{1}{c}{This work} \\
\midrule
Constant or adaptive $\beta$. & Analytically drive $\beta$ to 0. & Analytically drive $\beta$ to 0. \\
\midrule
Update the dual center. & Do not move the dual center. &  Update the dual center. \\
\midrule
Theory is driven  by the convergence in the dual. & Do not analyze the convergence of the dual.  & Only analyze the stability of the primal-dual sequence. \\
\midrule
Inner problems are solved inexactly. & Inner problems are solved inexactly. & Only ensure stability for the number of inner iterations and smoothness parameter.\\
\bottomrule
\end{tabular}
\vspace{-6ex}
\end{center}
\end{table}

\beforesec
\section{Convergence analysis: The proof of Theorems~\ref{th:convergence_of_A2} and~\ref{th:convergence_A2b}}
\aftersec
We present the full proof of Theorems~\ref{th:convergence_of_A2} and~\ref{th:convergence_A2b} in this section.

\beforesubsec
\subsection{\textbf{The proof of Theorem~\ref{th:convergence_of_A2}: Convergence of Algorithm \ref{alg:A2} for \eqref{eq:ns_cvx_prob}}}
\aftersubsec
%
With the same argument as in \cite[Lemma 11]{TranDinh2015b}, we can prove the following estimate at the $k$-th iteration at the state $s$ of the outer loop, i.e., $K_s \leq k < K_{s+1}  := K_s +m_s$, of Algorithm~\ref{alg:A2}:
\begin{align}\label{eq:thm2_proof1}
&S_{\beta_s}(\bar{x}^{k+1};\dot{y}^s) + \tfrac{\tau_k^2 \Vert A\Vert^2 }{\beta_s} d_\Xc(x^\star, \hat{x}^{k+1}) \leq (1-\tau_k)S_{\beta_s}(\bar{x}^{k};\dot{y}^s) + \tfrac{\tau_k^2\Vert A\Vert^2}{\beta_s}d_\Xc(x^\star, \hat{x}^{k}) \vspace{1ex}\notag\\
&- \tau_k\beta_s b_\Yc (y^{\ast}_{\beta_s}(A\tilde{x}^k;\dot{y}^s), \dot{y}^s) - \tfrac{(1-\tau_k)\beta_s}{2}\Vert y^{\ast}_{\beta_s}(A\tilde{x}^k;\dot{y}^s) - y^{\ast}_{\beta_s}(A\bar{x}^k;\dot{y}^s)\Vert_\Yc^2,
\end{align}
where $S_{\beta}(\bar{x};\dot{y}) := P_{\beta}(\bar{x};\dot{y}) - P(x^{\star})$.
Note that this estimate remains true if we use APG instead of FISTA, and APG with Option 2.

Next, by strong convexity of $b_\Yc(\cdot, \dot{y} )$, the optimality condition of $g_{\beta}$-subproblem and convexity of $g^\ast(\cdot)$, we have
\begin{align}
g_{\beta}(A\bar{x};\dot{y}) &= \displaystyle\max_{y \in \R^n}\set{\iprods{A\bar{x}, y} - g^{\ast}(y) - \beta b_\Yc(y, \dot{y}) } \vspace{1ex} \notag\\
& \geq \langle A\bar{x}, y^\star \rangle - g^\ast (y^\star) - \beta b_\Yc(y^\star, \dot{y} ) + \beta b_\Yc(y^{\star}, y^\ast _{\beta} (A\bar{x};\dot{y}) ) \label{eq:gbeta_lb}.
\end{align}
Now, from the optimality condition of \eqref{eq:ns_cvx_prob}, we have $-A^\top y^\star \in \partial f(x^\star)$. 
Using this inclusion and convexity of $f$, we can derive
\begin{equation}\label{eq:f_conv} 
f(\bar{x}) \geq f(x^\star) + \langle -A^\top y^\star, \bar{x}-x^\star \rangle.
\end{equation}
Combining~\eqref{eq:gbeta_lb} and~\eqref{eq:f_conv}, we get
\begin{equation}\label{eq:thm2_proof2}
\begin{array}{ll}
S_\beta(\bar{x}; \dot{y}) & = P_{\beta}(\bar{x};\dot{y}) - P(x^{\star}) = f(\bar{x}) + g_{\beta}(A\bar{x}; \dot{y}) - \left(f(x^\star) + g(Ax^\star)\right) \vspace{1ex}\\
&\geq -\beta b_\Yc(y^\star, \dot{y}) + \beta b_\Yc(y^\star, y^\ast _{\beta} (A\bar{x};\dot{y})).
\end{array}
\end{equation}
From \eqref{eq:thm2_proof1}, if we ignore the two last terms, which are nonpositive, then for $K_s \leq k \leq K_s + m_s-1$ we obtain
\begin{align}\label{eq:thm3ns_proof1}
{\!\!\!}\tfrac{1}{\tau_k^2}S_{\beta_s}(\bar{x}^{k\!+\!1};\dot{y}^s) + \tfrac{\Vert A\Vert^2}{\beta_s} d_\Xc(x^\star, \hat{x}^{k\!+\!1}) \leq \tfrac{1-\tau_k}{\tau_k^2}S_{\beta_s}(\bar{x}^{k};\dot{y}^s)  + \tfrac{\Vert A\Vert^2}{\beta_s}d_\Xc(x^\star, \hat{x}^{k}).{\!\!\!}
\end{align}
Let us define $D^s_k := S_{\beta_s}(\bar{x}^k;\dot{y}^s) + \beta_s D_{\Yc} \geq P(\bar{x}^k) - P(x^{\star}) \geq 0$.
By adding $\frac{1}{\tau_k^2}\beta_s D_{\Yc}$ to both sides of \eqref{eq:thm3ns_proof1} and using the definition of $D^s_k$, we obtain
\begin{equation}\label{eq:thm3ns_proof1a}
\tfrac{1}{\tau_k^2}D_{k+1}^s + \tfrac{\Vert A\Vert^2}{\beta_s} d_\Xc(x^\star, \hat{x}^{k+1}) \leq \tfrac{(1-\tau_k)}{\tau_k^2}D_k^s + \tfrac{\Vert A\Vert^2}{\beta_s} d_\Xc(x^\star, \hat{x}^{k}) + \tfrac{\beta_s}{\tau_k}D_{\Yc}.
\end{equation}
Let us choose $\tau_k = \frac{2}{k - K_s +2}$. Then, it is clear that $\tau_{K_s} = 1$.
Moreover, $\frac{1-\tau_k}{\tau_k^2} = \frac{(k-K_s + 2)(k - K_s)}{4} \leq \frac{(k-K_s+1)^2}{4} = \frac{1}{\tau_{k-1}^2}$.
In this case, we can overestimate \eqref{eq:thm3ns_proof1a} as
\begin{equation}\label{eq:thm3ns_proof1b}
\tfrac{1}{\tau_k^2}D_{k+1}^s + \tfrac{\Vert A\Vert^2}{\beta_s}d_\Xc(x^\star, \hat{x}^{k+1}) \leq \tfrac{1}{\tau_{k-1}^2}D_k^s + \tfrac{\Vert A\Vert^2}{\beta_s}d_\Xc(x^\star, \hat{x}^{k}) + \tfrac{\beta_s}{\tau_k}D_{\Yc}.
\end{equation}
Taking a telescope from $k = K_s+1$ to $k = K_{s+1} - 1 = K_s + m_s - 1$ of \eqref{eq:thm3ns_proof1b} and reuse \eqref{eq:thm3ns_proof1a} for $k = K_s$, we obtain
\begin{align*} 
\begin{array}{ll}
D_{K_{s+1}}^s &+ \frac{\tau_{K_{s+1}-1}^2\Vert A\Vert^2}{\beta_s} d_\Xc(x^\star, \hat{x}^{K_{s+1}}) \leq \frac{\tau_{K_{s+1}-1}^2(1-\tau_{K_s})}{\tau_{K_s}^2} D_{K_s}^s  \vspace{0.75ex}\\ 
&+ \frac{\tau_{K_{s+1}-1}^2 \Vert A\Vert^2}{\beta_s} d_\Xc(x^\star, \hat{x}^{K_s}) + \beta_s\tau_{K_{s+1}-1}^2D_{\Yc}\sum_{j=K_s}^{K_s+m_s-1}\frac{1}{\tau_j} \vspace{0.75ex}\\
&\overset{\mathclap{(i)}}{\leq} \frac{\tau_{K_{s+1}-1}^2\Vert A\Vert^2}{\beta_s}d_\Xc(x^\star, \hat{x}^{K_s}) + \beta_s\tau_{K_{s+1}-1}^2D_{\Yc}\sum_{j=K_s}^{K_s+m_s-1}\frac{1}{\tau_j},
\end{array}
\end{align*}
where $(i)$ holds since $\tau_{K_s}=1$.
Since $\tau_k = \frac{2}{k - K_s +2}$, we have $\tau_{K_{s+1}-1} = \frac{2}{m_s+1}$ and $\sum_{j=K_s}^{K_s+m_s-1}\frac{1}{\tau_j} = \frac{m_s(m_s+3)}{4}$.
Using this relation, the last estimate leads to
\begin{align*} 
D_{K_{s+1}}^s + \tfrac{4\Vert A\Vert^2}{(m_s+1)^2\beta_s}d_\Xc(x^\star, \hat{x}^{K_{s+1}}) &\leq \tfrac{4\Vert A\Vert^2}{(m_s+1)^2\beta_s}d_\Xc(x^\star, \hat{x}^{K_s})  + \tfrac{\beta_sm_s(m_s+3)}{(m_s+1)^2}D_{\Yc}.
\end{align*}
Since $D_{K_{s+1}}^s = S_{\beta_s}(\bar{x}^{K_{s+1}};\dot{y}^s) + \beta_s D_{\Yc}$, the last estimate leads to
\begin{align}\label{eq:thm3ns_proof1d}
S_{\beta_s}(\bar{x}^{K_{s+1}};\dot{y}^{s}) + \tfrac{4\Vert A\Vert^2}{(m_s+1)^2\beta_s}d_\Xc(x^\star, \hat{x}^{K_{s+1}}) &\leq  \tfrac{4\Vert A\Vert^2}{(m_s+1)^2\beta_s}d_\Xc(x^\star, \hat{x}^{K_s}) \notag\\
& + \tfrac{\beta_s(m_s-1)}{(m_s+1)^2}D_{\Yc}.
\end{align}
Here, we note that $ \frac{\beta_sm_s(m_s+3)}{(m_s+1)^2}D_{\Yc}- \beta_sD_{\Yc} =  \frac{\beta_s(m_s-1)}{(m_s+1)^2}D_{\Yc}$.

Next, from \eqref{eq:thm2_proof2}, we have
\begin{equation}\label{eq:thm3_proof2}
S_{\beta_s}(\bar{x}^{K_{s+1}}; \dot{y}^s)  \geq \beta_s b_\Yc(y^\star, y^\ast _{\beta_{s}} (A\bar{x}^{K_{s+1}};\dot{y}^s)) - \beta_s b_\Yc(y^\star, \dot{y}^s ).
\end{equation}
Then, combining \eqref{eq:thm3ns_proof1d} and \eqref{eq:thm3_proof2} with the fact that $\dot{y}^{s+1} \leftarrow y^\ast _{\beta_s} (A\bar{x}^{K_{s+1}};\dot{y}^s)$, we can show that
\begin{align*}
\tfrac{4\Vert A \Vert^2}{(m_s+1)^2} d_\Xc(x^\star, \hat{x}^{K_{s+1}}) + \beta_s^2 b_\Yc( y^\star, \dot{y}^{s+1}) &\leq \tfrac{4\Vert A\Vert^2}{(m_s+1)^2}d_\Xc(x^\star, \hat{x}^{K_s}) \\
&+ \beta_s^2 b_\Yc (y^\star, \dot{y}^s) +  \tfrac{\beta_s^2(m_s-1)}{(m_s+1)^2}D_{\Yc}.
\end{align*}
Using the update rule \eqref{eq:update_param1a} $m_{s+1} \leftarrow \lfloor \omega(m_s+1) + 1\rfloor - 1$ , we have 
\begin{equation}\label{eq:m_s_ineq}
\omega(m_s + 1)  \leq m_{s+1} + 1 \leq \omega(m_s + 1) + 1.
\end{equation}
Define $q_s := \frac{\beta_s^2(m_s-1)}{(m_s+1)^2}D_{\Yc}$, then using $\beta_{s+1} \leftarrow \frac{\beta_s}{\omega}$ from \eqref{eq:update_param1a} and $\omega(m_s + 1)  \leq m_{s+1} + 1$, we obtain
\begin{align*}
\tfrac{4\Vert A \Vert^2}{(m_{s+1}+1)^2}d_\Xc(x^\star, \hat{x}^{K_{s+1}}) + \beta_{s+1}^2 b_\Yc(y^\star, \dot{y}^{s+1} ) &\leq \tfrac{1}{\omega^2}\big[ \tfrac{4\Vert A\Vert^2}{(m_s+1)^2} d_\Xc(x^\star, \hat{x}^{K_s}) \\
&+ \beta_s^2 b_\Yc( y^\star, \dot{y}^s )\big] + \tfrac{q_s}{\omega^2}.
\end{align*}
Telescoping this inequality from $s \leftarrow 0$ to $s \leftarrow s-1$, we finally obtain 
\begin{align}\label{eq:thm3_proof4a}
\tfrac{4\Vert A\Vert^2}{(m_s+1)^2} d_\Xc(x^\star, \hat{x}^{K_s}) + \beta_s^2 b_\Yc( y^\star, \dot{y}^s ) &\leq \tfrac{1}{\omega^{2s}} \big[ \tfrac{4\Vert A\Vert^2}{(m_0+1)^2}d_\Xc(x^\star, \hat{x}^{0}) \notag \\
&+ \beta_0^2 b_\Yc( y^\star, \dot{y}^0) \big] + Q_s,
\end{align}
where $Q_s := \frac{q_{s-1}}{\omega^2} + \frac{q_{s-2}}{\omega^{4}} + \cdots + \frac{q_0}{\omega^{2s}}$.
If we ignore the second term on the left-hand side, which is nonnegative, and use the equality $\hat x^0 = \bar x^0$, we obtain
\begin{equation}\label{eq:thm3_proof4}
{\!\!\!}\tfrac{4\Vert A\Vert^2}{(m_s + 1)^2}d_\Xc(x^\star, \hat{x}^{K_s}) \leq \tfrac{1}{\omega^{2s}} \big[ \tfrac{4\Vert A\Vert^2}{(m_0+1)^2}d_\Xc(x^\star, \bar{x}^{0}) + \beta_0^2 b_\Yc(y^\star, \dot{y}^0 ) \big]  + Q_s.{\!\!\!}
\end{equation}
Lower bounding the second term on the left-hand side in \eqref{eq:thm3ns_proof1d} by $0$, and combining the result with \eqref{eq:thm3_proof4} we obtain
\begin{equation}\label{eq:thm3_proof5} 
\begin{array}{ll}
S_{\beta_s}(\bar{x}^{K_{s+1}};\dot{y}^s)  &\leq \frac{4\Vert A\Vert^2}{(m_s+1)^2\beta_s}d_\Xc(x^\star, \hat{x}^{K_s}) + \frac{\beta_s(m_s-1)}{(m_s+1)^2}D_{\Yc} \vspace{-1ex}\\
&\leq  \frac{1}{\beta_s\omega^{2s}} \left[ \frac{4\Vert A\Vert^2}{(m_0+1)^2}d_\Xc(x^\star, \bar{x}^{0}) + \beta_0^2 b_\Yc(y^\star, \dot{y}^0 ) \right] + \overbrace{\tfrac{Q_s}{\beta_s} +  \tfrac{q_s}{\beta_s}}^{\hat{Q}_s}.
\end{array}
\end{equation}
Note from \eqref{eq:update_param1a} that $\beta_s = \frac{\beta_0}{\omega^s}$ implies $\frac{1}{\beta_s\omega^{2s}} = \frac{1}{\beta_0\omega^s}$.
By induction of~\eqref{eq:m_s_ineq}, yields
\begin{equation}\label{eq:condition_omega}
m_0\omega^s < (m_0 + 1)\omega^s - 1 \leq  m_{s} \leq \big(m_0 + \tfrac{\omega}{\omega - 1}\big)\omega^s - \tfrac{\omega}{\omega - 1} < \kappa_0\omega^s,
\end{equation}
where $\kappa_0 := m_0 + \frac{\omega}{\omega - 1} \!>\! 0$ for $\omega > 1$ and $m_0 \geq 1$.
Using these bounds, one has
\begin{equation*}
q_s = \tfrac{\beta_s^2(m_s - 1)}{(m_s + 1)^2}D_{\Yc} \leq \tfrac{\beta_0^2}{\omega^{2s}(m_s + 1)}D_{\Yc}  \leq  \tfrac{\beta_0^2}{m_0\omega^{3s}}D_{\Yc}.
\end{equation*}
Substituting this inequality into $\hat{Q}_s$, we can bound
\begin{align}\label{eq:thm3_proof6}
\hat{Q}_s &:= \tfrac{1}{\beta_s}\left( q_s + \tfrac{q_{s-1}}{\omega^2} + \cdots + \tfrac{q_{0}}{\omega^{2s}}\right)  \leq \tfrac{\beta_0\omega^sD_{\Yc}}{m_0}\left(\tfrac{1}{\omega^{3s}} + \tfrac{1}{\omega^2\omega^{3(s-1)}} + \cdots + \tfrac{1}{\omega^{2s}} \right) \nonumber\\
&\leq \tfrac{\beta_0D_{\Yc}}{m_0\omega^s}\left(\tfrac{1}{\omega^s} + \tfrac{1}{\omega^{s-1}} + \cdots + \tfrac{1}{\omega}  + 1\right)  \leq \tfrac{\beta_0\omega D_{\Yc}}{(\omega-1)m_0\omega^s}.
\end{align}
Using $m_s \leq \kappa_0\omega^s$ in \eqref{eq:condition_omega} to estimate the total number of iterations $K_{s+1}$ as
\begin{equation*}
K_{s+1} = \sum_{i=0}^{s}m_i \leq \kappa_0\sum_{i=0}^{s}\omega^i = \kappa_0\left(\tfrac{\omega^{s+1} - 1}{\omega - 1}\right).
\end{equation*}
This condition leads to $\omega^s \geq \frac{(\omega-1)K_{s+1} + \kappa_0}{\omega \kappa_0}$. 
Using this estimate, $\beta_s=\tfrac{\beta_0}{\omega^s}$ and \eqref{eq:thm3_proof6} into \eqref{eq:thm3_proof5}, we obtain
\begin{align*} 
S_{\beta_s}(\bar{x}^{K_{s+1}};\dot{y}^s)  \leq \tfrac{1}{\beta_0\omega^{s}} \big[ \tfrac{4\Vert A\Vert^2}{(m_0+1)^2} d_\Xc(x^\star, \bar{x}^{0}) &+ \beta_0^2 b_\Yc( y^\star, \dot{y}^0 ) + \tfrac{\beta_0^2\omega D_{\Yc}}{(\omega-1)m_0}\big].
\end{align*}
Finally, by the fact that $P(\bar{x}^{K_{s+1}}) - P^{\star} \leq S_{\beta_s}(\bar{x}^{K_{s+1}};\dot{y}^s) + \beta_sD_{\Yc}$, and $\beta_s  = \frac{\beta_0}{\omega^s} \leq \frac{\omega \beta_0\kappa_0}{(\omega-1)K_{s+1} + \kappa_0}$,  the last estimate implies \eqref{eq:A2_key_est2a}.
\Eproof

\beforesubsec
\subsection{\textbf{The proof of Theorem~\ref{th:convergence_A2b}: Convergence of Algorithm \ref{alg:A2} for \eqref{eq:constr_cvx}}}
\aftersubsec
%
%
By Lemma~\ref{le:opt_cond_constr_cvx2}, we have $\beta b_{\Yc} (y^{\star}, \dot{y}) + S_{\beta}(\bar{x};\dot{y}) \geq 0$.
Let us define $\hat{D}^s_k := S_{\beta_s}(\bar{x}^k;\dot{y}^s) + \beta_s b_{\Yc}(y^{\star}, \dot{y}^s )$.
Then, using the same proof as \eqref{eq:thm3ns_proof1a} in Theorem~\ref{th:convergence_of_A2} by replacing $D_{\Yc}$ with $b_{\Yc}(y^\star, \dot{y}^s)$.
\begin{equation}\label{eq:thm4ns_est1}
{\!\!\!}\tfrac{1}{\tau_k^2}\hat{D}_{k+1}^s + \tfrac{\norms{A}^2}{\beta_s} d_{\Xc}(x^\star, \hat{x}^{k+1}) \leq \tfrac{(1-\tau_k)}{\tau_k^2}\hat{D}_k^s + \tfrac{\norms{A}^2}{\beta_s}d_{\Xc}(x^\star, \hat{x}^{k})  + \tfrac{\beta_s}{\tau_k}b_{\Yc}(y^\star, \dot{y}^s).{\!\!\!}
\end{equation}
From this estimate, with the same proof as \eqref{eq:thm3ns_proof1d}, we obtain
\begin{equation}\label{eq:thm4ns_proof10} 
\begin{array}{ll}
S_{\beta_s}(\bar{x}^{K_{s+1}};\dot{y}^{s}) + \frac{4\norms{A}^2}{(m_s+1)^2\beta_s} d_{\Xc}(x^\star, \hat{x}^{K_{s+1}}) &\leq  \frac{4\norms{A}^2}{(m_s+1)^2\beta_s}d_{\Xc}(x^\star, \hat{x}^{K_s}) \vspace{1ex}\\
& + \frac{\beta_s(m_s-1)}{(m_s+1)^2} b_{\Yc}(y^{\star}, \dot{y}^s).
\end{array}
\end{equation}
Combining this, \eqref{eq:thm3_proof2}, and  the fact that $\dot{y}^{s+1} \leftarrow y^\ast _{\beta_s} (A\bar{x}^{K_{s+1}}; \dot{y}^s)$, we obtain
\begin{equation}\label{eq:thm4ns_est10}
\begin{array}{ll}
\frac{4\norms{A}^2}{(m_s+1)^2}d_{\Xc}(x^\star, \hat{x}^{K_{s+1}}) + \beta_s^2 b_{\Yc}( y^\star, \dot{y}^{s+1} )&\leq \frac{4\norms{A}^2}{(m_s+1)^2}d_{\Xc}(x^{\star}, \hat{x}^{K_s}) \vspace{1ex}\\
 & + \frac{\beta_s^2m_s(m_s+3)}{(m_s+1)^2} b_{\Yc}( y^\star, \dot{y}^s).
 \end{array}
\end{equation}
Let us choose $m_{s+1} := \lfloor\omega(m_s+1) + 1\rfloor - 1$ and $\beta_{s+1} := \frac{\beta_s(m_{s+1}+1)}{\omega\sqrt{m_{s+1}(m_{s+1}+3)}}$ as in \eqref{eq:update_param2a}.
Then, similar to the proof of \eqref{eq:condition_omega}, we have 
\begin{equation}\label{eq:upper_bound_of_params}
m_0\omega^s \leq m_s \leq \kappa_0\omega^s, ~~~\text{and}~~~\beta_{s+1} \leq \tfrac{\beta_s}{\omega} \leq \tfrac{\beta_0}{\omega^{s+1}},
\end{equation}
where $\kappa_0 := m_0 + \frac{\omega}{\omega - 1} > 0$.

\noindent Next, we need to lower bound $\beta_s$.
We can show that, for $m_s \geq 1$, we have
\begin{equation*}
\tfrac{m_{s+1}+1}{\sqrt{m_{s+1}(m_{s+1}+3)}} \geq 1 - \tfrac{1}{m_{s+1}} \geq 0.
\end{equation*}
In this case, we can estimate $\beta_{s+1} = \tfrac{\beta_s(m_{s+1}+1)}{\omega\sqrt{m_{s+1}(m_{s+1}+3)}} \geq \tfrac{\beta_s}{\omega}\left(1 - \tfrac{1}{m_{s+1}}\right) = \tfrac{\beta_s}{\omega} - \tfrac{\beta_s}{m_{s+1}\omega}$.
Substituting \eqref{eq:upper_bound_of_params} on $m_{s+1}$ and $\beta_s$ into this inequality, we obtain
\begin{equation*}
\beta_{s+1} \geq \tfrac{\beta_s}{\omega} - \tfrac{c_0}{\omega^{2s+1}},~~~\text{where}~~~c_0 := \tfrac{\beta_0}{\omega m_0}.
\end{equation*}
This condition leads to $\omega \beta_{s+1} + \frac{c_0}{\omega^{2s}} \geq \beta_s$.
By induction, we can show that $\omega^s\beta_s + c_0\sum_{j=0}^{s-1}\frac{1}{\omega^j} \geq \beta_0$, which leads to
\begin{equation}\label{eq:lower_bound_of_params}
\beta_s  \geq \tfrac{1}{\omega^s}\left(\beta_0 - \tfrac{c_0\omega(\omega^s - 1)}{(\omega - 1)\omega^s}\right) \geq \beta_0\left(1 - \tfrac{1}{m_0(\omega-1)}\right) \tfrac{1}{\omega^s}.
\end{equation}
Here, we use the fact that 
\[\rho_0:=\beta_0 - \frac{c_0\omega(\omega^s - 1)}{(\omega - 1)\omega^s} \geq \beta_0 - \frac{c_0\omega}{\omega-1} = \beta_0\left(1 - \frac{1}{m_0(\omega-1)}\right) > 0\]
since $m_0 > \frac{1}{\omega-1}$.
This condition gives us a lower bound on $\beta_s$.

Now, using \eqref{eq:update_param2a}, we have $\frac{\omega^2\beta^2_{s+1}m_{s+1}(m_{s+1}+3)}{(m_{s+1}+1)^2} = \beta_s^2$ and $\frac{\omega^2}{(m_{s+1}+1)^2} \leq \frac{1}{(m_s+1)^2}$ as in~\eqref{eq:m_s_ineq}. 
Plugging these estimates into \eqref{eq:thm4ns_est10}, we obtain
\begin{align}\label{eq:thm4ns_proof10b} 
\tfrac{4\norms{A}^2}{(m_{s+1}+1)^2}&d_{\Xc}(x^\star, \hat{x}^{K_{s+1}}) + \tfrac{\beta_{s+1}^2m_{s+1}(m_{s+1}+3)}{(m_{s+1}+1)^2} b_{\Yc}( y^\star, \dot{y}^{s+1} ) \leq \notag \\
&\tfrac{1}{\omega^2}\Big[\tfrac{4\norms{A}^2}{(m_s+1)^2}d_{\Xc}(x^\star, \hat{x}^{K_s}) + \tfrac{\beta_s^2m_s(m_s +3)}{(m_s + 1)^2}b_{\Yc}(y^{\star},  \dot{y}^s)\Big].
\end{align}
By induction and using that $\hat{x}^0 = \bar{x}^0$, we obtain
\begin{equation}\label{eq:thm4ns_proof11}
\begin{array}{ll}
 \frac{4\norms{A}^2}{(m_{s}+1)^2}d_{\Xc}(x^\star, \hat{x}^{K_{s}}) + \frac{\beta_{s}^2m_{s}(m_{s}+3)}{(m_{s} + 1)^2} b_{\Yc}( y^\star, \dot{y}^{s} ) &\leq \frac{1}{\omega^{2s}}\Big[\frac{4\norms{A}^2}{(m_0+1)^2}d_{\Xc}(x^\star, \bar{x}^{0}) \vspace{1ex}\\
 &  + \frac{\beta_0^2m_0(m_0+3)}{(m_0+1)^2} b_{\Yc}( y^\star, \dot{y}^0 )\Big].
 \end{array}
\end{equation}
Since $m_s(m_s+3) \geq m_s - 1$, combining \eqref{eq:thm4ns_proof11} and \eqref{eq:thm4ns_proof10}, we obtain
\begin{align}\label{eq:thm4ns_proof12}
S_{\beta_s}(\bar{x}^{K_{s+1}};\dot{y}^{s}) &\leq \tfrac{1}{\beta_s\omega^{2s}}\left[\tfrac{4\norms{A}^2}{(m_0+1)^2}d_{\Xc}(x^\star, \bar{x}^{0}) + \tfrac{\beta_0^2m_0(m_0+3)}{(m_0+1)^2} b_{\Yc}( y^\star, \dot{y}^0 )\right] \notag \\
&\leq \tfrac{R_0^2}{\beta_s\omega^{2s}},
\end{align}
where $R_0 := \left[\frac{4\norms{A}^2}{(m_0+1)^2} d_{\Xc}(x^\star, \bar{x}^{0})  + \frac{\beta_0^2m_0(m_0+3)}{(m_0+1)^2} b_{\Yc}( y^\star, \dot{y}^0 )\right]^{1/2}$.

Using \eqref{eq:upper_bound_of_params} and \eqref{eq:lower_bound_of_params} of $\beta_s$ and $m_s$ into \eqref{eq:thm4ns_proof12}, we obtain
\begin{equation}\label{eq:thm4ns_proof12/2}
S_{\beta_s}(\bar{x}^{K_{s+1}};\dot{y}^{s})  \leq \tfrac{R_0^2}{\rho_0\omega^{s}} \leq \tfrac{\omega\kappa_0 R_0^2}{\rho_0\left[(\omega-1)K_{s+1} + \kappa_0\right]}.
\end{equation}
Here, we use the same argument as in Theorem~\ref{th:convergence_of_A2} to bound $\omega^s$ via the number of iterations $K_{s+1}$ as $\omega^s \geq \frac{(\omega-1)K_{s+1} + \kappa_0}{\omega \kappa_0}$, and $\rho_0 := \beta_0\big(1 - \frac{1}{m_0(\omega-1)}\big) > 0$.

Our next step is using \eqref{eq:thm4ns_proof11} to bound $\Vert \dot{y}^s - y^{\star}\Vert_{\Yc}$.
Clearly, $\frac{\beta_{s}^2m_{s}(m_{s}+3)}{(m_{s} + 1)^2} = \frac{\beta_{s-1}^2}{\omega^2} \geq \frac{\rho_0^2}{\omega^{2s}}$ by \eqref{eq:lower_bound_of_params}.
Using \eqref{eq:thm4ns_proof11}, and strong convexity of $b_{\Yc}$ with respect to the given norm, we can show that
\begin{equation}\label{eq:thm4ns_proof13}
\tfrac{1}{2} \Vert \dot{y}^s - y^{\star} \Vert_{\Yc} ^2 \leq b_{\Yc}(y^\star, \dot{y}^s) \leq \tfrac{R_0^2}{\rho_0^2}.
\end{equation}

\beforesubsubsec
\subsubsection{The first estimate of \eqref{eq:key_estimate3b}.}
\aftersubsubsec
First, using Lemma~\ref{le:opt_cond_constr_cvx2}, and by defining $\bar{\beta}_s:=\beta_s L_{b_{\Yc}}$ we write
\begin{align*}
f(\bar{x}^{K_{s+1}}) -& f^\star \geq \bar{\beta}_s\langle \dot{y}^s, y^\star \rangle - \Vert y^\star \Vert \kdist{\Kc}{A\bar{x}^{K_{s+1}} -b + \bar{\beta}_s \dot{y}^s} \\
&\overset{\mathclap{{\eqref{eq:constr_san}}}}{\geq} \bar{\beta}_s\langle \dot{y}^s, y^\star \rangle - \Vert y^\star \Vert \kdist{\Kc}{A\bar{x}^{K_{s+1}}-b} - \Vert y^\star\Vert \bar{\beta}_s (\Vert \dot{y}^s-y^\star \Vert + \Vert y^\star \Vert) \\
&\geq -2 \bar{\beta}_s \Vert y^\star \Vert \Vert \dot{y}^s - y^\star \Vert - \Vert y^\star \Vert \kdist{\Kc}{A\bar{x}^{K_{s+1}}-b} 
\end{align*}
By using the bound of $\Vert \dot{y}^s - y^\star \Vert$ from~\eqref{eq:thm4ns_proof13} along with $\beta_s \leq \frac{\beta_0}{\omega_s} \leq \frac{\omega\beta_0\kappa_0}{(\omega-1)K_{s+1} + \kappa_0}$, we conclude that
\begin{equation}
f(\bar{x}^{K_{s+1}} )- f^\star \geq - \Vert y^\star \Vert \kdist{\Kc}{A\bar{x}^{K_{s+1}}-b} - \tfrac{2\sqrt{2}\omega\beta_0L_{b_{\Yc}}\kappa_0\Vert y^\star \Vert R_0}{\rho_0\left[(\omega-1)K_{s+1} + \kappa_0\right]}.
\end{equation}

\beforesubsubsec
\subsubsection{The second estimate of~\eqref{eq:key_estimate3b}.}
\aftersubsubsec
Using Lemma~\ref{le:opt_cond_constr_cvx2}, we have
\begin{align}\label{eq:constr_2}
f(\bar{x}^{K_{s+1}}) - f^\star &\leq S_{\beta_s}(\bar{x}^{K_{s+1}};\dot{y}^s) + \tfrac{\bar{\beta}_s}{2} \Vert \dot{y}^s\Vert ^2\nonumber \\
&\leq S_{\beta_s}(\bar{x}^{K_{s+1}};\dot{y}^s) + \tfrac{\bar{\beta}_s}{2} \left( \Vert \dot{y}^s - y^\star \Vert ^2 + \Vert y^\star \Vert ^2 \right).
\end{align}
Combining this bound,~\eqref{eq:thm4ns_proof12/2} and~\eqref{eq:thm4ns_proof13}  into~\eqref{eq:constr_2} gives the second bound of~\eqref{eq:key_estimate3b}.

\beforesubsubsec
\subsubsection{Third estimate of~\eqref{eq:key_estimate3b}.}
\aftersubsubsec
Finally, we note that, by using \eqref{eq:lower_bound_of_params} and \eqref{eq:thm4ns_proof12}, we can bound $\frac{2S_{\beta_s}(\bar{x}^{K_{s+1}};\dot{y}^s)}{\beta_s} \leq \frac{2R_0^2}{\rho_0^2}$.
Using this upper bound,~\eqref{eq:constr_san} and \eqref{eq:thm4ns_proof13} into the third estimate of \eqref{eq:approx_opt_cond}, we obtain the third bound  of \eqref{eq:key_estimate3b}.
\Eproof

\beforesec
\section{Appendix: The proof of technical results}
\aftersec
%
This appendix provides the missing proof of the results in the main text.
\beforesubsec
\subsection{\textbf{The proof of Example~\ref{ex:cone_constr_l2}}.}
\aftersubsec
In this example, we have $b_\Yc(y,\dot{y})=\frac{1}{2}\Vert y- \dot{y}\Vert^2$.
First, from the definition \eqref{eq:gAx3_beta} of $g_{\beta}(Ax;\dot{y})$, by using  the definition of $s_{\Kc}$, we write
\begin{align*} 
g_{\beta}(Ax;\dot{y}) &= \displaystyle\min_{u\in\Kc}\displaystyle\max_{y\in\R^n}\set{ \iprods{Ax - b - u, y} - \beta b_\Yc(y, \dot{y})} \\
&= \displaystyle\min_{u\in\Kc}\displaystyle\max_{y\in\R^n}\set{ \iprods{Ax - b - u, y} - \tfrac{\beta}{2}\Vert y - \dot{y}\Vert^2}.
\end{align*}
The optimality condition of the $\max$ problem on the right hand side of the previous inequality is $Ax - b - u - \beta (y - \dot{y}) = 0$, which implies $y = \dot{y} + \frac{1}{\beta }(Ax - b - u)$.
In this case, $\iprods{Ax - b - u, y} - \tfrac{\beta }{2}\Vert y - \dot{y}\Vert^2 = \frac{1}{2\beta }\Vert Ax - b - u\Vert^2 + \iprods{\dot{y}, Ax - b - u} = \frac{1}{2\beta }\Vert Ax - b - u + \beta\dot{y}\Vert^2 - \frac{\beta}{2}\Vert\dot{y}\Vert^2$.
Hence, we obtain 
\begin{align*}
g_{\beta}(Ax;\dot{y}) &= \displaystyle\min_{u\in\Kc}\set{\tfrac{1}{2\beta}\Vert u - (Ax - b + \beta\dot{y})\Vert^2} - \tfrac{\beta}{2}\Vert\dot{y}\Vert^2 \\
&= \tfrac{1}{2\beta}\kdist{\Kc}{Ax - b + \beta\dot{y}}^2 - \tfrac{\beta}{2}\Vert\dot{y}\Vert^2,
\end{align*}
which is \eqref{eq:gAx3_beta_b}.
In addition, this implies $u = \kproj{\Kc}{Ax - b + \beta\dot{y}}$.
Hence, we obtain $y^{\ast}_{\beta}(Ax;\dot{y}) = \dot{y} + \frac{1}{\beta}(Ax - b - u) = \dot{y} + \frac{1}{\beta}\left(Ax - b - \kproj{\Kc}{Ax - b + \beta\dot{y}}\right)$, which is exactly \eqref{eq:ystar_beta}.

If $\Kc$ is a cone, then using Moreau's decomposition \cite[Theorem 6.30]{Bauschke2011}, we can show that 
\begin{equation*}
Ax - b + \beta\dot{y} - \kproj{\Kc}{Ax - b + \beta\dot{y}} = \kproj{\Kc^{\circ}}{Ax - b + \beta\dot{y}}, 
\end{equation*}
where $K^{\circ}$ is the polar set of $\Kc$.
Since $\Kc$ is a cone, $\Kc^{\circ} = -\Kc^{\ast}$, where $\Kc^{\ast}$ is the dual cone of $\Kc$.
Hence, we have $y^{\ast}_{\beta}(Ax;\dot{y}) = \kproj{-\Kc^{\ast}}{\dot{y} + \frac{1}{\beta}(Ax - b)}$.
\Eproof

\beforesubsec
\subsection{\textbf{The proof of Lemma~\ref{le:opt_cond_constr_cvx2}: Optimality bounds for \eqref{eq:constr_cvx}}.}\label{apdx:le:opt_cond_constr_cvx2}
\aftersubsec
Using the property of distance function, we can derive
\begin{equation}\label{eq:thm4ns_estmate0}
\begin{array}{ll}
\kdist{\Kc}{A\bar{x} - b} &\leq \kdist{\Kc}{A\bar{x} - b + \beta\dot{y}} + \kdist{\Kc}{\beta\dot{y}} \vspace{1ex}\\
&= \kdist{\Kc}{A\bar{x} - b + \beta\dot{y}}  + \min_{u\in\Kc}\Vert u - \beta\dot{y}\Vert\\
& \overset{\mathclap{(i)}}{\leq} \kdist{\Kc}{A\bar{x} - b + \beta\dot{y}}  + \beta\Vert\dot{y}\Vert \vspace{1ex}\\
& \leq \kdist{\Kc}{A\bar{x} - b + \beta\dot{y}}  + \beta\left(\Vert\dot{y} - y^{\star}\Vert + \Vert y^{\star}\Vert\right),
\end{array}
\end{equation}
where $(i)$ holds since $\boldsymbol{0}^n\in\Kc$. Similarly, we can start from $\kdist{\Kc}{A\bar{x}-b + \beta \dot{y}} \leq \kdist{\Kc}{A\bar{x}-b}+\kdist{\Kc}{\beta \dot{y}}$, to get the similar bound
\begin{equation}\label{eq:dist_bound_opposite}
\kdist{\Kc}{A\bar{x}-b+\beta\dot{y}} \leq \kdist{\Kc}{A\bar{x}-b}+\beta\left(\Vert \dot{y} - y^\star \Vert + \Vert y^\star \Vert \right).
\end{equation}
Let $y^{\star}$ be an arbitrary optimal solution of the dual problem \eqref{eq:constr_cvx_dual}.
Applying the strong duality condition for \eqref{eq:constr_cvx} and \eqref{eq:constr_cvx_dual}, we have
\begin{equation*}
\begin{array}{ll}
-f(x^{\star}) &= -f^{\star} = D^{\star} = f^{\ast}(-A^Ty^{\star}) + \iprods{b, y^{\star}} + s_{\Kc}(y^{\star}) \\
&\overset{\mathclap{(i)}}{\geq} \iprods{b-A\bar{x}, y^{\star}} - f(\bar{x}) +  s_{\Kc}(y^{\star}) \overset{\mathclap{(ii)}}{=} \displaystyle\max_{u\in\Kc}\set{\iprods{y^{\star}, b - A\bar{x} + u}} - f(\bar{x})  \\
&= \displaystyle\max_{u\in\Kc}\set{ - \iprods{y^{\star}, A\bar{x} - b - u}} - f(\bar{x}),
\end{array}
\end{equation*}
where $(i)$ follows by the definition of conjugate function $f^\ast(\cdot)$, and $(ii)$ follows by the definition of support function $s_\Kc(\cdot)$.
By rearranging, we get the following relation:
\begin{equation*}
f(\bar{x}) - f(x^\star) \geq \displaystyle\max_{u\in\Kc}\set{ - \iprods{y^{\star}, A\bar{x} - b - u}}.
\end{equation*}
Now, since $\iprods{y^{\star}, A\bar{x} - b - u + \beta\dot{y}} \leq \Vert y^{\star}\Vert \Vert A\bar{x} - b - u + \beta\dot{y}\Vert$, we have
\begin{equation}\label{eq:lmA1_est1a}
\begin{array}{ll}
\displaystyle\max_{u\in\Kc}\set{ - \iprods{y^{\star}, A\bar{x} - b - u}} - \beta\iprods{y^{\star}, \dot{y}} &\geq \displaystyle\max_{u \in\Kc}\set{ - \Vert y^{\star}\Vert \Vert A\bar{x} - b - u + \beta\dot{y}\Vert} \vspace{1ex}\\
& = -\Vert y^{\star}\Vert \displaystyle\min_{u\in\Kc}\Vert A\bar{x} - b + \beta\dot{y} - u\Vert \vspace{1ex}\\
& = -\Vert y^{\star}\Vert\kdist{\Kc}{A\bar{x}-b + \beta\dot{y}}.
\end{array}
\end{equation}
Combining the above two inequalities, we obtain
\begin{equation}\label{eq:lm_f_lower_bound}
f(\bar{x}) - f(x^{\star}) \geq \beta\iprods{y^{\star},\dot{y}} -\Vert y^{\star}\Vert\kdist{\Kc}{A\bar{x} - b + \beta\dot{y}},
\end{equation}
which is the first estimate of \eqref{eq:approx_opt_cond}.

Since we assume  $b_\Yc(y,\dot{y})$ with Lipschitz gradients, we can write
\begin{equation} \label{eq:gb_bound}
\begin{array}{ll}
g_{\beta}(Ax;\dot{y}) &= \displaystyle\min_{u\in\Kc}\displaystyle\max_{y\in\R^n}\{ \iprods{Ax - b - u, y} - \beta b_\Yc(y, \dot{y})\} \vspace{0.0ex}\\
&\overset{\tiny\mathclap{(ii)}}{\geq} \displaystyle\min_{u\in\Kc}\displaystyle\max_{y\in\R^n}\set{ \iprods{Ax - b - u, y} - \tfrac{\beta L_{b_\Yc}}{2}\Vert y - \dot{y}\Vert^2} \vspace{0.0ex}\\
& \overset{\mathclap{(ii)}}{=} \tfrac{1}{2 \beta_b}\kdist{\Kc}{Ax - b +  \beta_b \dot{y}}^2 - \tfrac{\beta_b}{2}\Vert\dot{y}\Vert^2,
\end{array}
\end{equation}
where $(i)$ holds because of the Lipschitz gradient assumption on $b_\Yc$, and $(ii)$ follows from Example~\ref{ex:cone_constr_l2} and defining $\beta_b := \beta L_{b_\Yc}$.

Invoking~\eqref{eq:lmA1_est1a} with $\beta=\beta_b$ and combining it with~\eqref{eq:lm_f_lower_bound} and $S_{\beta}(\cdot;\dot{y}) = f(\bar{x}) + g_{\beta}(\bar{x};\dot{y}) - f(x^{\star})$ gives
\begin{equation}\label{eq:lmA1_est1}
\beta_b\iprods{y^{\star},\dot{y}} \!-\! \Vert y^{\star}\Vert\kdist{\Kc}{A\bar{x} \!-\! b \!+\! \beta_b \dot{y}} \leq f(\bar{x}) \!-\! f(x^{\star})  = S_{\beta}(\bar{x};\dot{y}) - g_{\beta}(A\bar{x};\dot{y}).
\end{equation}
Using \eqref{eq:gb_bound} into this inequalty, we have
\begin{align*}
\beta_b \iprods{y^{\star},\dot{y}} -\Vert y^{\star}\Vert\kdist{\Kc}{A\bar{x} - b + \beta_b \dot{y}} &\leq S_{\beta}(\bar{x};\dot{y}) - \tfrac{1}{2\beta_b}\kdist{\Kc}{Ax - b + \beta_b \dot{y}}^2 \\
&+ \tfrac{\beta_b }{2}\Vert\dot{y}\Vert^2.
\end{align*}
Let $t := \kdist{\Kc}{A\bar{x} - b + \beta_b \dot{y}}$. Then, this inequality becomes 
\begin{equation}\label{eq:t_eq}
\tfrac{1}{2\beta_b}t^2 - \Vert y^{\star}\Vert t - \left(S_{\beta}(\bar{x};\dot{y}) + \tfrac{\beta_b}{2}\Vert\dot{y}\Vert^2 - \beta_b \iprods{y^{\star},\dot{y}}\right) \leq 0.
\end{equation}
By using the strong convexity of $b_\Yc(\cdot, \dot{y})$ with respect to the corresponding norm, one can plug in the optimality condition of the maximization problem~\eqref{eq:smooth_g} to derive
\begin{equation*}
2b_\Yc(y^{\star}, \dot{y}) + \tfrac{2}{\beta}S_{\beta}(\bar{x};\dot{y}) \geq \Vert y^{\ast}_{\beta}(A\bar{x};\dot{y}) - \dot{y}\Vert^2 \geq 0.
\end{equation*}
By using this inequality and strong convexity of $b_\Yc(\cdot, \dot{y})$, we conclude that $\Vert y^\star - \dot{y} \Vert ^2 + \frac{2}{\beta}S_\beta(\bar{x};\dot{y}) \geq 0$, therefore the inequation~\eqref{eq:t_eq} has solution. Consequently, we can write that 
\begin{equation*}
t^{\ast} = \kdist{\Kc}{A\bar{x} - b + \beta_b \dot{y}} \leq \beta_b \Big[ \Vert y^{\star}\Vert + \left(\Vert y^{\star} - \dot{y}\Vert^2 + \tfrac{2}{\beta_b}S_{\beta}(\bar{x};\dot{y})\right)^{1/2}  \Big],
\end{equation*} 
which is the third estimate of \eqref{eq:approx_opt_cond}.

\noindent Finally, plugging~\eqref{eq:gb_bound} into~\eqref{eq:lmA1_est1}, we obtain the second estimate of~\eqref{eq:approx_opt_cond}.
\Eproof

\bibliographystyle{plain}

\end{document}